\documentclass[11pt]{article}
\usepackage[latin1]{inputenc}
\usepackage[a4paper,width=18cm,height=25cm]{geometry}
\usepackage{amsmath,amsfonts,amssymb,indentfirst,bbm}
\usepackage{graphicx,multicol}
\newcommand{\R}{\mathbb{R}}
\renewcommand{\div}{\operatorname{div}}
\newcommand{\norme}[2][]{\left\|#2\right\|_{#1}}
\usepackage{theorem}
{\theorembodyfont{\slshape}
  \newtheorem{thm}{Theorem}
  \newtheorem{cor}[thm]{Corollary}
  \newtheorem{prop}[thm]{Proposition}
  \newtheorem{dfn}{Definition}
}

\newcommand{\Proof}[1][]{{\par\smallskip\noindent{\bf Proof#1.\enspace }}}
\def\endProof{\par{\medskip\noindent}}
\newcommand{\cqfd}{\hfill\rule{0.35em}{0.35em}}
\DeclareMathOperator{\supp}{supp}

\newcommand{\eg}{\textit{e.g.} }

\title{Free Turbulence 
on $\R^3$ and $\mathbb{T}^3$.}
\author{Fran\c{c}ois Vigneron\footnotemark{}}

\begin{document}
\maketitle
\begin{center}\small\it
${}^\ast$Université Paris-Est,
Laboratoire d'Analyse et de Mathématiques Appliquées, UMR 8050 du CNRS\\
61, avenue du Général de Gaulle,
F-94010 Créteil -- France
\end{center}

\begin{abstract}
The hydrodynamics of Newtonian fluids has been the subject of a tremendous amount
of work over the past eighty years, both in physics and mathematics. Sadly, however, 
a mutual feeling of incomprehension has often hindered scientific contacts. 

This article provides a dictionary that allows mathematicians (including the author)
to define and study the spectral properties of Kolmogorov-Obukov turbulence in a simple
deterministic manner. 
In other words, this approach fits turbulence into the mathematical framework of studying
the qualitative properties of solutions of PDEs, independently from any a-priori
model of the structure of the flow.

To check that this approach is correct, this article proves some of the classical
statements that can be found in physics textbooks.
This is followed by an investigation of the compatibility between turbulence and
the smoothness of solutions of Navier-Stokes in 3D, which was the initial motivation
of this study.
\end{abstract}

\bigskip
The simplest  model of a Newtonian fluid is an incompressible
flow evolving freely with constant density and temperature.
Let us therefore consider the incompressible Navier-Stokes system on $\R_+\times \Omega$
with either $\Omega\subseteq\R^3$ or~$\Omega=\mathbb{T}^3$ :
\begin{equation}\label{NS}
\begin{cases}\displaystyle
\partial_t u - \nu\Delta u + (u\cdot\nabla)u = -\nabla\Pi\\
\div u =0\\
u_{\vert t=0} = u_0, \quad u_{\vert\partial\Omega}=0.
\end{cases}
\end{equation}
Here $u$ represents the velocity field and
$\Pi=p/\rho$, where $p$ is the pressure and $\rho$ is the density assumed 
to be constant and $\nu>0$ is the kinematic viscosity.
This equation has weak solutions (called Leray solutions~\cite{L34}) in the Leray space
$$u\in \mathcal{L}(\Omega) =L^\infty(\R_+;L^2(\Omega))\cap L^2(\R_+;\dot{H}^1(\Omega)),
\qquad \Pi\in L^{5/3}_\text{loc}(\R_+\times\Omega).$$
Pressure can be computed by solving
$-\Delta \Pi = \operatorname{Tr}({}^t\nabla u\cdot\nabla u)$.
In general, the smoothness of $(u,\Pi)$ is an open problem that has been
extensively studied. Among historical landmarks regarding
smoothness, one must cite the works \cite{Serrin}, \cite{FujitaKato} and \cite{KT01}
for the point of view of partial differential equations
and \cite{Cafarelli-Kohn-Nirenberg} for an approach based on
geometric-measure theory.
To get a more comprehensive survey of what is currently known about~\eqref{NS},
one should check  \textit{e.g.} \cite{Lemarie}, 
\cite{Lions}, \cite{Temam} and the references therein.

\bigskip
Since the seminal works of A.N.~Kolmogorov \cite{K41a}, \cite{K41b}, \cite{K41c},
\cite{K41d} and A.M.~Obukhoff \cite{Obukov},
a vast amount of effort has been put into understanding turbulence. 
In physics, one should definitely quote \cite{Batchelor} and
\cite{Frish} as major reference handbooks. 
Personally, I was very impressed by experiment \cite{LA06}, in which immersed floats equipped with GPS devices were allowed to drift in a Canadian river; the speed could be measured directly and a sufficient amount of data could be gathered to check the spectral 5/3 law (see Prop.~\ref{K41b} below) with striking precision.
I discovered the point of view of engineers in \cite{Davidson} and was pleased to
realized that they pay great attention to mathematical rigor because disregarding the
divergence of an integral can trigger a catastrophe in real life.
The following books and articles helped me acquire the experimental background
necessary to write this article
 \cite{BMONMF84}, \cite{BP97}, \cite{BO97}, 
  \cite{brachet2}, \cite{McComb},
  \cite{Sreenivasan99}, \cite{Majda},  \cite{Hopf}.
In mathematics, the question of turbulence often raises cynical reactions.
However, the books \cite{FMRT}, \cite{Lesieur} and the following
works were a valuable source of inspiration~:
\cite{Constantin97}, \cite{LJS97}, \cite{CNT99}, \cite{Kuksin},  \cite{Bardos01}, \cite{FMRT01},
\cite{Constantin06}, \cite{Cheverry} (by publication date).

\bigskip
On the question of whether the spectral properties of turbulence
are compatible with Navier-Stokes, there is a very interesting recent paper
\cite{Craig-Biryuk} dealing with weak solutions. The authors obtain compatibility conditions
that they themselves qualify as being reasonably satisfied (the upper-bounds on the
inertial range largely exceeds the range predicted by physics and observed experimentally).
The present article is an independent work, though
my motivations are similar. 
Starting from a more precise definition of turbulence makes it possible
to recover the real inertial range. Later on I focus on smooth solutions
and prove the exponential decay of the spectrum, 
which will lead to much stronger restrictions that
still allow smooth turbulence to exist.

\section*{Structure of the article, ideas and main results.}
Section \S\ref{BASIC} contains useful definitions and notations. In particular
\S\ref{DEF} provides a deterministic definition of K41-turbulence and a dictionary
to translate physics claims into mathematical statements.
Section \S\ref{RANGE} checks that it indeed leads to the classical properties
of the inertial range and of the energy spectrum that one can find in physics
textbooks.

The next step is to investigate a-priori bounds of the energy spectrum
respectively for low frequencies in~\S\ref{LOW} and high frequencies in~\S\ref{HIGH}.
The low-frequency bounds happens to answer a physical conjecture on what triggers
Batchelor's and Saffman's spectra ; the answer is strongly connected with the spatial
localization and temporal decay properties that have been extensively studied by
mathematicians in the past decade.

Section \S\ref{IFTHM}  contains the discussion of whether
solutions of \eqref{NS} in 3D can be turbulent and what conditions must be satisfied.
A general restriction applies to the time range on which averages are taken.
For smooth solutions, one can prove a lower bound of the fluctuation between
the dissipation rate and its average (a physical phenomenon known as intermittency).
There is also a new formula relating the analyticity radius to the size of the finest
scales in the inertial range.

As some other important statements of physics textbooks have not yet been rigorously
established, section~\S\ref{INTERACT} collects some open problems and hints on
how they could be tackled. 

\bigskip
The first main idea carried by this paper is that qualitative properties of turbulent
flows can be studied with the \textsl{deterministic} tools of PDEs. A probabilistic
approach might still be necessary later on to prove that ``most'' flows are turbulent,
but that problem should be addressed separately.

The second idea is that turbulence is not based on the failure of smoothness
because one can prove it to be \textsl{compatible}
with the analytic regularity of solutions. Turbulence is a specific mathematical
structure of solutions that local regularity methods fail to capture even though
they might be troublingly close (see~\S\ref{parTROUBLING}).
In all likelihood, turbulence will prove to be the key to understanding global smoothness.

The third and more concrete contribution of this article is the introduction of
the ``volume'' function $\operatorname{Vol}(u;[T_0,T_1])$ at the foundations of the theory.
This quantity describes the large scales involved in turbulence. Later on, the new
time scale $\mathcal{T}(u_0;\omega)$ is shown to be characteristic of time intervals on
which free turbulence can be observed. Their ubiquitous nature in this article
suggests that they should be given some attention.

Finally, this article sheds light on some \textsl{subtleties} related to the definition of the
energy spectrum in the discrete case. Luckily enough, spectral computations could
be carried out explicitly on $\mathbb{T}^3$ but the generalization to other
domains should be done very carefully as it might be responsible for a
substantial part of the troubles of ``real life'' turbulence.

\setcounter{tocdepth}{2}
\tableofcontents

\section*{Thanks.}
I would like to express my gratitude to F.~Debbasch and M.~Paicu for fruitful
discussions and the participants of Paris-Sud Orsay's PDE seminar of 01/28/2010, especially
S.~Alhinac, for valuable remarks that helped me improve the final redaction.
A special mention to C.~Bardos, W.~Craig and V.~Sverak for their interest in this paper
and their comments on the preprint.

\section{Definitions and basic properties}
\label{BASIC}

Let us recall some mathematical notations and physical definitions
regarding \eqref{NS}.

\subsection{Kinetic energy and energy spectrum}
The \textbf{kinetic energy} at time $t>0$ is
\begin{equation}
E(t)=\rho \norme[L^2]{u(t)}^2.
\end{equation}

\begin{subequations}
The \textbf{energy spectrum}  represents the contribution to the total kinetic energy
of the frequency~$K$. It is defined rigorously by the spectral resolution of the
Stokes operator $A=-\mathbb{P}\Delta$ which is a positive self-adjoint operator with
domain $D(A)=\left\{u\in H^2(\Omega)\cap H^1_{0}(\Omega) \,;\, \div u=0\right\}$
if $\Omega$ is smooth
(for non-smooth domains, see \eg \cite[p.~7, p.~128]{S01} or  \cite{Galdi}).
Here $\mathbb{P}$ denotes the
orthogonal projection on divergence-free vector fields in $L^2(\Omega)$.

\subsubsection{Case of a continuous Stokes spectrum}

On $\Omega=\R^3$ one has $\sigma(A^{1/2})=[0,\infty)$.
Thus the (isotropic) energy spectrum of a function $u\in\mathcal{L}(\R^3)$ is defined by~:
\begin{equation}\label{DEFSPECTRUM}
E^\ast(K,t)=\frac{d}{dK}\left(\rho \norme[L^2]{\chi\left(\frac{A^{1/2}}{K}\right)u(t,\cdot)}^2\right)
\qquad\text{with}\qquad
\begin{cases}
\chi\in C^\infty(\R_+;\R_+),\\
\chi(r)=1 & r<1/2,\\
\chi(r)=0 & r>2.
\end{cases}
\end{equation}
We use a smooth cut-off  to ensure that $E^\ast$ exists for any $u(t)\in L^2(\R^3)$.
The non smooth cut-off is defined as the limit in the distribution sense
$\chi^2(r) \rightharpoonup  \mathbbm{1}_{r<1}$.
The spectrum satisfies the fundamental property~:
\begin{equation}\label{ENERGY}
E(t)=\int_0^\infty E^\ast(K,t) \,dK
\end{equation}
\end{subequations}
The projector $\mathbb{P}$
commutes with derivations, which allows one to compute the spectrum
explicitly in Fourier variables (see \eg \cite[pp.~38-40]{CDGG06}).
\begin{prop}\label{EXPLICITSPECTRUM}
Let us compute the Fourier transform with $\hat{u}(t,\xi)=\int_{\R^3}e^{-ix\cdot\xi}u(t,x)dx$.
\begin{enumerate}
\item
If $u$ is a divergence-free vector field in $\mathcal{\R}^3$~:
\begin{equation}\label{E*}
E^\ast(K,t)= \frac{\rho}{K}\int_{\R^3}\psi\left(\frac{|\xi|}{K}\right) |\hat{u}(t,\xi)|^2 \frac{d\xi}{(2\pi)^3} 
\end{equation}
where $\psi$ is a positive smooth bump function,
supported on $[2^{-1},2]$ and such that $\int_0^\infty \psi(r)\frac{dr}{r}=1$.
\item
In the limit of non smooth cut-off, one has  $\psi(r)\rightharpoonup\delta_{r=1}$ and
\begin{equation}\label{EdaggerR3}
E^\ast(K,t)\to E^\dagger(K,t) \underset{\text{def}}{=}
(2\pi)^{-3}\rho K^2\int_{\mathbb{S}^2} |\hat{u}(t,K\vartheta)|^2 d\vartheta.
\end{equation}
Conversely, if one defines the ``experimental'' value of $E^\dagger(K,t)$ as the average
over a spherical shell of relative amplitude $\delta\in]0,1[$, one finds 
$$\displaystyle
\frac{1}{2\delta K}\int_{(1-\delta)K}^{(1+\delta)K} |E^\dagger(\kappa,t)|^2d\kappa=E^\ast(K,t)$$
for~${\psi(r)=\frac{1}{2\delta} \mathbf{1}_{[1-\delta,1+\delta]}(r)}$ and
$\psi(r)\rightharpoonup\delta_{r=1}$ as $\delta\to0$.

\end{enumerate}
\end{prop}

\Proof
Applying the spectral theorem, one has $A=\int_0^\infty \lambda dP_\lambda$ 
with $P_\lambda = \mathbb{P}\circ P^D_\lambda$ where $P^D_\lambda$
is the spectral projector associated to the Dirichlet operator $(-\Delta)_D$ on $\Omega$.
For $\R^3$ and $\mathbb{T}^3$, one has $P_\lambda^D=\mathcal{F}^{-1}
\circ\widehat{P}^D_\lambda\circ\mathcal{F}$ where $\mathcal{F}$ is the
Fourier transform and respectively ($\xi\in\R^3$
and $k\in\mathbb{Z}^3$)~:
$$\widehat{P}^D_\lambda = \mathbbm{1}_{|\xi|^2\leq \lambda} \qquad\text{or}
\qquad \widehat{P}^D_\lambda = \mathbbm{1}_{|k|^2\leq \lambda}.$$
In both cases, one has $\mathbb{P}\circ P^D_\lambda = P^D_\lambda\circ\mathbb{P}$ hence~:
$$\chi\left(\frac{A^{1/2}}{K}\right) u=\mathcal{F}^{-1}\circ
\int_0^\infty \chi\left(\frac{\lambda^{1/2}}{K}\right) d\widehat{P}_\lambda\circ\mathcal{F}(u)=
\mathcal{F}^{-1}\left[\chi\left(\frac{|\cdot|}{K}\right)\widehat{\mathbb{P}u}\right].$$
If $\div u=0$, one has $\mathbb{P}u=u$.
Then \eqref{E*} follows from Parseval identity $\norme[L^2(\R^3)]{v}^2=(2\pi)^{-3}\norme[L^2]{\hat{v}}^2$.
The rest of the statement is obvious with $\psi(r)=-2r\chi(r)\chi'(r)$.
\cqfd\endProof
\paragraph{Remarks}
\begin{enumerate}
\item
To unify notations with the case of a discrete Stokes spectrum, let us
state~\eqref{ENERGY} as
\begin{equation}
E(t)=\int_{\sigma(A^{1/2})} E^\dagger(K,t) \,d\mu(K)
\end{equation}
where $E^\dagger$ is defined by~\eqref{EdaggerR3} and
$\mu$ is the Lebesgue measure on $\sigma(A^{1/2})=[0,\infty)$.
Note that by Fubini's Theorem, the spectrum $E^\dagger(K,t)$ exists
in $L^1(\R_+)$ for any $u\in L^2(\R^3)$.
\item
In the following, one should not rely on any other norm
than $\norme[L^1]{\psi}\leq 2$ because other $L^p$ norms of $\psi$ are unbounded in the
limit $\psi(r)\rightharpoonup\delta_{r=1}$. Estimates like
$E^\ast(K,t)\leq \norme[L^\infty]{\psi} \frac{E(t)}{K}$ should be disregarded as
empty of any physical meaning and because
they do not correspond to any property of $E^\dagger$.
\end{enumerate}

\subsubsection{Case of a discrete Stokes spectrum}

On $\mathbb{T}^3=\R^3/(L\mathbb{Z})^3$ the spectrum of the Sokes operator is discrete
$\sigma(A) = \{ K^2 \,;\, K\in\Sigma\}$ where
\begin{equation}
\Sigma = \sigma(A^{1/2}) =\{|\mathbf{k}| \,;\, \mathbf{k}\in (2\pi L^{-1}\mathbb{Z})^3\}
= \{2\pi L^{-1}\sqrt{n} \,;\, n\in\Box_3\}
\end{equation}
and $\Box_3$ is the set of integers that are the sum of three squares.
Let us denote by $\Sigma^\ast = \Sigma\backslash\{0\}$.
By the  Gauss-Legendre three-squares theorem~:
\begin{equation}
\Box_3=\{n\in\mathbb{N}, \enspace n\neq 4^p(8q+7)\}.
\end{equation}
Writing an analogue to formula  \eqref{DEFSPECTRUM} must be done
carefully in order to cope with the discrete differentiation.

\medskip
Let us denote by $\hat{u}(t,\mathbf{k})=\int_{\mathbb{T}^3} e^{-i\mathbf{k}\cdot x} u(t,x) dx$
the $k^{\text{th}}$ Fourier coefficient of $u(t)$ and by
$L=\operatorname{Vol}(\mathbb{T}^3)^{1/3}$ 
the characteristic length of $\mathbb{T}^3$. 

\begin{dfn}
On $\mathbb{T}^3$ the energy spectrum is defined by~:
\begin{equation}\label{E*PERIODIC}
\forall K\in \Sigma^\ast,\qquad
E^\dagger(K,t)= (2\pi)^{-2} \rho \left(\frac{K}{L}\right)
\sum_{\substack{\mathbf{k}\in(2\pi L^{-1}\mathbb{Z})^3\\|\mathbf{k}|= K}} |\hat{u}(t,\mathbf{k})|^2.
\end{equation}
One defines the following measure on $\Sigma^\ast$~:
\begin{equation}
\mu = 2\pi L^{-1}\sum_{K\in\Sigma^\ast} \left(\frac{2\pi}{K L}\right) \delta_K.
\end{equation}
For  $\delta\ll K$, the corresponding average value on the spherical shell
$\Sigma_\delta(K)=\left\{\kappa \in \Sigma^\ast \,;\,
\left|\frac{|\kappa|}{K} -1\right|\leq \frac{\delta}{K}\right\}$ is~:
\begin{equation}\label{EEXPPERIODIC}
E^\ast_\delta(K,t)=
\frac{\displaystyle \int_{\Sigma_{\delta}(K)} E^\dagger(\kappa,t) d\mu(\kappa)}
{\mu(\Sigma_{\delta}(K))}
=\frac{\displaystyle 2\pi L^{-1}\sum_{\kappa\in\Sigma_\delta(K)}
\left(\frac{2\pi}{\kappa L}\right) E^\dagger(\kappa,t)}
{\displaystyle 2\pi L^{-1}\sum_{\kappa\in\Sigma_\delta(K)} \left(\frac{2\pi}{\kappa L}\right)}
\cdotp
\end{equation}
One will call $E^\ast_\delta(K,t)$ the ``experimental'' value of $E^\dagger(K,t)$.
\end{dfn}
\begin{prop}
The following statements hold.
\begin{enumerate}
\item
To compute the total energy, formula \eqref{ENERGY} is replaced by~:
\begin{align}\label{ENERGYbis}
E(t)&=\left(\int_{\mathbb{T}^3} \rho u_0(x) dx\right)^2 +
\int_{\sigma(A^{1/2})\backslash\{0\}} E^\dagger(\kappa,t) d\kappa
\\\notag&=\left(\int_{\mathbb{T}^3} \rho u_0(x) dx\right)^2 +
2\pi L^{-1}\sum_{K\in\Sigma^\ast} \left(\frac{2\pi}{K L}\right)E^\dagger(K,t).
\end{align}
\item
There exist $\beta\in]0,1[$ and $C\geq1$ such that
\begin{equation}\label{EXPERIMENTAL}
0<\delta \leq \beta K  \quad\Longrightarrow\quad
C^{-1}E^\ast_\delta(K,t) \leq
\frac{\rho}{4 \delta}\left(\frac{1}{L^3} \sum_{k\in S_{\delta}(K)} |\hat{u}(t,k)|^2\right)
\leq C E^\ast_\delta(K,t)
\end{equation}
where $S_\delta(K)=\left\{k\in (2\pi L^{-1}\mathbb{Z})^3 \,;\, \left|\frac{|k|}{K} -1\right|\leq \frac{\delta}{K}\right\}$ is the spherical shell of frequencies $K\pm \delta$.
\end{enumerate}
\end{prop}
Identity \eqref{EXPERIMENTAL} is crucial to match the theory to real-world experiments.
Indeed, the sum of the squares of Fourier coefficients on spherical shells \textbf{is} the
energy spectrum of all numerical and physical experiments. Therefore, the universal
behavior observed for the $E^\ast_\delta$ of turbulent flows can be addressed
mathematically by investigating the corresponding property on~$E^\dagger$.
\paragraph{Remarks}
\begin{enumerate}
\item\label{RMKDISCRETESPECTRUM}
One cannot emphasize enough that \eqref{EXPERIMENTAL} does \textbf{not} hold
for any other normalization than \eqref{E*PERIODIC}--\eqref{ENERGYbis}.
For example, using the analogy with \eqref{EdaggerR3}
one could be tempted to replace \eqref{E*PERIODIC} by~:
$$(2\pi)^{-3} \rho K^2\: \sum_{\substack{L k\in2\pi\mathbb{Z}^3\\|k|= K}} |\hat{u}(t,k)|^2$$
This choice would lead to a catastrophe. First, the averages on
spherical shells would be equivalent to
$$K\times \left(\frac{\rho}{\delta L^2} \sum_{k\in S_{\delta}(K)} |\hat{u}(t,k)|^2\right)$$
which is not the usual normalization of experimental spectra.
If this fact remains unnoticed and one develops the rest of the theory,
then in Theorem~\ref{RANGETHM} one would not recover the usual
Kolmogorov dissipation frequency $K_d$ (defined below) but instead $K_d'=\frac{\bar{\varepsilon}\operatorname{Vol}(\mathbb{T}^3)}{(\alpha\nu)^3}$.
Physics textbooks usually dodge the subtlety of \eqref{E*PERIODIC}--\eqref{ENERGYbis} ;
at best, they use \eqref{EdaggerR3} to define~$E^\ast(K,t)$ in the case of~$\R^3$ and
rely on~\eqref{EXPERIMENTAL} for any practical purposes.
This subtlety is however of great practical importance
as most experimental data is obtained in a situation where the spectrum is discrete\ldots
\item\label{RMKDISCRETESPECTRUM2}
The key to this computation is the asymptotic of $\Sigma^\ast$ \textit{i.e.}~of
eigenvalues counted without multiplicity.
Let us reorder the Stokes spectrum in an increasing
sequence~$\sigma(A)=\{K_j^2\,;\,j\in\mathbb{N}\}$ with $K_j<K_{j+1}$.
One can prove the existence of $C,C'\geq1$ such that~:
\begin{equation}\label{SPECTRALASSYMPTOTIC}
\frac{C^{-1}}{L^2} \leq K_{j+1}^2-K_j^2 \leq \frac{C}{L^2}
\qquad\text{or equivalently}\qquad
\frac{{C'}^{-1} L^{-2}}{K_j} \leq K_{j+1}-K_j \leq \frac{C' L^{-2}}{K_j}
\end{equation}
This asymptotic is responsible for the fact that~$E^\dagger$
(or more generally any function on $\sigma(A^{1/2}))$
has a different normalization than its shell averages $E^\ast_\delta$.
Describing which domains satisfy \eqref{SPECTRALASSYMPTOTIC} on the
pertinent range of frequencies is a widely open problem whose answer might
actually provide a hint as to why physicists have trouble unifying the
description of all turbulences.
Conversely, it might also explain why some common patterns have been found for
each of the ``families'' of turbulence
(\textit{e.g.}~grid turbulence, wake turbulence, jet turbulence,
boundary layer turbulence or turbulence in bounded domains).
\item The experimental spectrum is defined by averages on
spherical shells of frequencies $K\pm \delta$ with~$\delta\ll K$. For dyadic shells
\textit{i.e.} frequencies $(1\pm \epsilon)K$ the corresponding average
is $$E^\ast_{\epsilon K}(K,t)\simeq
\frac{\rho}{\epsilon K}\left(\frac{1}{L^3} \sum_{1-\epsilon<|k|/K<1+\epsilon} |\hat{u}(t,k)|^2\right).$$
\end{enumerate}

\begin{figure}[ht!]\sf
\begin{minipage}{\textwidth}
\begin{multicols}{2}
\begin{center}
\includegraphics[width=.37\textwidth]{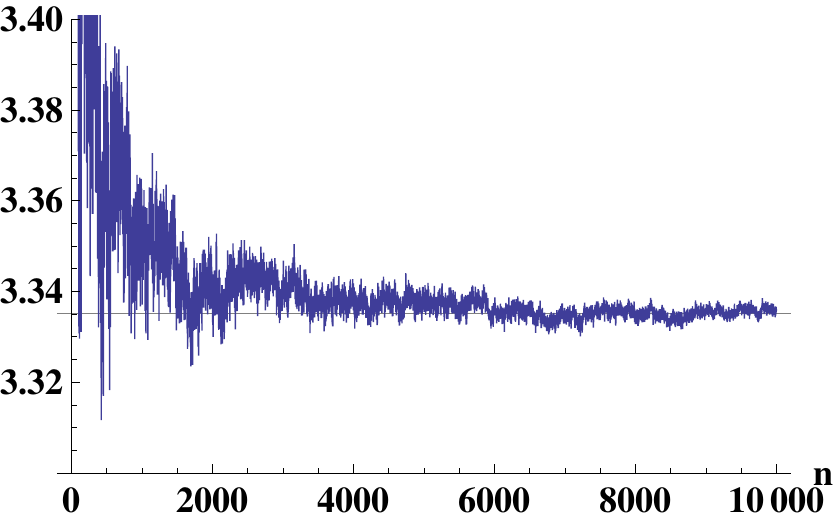}
\end{center}
\caption{\label{NBTHEXO}}
Plot of $\frac{10}{\sqrt{n}}\sum \frac{1}{\sqrt{n'}}$ for $n'\in\Box_3$ such that~:
$$|\sqrt{n'}-\sqrt{n}|\leq \frac{1}{10}\sqrt{n}.$$
On this range, one has~:
$$3.3\leq \frac{10}{\sqrt{n}}\sum \frac{1}{\sqrt{n'}}\leq 16.2$$
The numerical value of the asymptotic equivalent is $3.335$.
The corresponding sum without the restriction $n'\in\Box_3$ is equivalent to~4 as $n\to\infty$.
\end{multicols}
\end{minipage}
\end{figure}

\Proof
Parseval indentity  $\norme[L^2(\mathbb{T}^3)]{v}^2=L^{-3}\sum |\hat{v}(k)|^2$ dictates
the compatibility between formulas \eqref{E*PERIODIC} and \eqref{ENERGYbis}.
Note that \eqref{NS} implies $\displaystyle
\int_{\mathbb{T}^3} \rho u(t,x) dx = \int_{\mathbb{T}^3} \rho u_0(x) dx$ but
the corresponding energy is not seen by the spectrum.
Parseval identity also ensures that the numerator of $E^\ast_\delta(K,t)$ is
$$2\pi L^{-1}\sum_{\kappa\in\Sigma_\delta(K)}
\left(\frac{2\pi}{\kappa L}\right) E^\dagger(\kappa,t)
=\rho\left(\frac{1}{L^3} \sum_{k\in S_{\delta}(K)} |\hat{u}(t,k)|^2\right).$$
To prove \eqref{EXPERIMENTAL} one has only to check that the denominator
$$\left(\frac{2\pi}{L}\right)^2\sum_{\kappa\in\Sigma_\delta(K)} \frac{1}{\kappa}=
\frac{2\pi}{L} \sum_{\substack{n'\in\Box_3\\
|\sqrt{n'}-\sqrt{n}|\leq \frac{\delta L}{2\pi}}} \frac{1}{\sqrt{n'}}
\qquad\text{where}\qquad K=\frac{2\pi}{L}\sqrt{n}$$
is equivalent to $\delta$, which is an easy exercise in number theory (see the proof of
Theorem~\ref{RANGETHM} below, where a similar computation is fully detailed).
The numerical value of the constant is illustrated in Figure~\ref{NBTHEXO}.
The key ingredient is that $\mathbb{N}\backslash\Box_3$ contains only integers
$n\equiv0,4 \text{ or }7 \text{ mod }8$, which in turn ensures that
$\sum_{n\in[1,N]\cap\Box_3} n^{-s}$ and
$\sum_{n=1}^N n^{-s}$ are equivalent up to a numerical factor.
This exact same property of congruence modulo 8 also implies that
\begin{equation}
n\in\Box_3 \quad\Longrightarrow\quad \Box_3\cap\{n+1,n+2,n+3\}\neq\emptyset.
\end{equation}
The asymptotic \eqref{SPECTRALASSYMPTOTIC} of the Stokes spectrum follows
immediately as $\sqrt{n+j}-\sqrt{n}=\frac{j}{2\sqrt{n}}+O\left(\frac{j^2}{n^{3/2}}\right)$.
\cqfd\endProof

\subsection{Dissipation of kinetic energy}
The \textbf{dissipation rate} at time $t$ is defined by~:
\begin{equation}\label{EPSILON-DEF}
\varepsilon(t) = 2\rho\nu \norme[L^2]{\nabla u(t)}^2.
\end{equation}
One can check immediately in Fourier variables that~:
\begin{subequations}\label{EPSILON}
\begin{align}
\frac{1}{2} \varepsilon(t)\leq 2\nu \int_0^\infty  K^2\,E^\ast(K,t)\, dK\leq 2\varepsilon(t)
\qquad&\text{on }\R^3,
\\
\varepsilon(t) \enspace= \enspace 2\nu \int_{\sigma(A^{1/2})} K^2E^\dagger(K,t) d\mu(K)
\qquad &\text{on }\mathbb{T}^3.
\end{align}
\end{subequations}
Note that $2\norme[L^2]{\nabla u(t)}^2 = \norme[L^2]{\omega(t)}^2$ if $u$ is a square integrable divergence-free vector field on $\R^3$ or $\mathbb{T}^3$ with
vorticity $\omega=(\partial_i u_j-\partial_ju_i)_{1\leq i,j\leq 3}$.
Because of the so called ``stretching term'' in the right-hand side of
$$\partial_t \omega - \nu \Delta \omega +
(u\cdot\nabla)\omega = 
(\omega\cdot\nabla)u
$$
no good a-priori estimate of $\omega$ or $\varepsilon$ is known
(check \cite{C09} for partial results).
One has only~:
\begin{equation}
\varepsilon(t)= \varepsilon(0)-2\rho\nu^2\int_0^t \norme[L^2]{\nabla \omega(\tau)}^2 d\tau
-4\rho\nu\sum_{i,j,k} \int_0^t\int_{\Omega} (\partial_i u_j)(\partial_k u_i )(\partial_k u_j).
\end{equation}
Gronwal inequality provides exponential estimates like : 
$\varepsilon(t)\leq C \varepsilon(0) \exp\left(\int_0^t\norme[L^\infty]{\nabla u(t')}dt'\right)$.

\paragraph{Global balance of energy.}
\begin{subequations}\label{ENERGYCONSERVATION}
If $u$ is a smooth solution of \eqref{NS}, one has
\begin{equation}
\varepsilon=-\frac{dE}{dt}
\end{equation}
which justifies the name of  ``energy dissipation rate'' for $\varepsilon$.
As $\varepsilon$ depends only on $\omega$
it means that the dissipation of kinetic energy occurs
exclusively through vortex structures.
For Leray solutions, one has only for a.e. $t>t'$~:
\begin{equation}\label{LERAYINEQ}
E(t)\leq E(t')-\int_{t'}^t \varepsilon(\tau)d\tau.
\end{equation}
As observed in \cite{DR00} the possible lack of smoothness (\textit{i.e.}~a strict inequality
in \eqref{LERAYINEQ})
would mean that an extraordinary dissipation
has occured between $t$ and $t'$.
The discussion of whether \eqref{NS} would remain a good physical model 
on such a~$[t,t']$ (or what model should replace it) is beyond the scope of this article.
\end{subequations}

\subsection{Scaling transformation to work per unit of mass}

One can construct two families of transformations that preserve \eqref{NS} without
changing the kinematic viscosity $\nu$.
\begin{itemize}
\item
If $(u,\Pi)$ is a solution of \eqref{NS}, then so is~:
\begin{subequations}\label{SCALING}
\begin{equation}\label{SCALINGa}
\forall\lambda>0, \qquad
u_\lambda(t,x)=\lambda u(\lambda^2 t,\lambda x),
\qquad
\Pi_\lambda(t,x)=\lambda^2 \Pi(\lambda^2 t,\lambda x).
\end{equation}
In the case of $\mathbb{T}^3$, the new domain becomes $\mathbb{T}_\lambda^3=\R^3/(\lambda^{-1} \mathbb{Z})^3$ and to ensure the conservation of the total mass,
the new density must be~:
\begin{equation}\label{SCALINGb}
\rho_\lambda=\lambda^3 \rho.
\end{equation}
\end{subequations}
\item
$(u,p,\rho)\mapsto (u,\mu p,\mu\rho)$  for $\mu>0$ is another family of solutions.
This transformation means that the mass of each particle is multiplied
by a factor $\mu$ without changing
the number of particles in the fluid\footnote{Like replacing a hydrogen flow by a
helium flow with the same velocity field and assuming that all other physical properties,
including~$\nu$, will remain identical.  However these particular gases are compressible
so \eqref{NS} does not describe them properly.}.
In the case of~$\mathbb{T}^3$, this transformation
would obviously change the total mass by a factor $\mu$.
\end{itemize}
\begin{dfn}
In the following, one shall work ``per unit of mass'' \textit{i.e.} for a given number of molecules.
For a bounded domain $\Omega$ or $\mathbb{T}^3$, it
means that one choses $\displaystyle\mu= \left(\int_{\Omega} \rho  dx\right)^{-1}$
and is left  with 
\begin{subequations}
\begin{equation}
\rho=\operatorname{Vol}(\Omega)^{-1}.
\end{equation}
For $\R^3$, one will extend this by convention by letting~:
\begin{equation}
\rho=1 \:[\text{length}]^{-3}.
\end{equation}
\end{subequations}
\end{dfn}
Once one works per unit of mass, one cannot apply the second family of
transformations anymore.
The first one  \eqref{SCALINGa}-\eqref{SCALINGb} remains
admissible with a dimensionless parameter $\lambda>0$. It acts in the following way on the energy-related quantities~:
\begin{equation}
E_\lambda(t) = \lambda^2 E(\lambda^2 t),\quad
E_\lambda^\ast(K,t)=\lambda E^\ast\left(\frac{K}{\lambda},\lambda^2 t\right),\quad
\varepsilon_\lambda(t)=\lambda^4 \varepsilon(\lambda^2 t)
\end{equation}
and  $\norme[Y]{u_\lambda}=\norme[Y]{u}$ 
for  $Y=L^\infty(\R_+;X)$ and e.g.~$X=L^3(\R^3)$, $\dot{H}^{1/2}(\R^3)$ or $BMO^{-1}(\R^3)$
or any of the so called ``scaling invariant'' function spaces.

\begin{figure}[ht!]\sf
\begin{minipage}{\textwidth}
\begin{multicols}{2}
\begin{center}
\includegraphics[width=.37\textwidth]{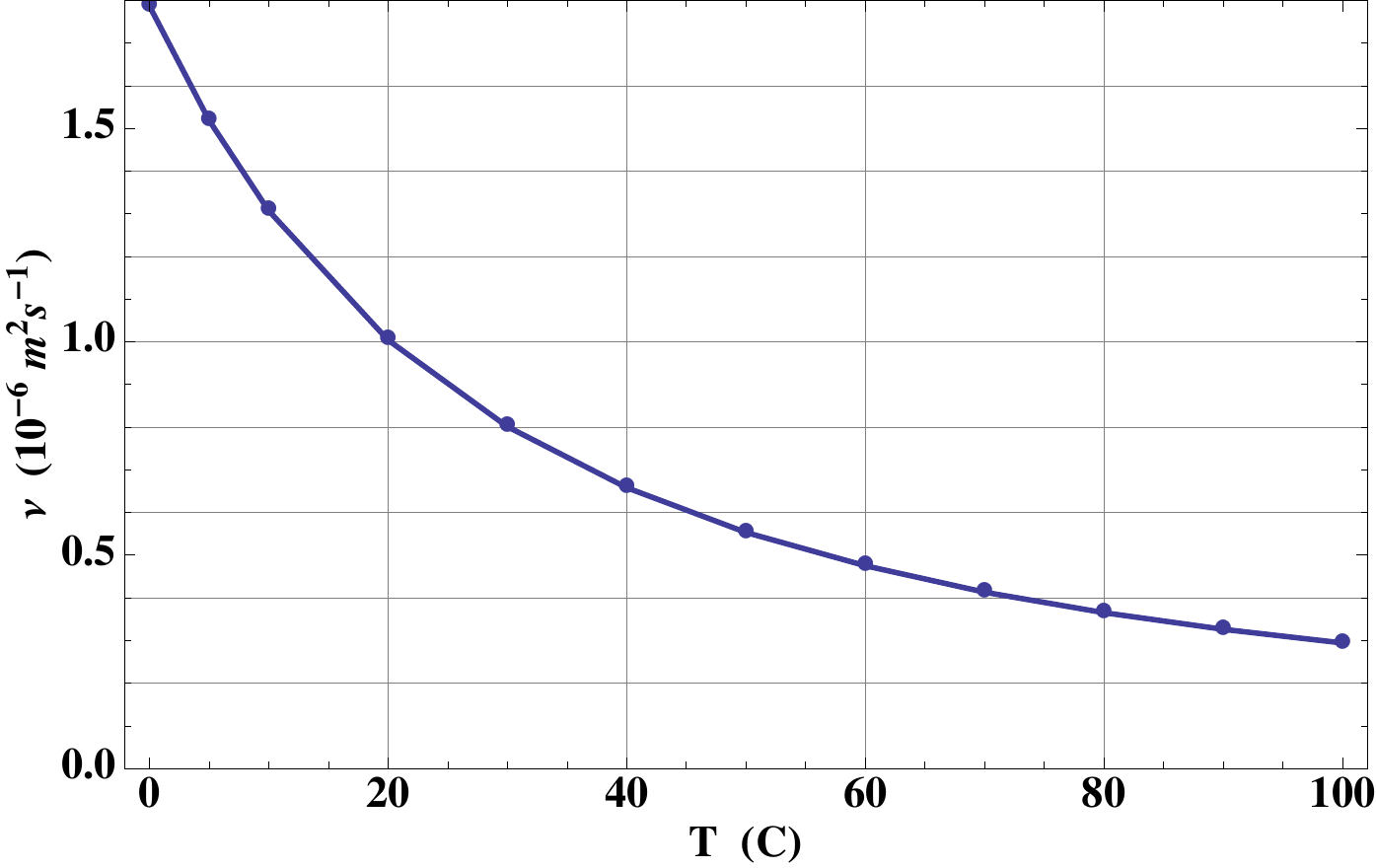}
\end{center}
\columnbreak
\vspace*{\stretch{1}}
\caption{\label{NUWATER} Real-world illustration -- kinematic viscosity of liquid water
as a function of temperature (see \textit{e.g.} \cite{viscosity} and the references therein)
; the order of magnitude is $\nu\simeq 10^{-6} \: \text{m}^2 \text{s}^{-1}$.}
\vspace*{\stretch{1}}
\phantom{.}
\end{multicols}
\end{minipage}
\end{figure}

\paragraph{Physical dimensions.}

Let's recall that $\nu=[\text{length}]^2\cdot[\text{time}]^{-1}$ is the kinematic viscosity.
For liquid water, the value of the kinematic viscosity is illustrated by Figure~\ref{NUWATER}.
One can easily check that the natural scaling of energy per unit of mass is
$$E(t)=[\text{length}]^2\cdot[\text{time}]^{-2}.$$
It follows that $E^\ast(t,K)=[\text{length}]^3\cdot[\text{time}]^{-2}$
and $\varepsilon(t)=[\text{length}]^{2}\cdot[\text{time}]^{-3}$.
With the notations of Proposition~\ref{EXPLICITSPECTRUM}, 
the Fourier transform is
$\hat{u}(t,\cdot)=[\text{length}]^4\cdot[\text{time}]^{-1}$ on both
$\R^3$ and $\mathbb{T}^3$,
the modulus of frequencies is $K=[\text{length}]^{-1}$ and the density is
$\rho=[\text{length}]^{-3}$. 
From a mathematical point of view,  working with physical dimensions is equivalent
to checking that each identity is scale-invariant under~\eqref{SCALING}.

\subsection{Time averages and intermittency}\label{TAI}
The function $E^\ast(K,t)$ and $\varepsilon(t)$ are to some extent
accessible to the experiments of fluid
dynamics. And it is a common observation that even if those functions fluctuate a lot,
time averages (over proper time intervals) display universal behaviors.
This observation is the experimental essence of turbulence theory.
Time averages of a function $f(t,\theta)$ on $[T_0,T_1]\times \Theta$  will be defined by~:
\begin{equation}\label{TAVERAGE}
\bar{f}(\theta)  = \frac{1}{\Delta}\int_{T_0}^{T_1} f(t,\theta) \,dt \quad\text{with}\quad
\Delta=T_1-T_0.
\end{equation}
We will also denote $\bar{f}$ by $\langle f \rangle$ if the expression of $f$
is too large and makes the first notation ambiguous.

\begin{dfn}
Once the time interval $[T_0,T_1]$ is given, one defines the mean energy $\bar{E}$,
the mean energy spectrum $\bar{E}^\ast(K)$ and the mean energy dissipation
rate $\bar{\varepsilon}$ according to \eqref{TAVERAGE}. 
\end{dfn}

\bigskip
Substantial fluctuations of $\varepsilon(t)$ away from $\bar\varepsilon$ are an interesting
phenomenon known as \textbf{intermittency}. More subtle definitions are possible~;
one could call this one ``temporal intermittency'' to differentiate it from ``spatial intermittency''
that deals with substantial spatial fluctuations, either at large or small scales.
The following result
provides a simple way to detect intermittency by comparing the average
energy $\bar{E}$ to the linear interpolation between the initial and final energy.
We will need this statement for Theorem~\ref{SMOOTHTURBULENCE}, which establishes
a subtle relationship between intermittency and the smoothness of turbulent flows.
\begin{prop} If $u$ is a smooth solution of~\eqref{NS} on $[T_0,T_1]$ then~:
\begin{equation}\label{SMOOTH-EPSILON}
\bar\varepsilon = \frac{E_0-E_1}{\Delta}
\end{equation}
with $E_i=E(T_i)$, $\Delta=T_1-T_0$ and
\begin{equation}\label{INTERMITTENCY}
\left|\bar{E}-\frac{E_0+E_1}{2}\right|\leq
\int_{T_0}^{T_1} |\varepsilon(t) - \bar{\varepsilon}| dt.
\end{equation}
\end{prop}
\Proof
As $u$ is smooth, the balance of energy reads $E(t) = E(t')-\int_{t'}^t \varepsilon(\tau)d\tau$
hence $\bar\varepsilon=\frac{E(T_0)-E(T_1)}{\Delta}$.
Let us now integrate this relation on $[T_0,T_0+t]$ for $t\in[0,\Delta]$~:
$$\bar{E}=\frac{1}{\Delta}\int_0^\Delta E(T_0+t)\,dt = E(T_0) -
\frac{1}{\Delta}\int_{T_0}^{T_1} (T_0+\Delta-\tau) \varepsilon(\tau) \, dt d\tau$$
hence~:
$$\bar{E}-\frac{E_0+E_1}{2}=\bar{E} - E(T_0) + \frac{\Delta \cdot  \bar\varepsilon}{2} =
\frac{1}{\Delta}\int_{T_0}^{T_1} (T_1-\tau) (\varepsilon(\tau)-\bar\varepsilon) \, dt d\tau.$$
One concludes using the $L^\infty\ast L^1\to L^\infty$ convolution property.
\cqfd\endProof

\subsection{Average ``volume'' of a function}

In naive terms, one can describe $\text{Vol}(u;[T_0,T_1])$ as an intrinsic
measure of $\{x\,;\,|u(t,x)|>\epsilon\}$ for ``adequate'' $\epsilon$ and
proper time average, \textit{i.e.}
the average volume of the region where $u$ is most intense.
\begin{dfn}
For any measurable function $u(t,x)\in L^2([T_0,T_1]\times\Omega)$,
one defines the \textbf{average volume} occupied by $u$ on $[T_0,T_1]$ by~:
\begin{equation}\label{VOL}
\operatorname{Vol}(u;[T_0,T_1]) = \frac{\langle\norme[L^1(\Omega)]{u}^2\rangle}{\langle\norme[L^2(\Omega)]{u}^2\rangle} \in \overline{\R}_+
\end{equation}
where $\langle\cdot\rangle$ refers to time averages defined by \eqref{TAVERAGE}.
\end{dfn}
The following sections will show the relevance of \eqref{VOL} for turbulence.
Definition \eqref{VOL} provides in a way a substitute to the probabilistic
assumption of ``spatial homogeneity of the flow'' (see Theorem~\ref{RANGETHM} and  \S\ref{PAR:INTERM}) and will also prove to be
well suited to unbounded domains (see \S\ref{LOW},
 Theorem~\ref{SPECTRUMFROML1}).

\begin{figure}[ht!]\sf
\begin{center}
\includegraphics[width=.75\textwidth]{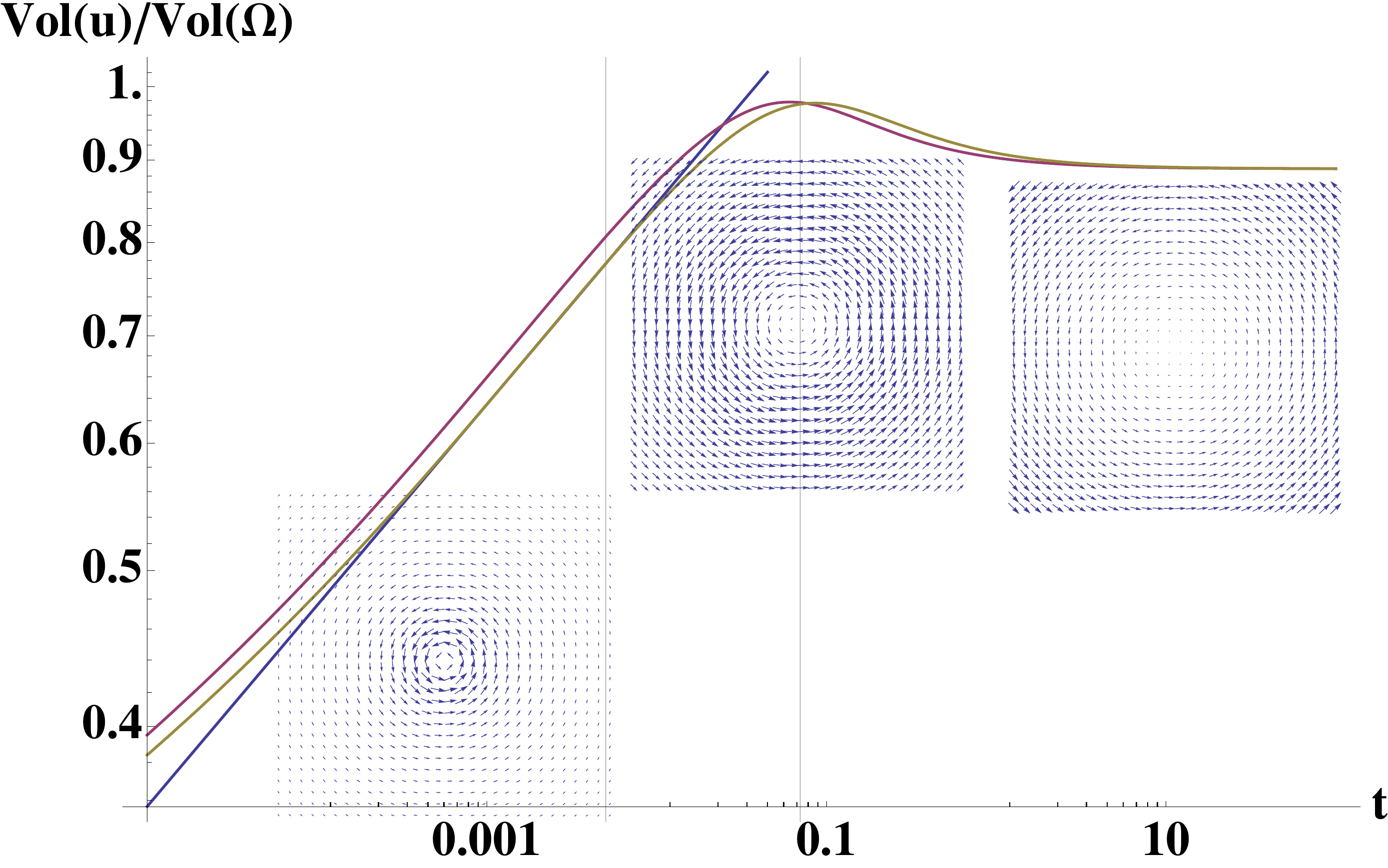}
\end{center}
\caption{\label{figoseen} Plot of $\norme[L^1]{u(t)}^2/\norme[L^2]{u(t)}^2$ and
of~$\operatorname{Vol}(u:[t/2,t])$ for $\Gamma=\nu=1$. The computation is done
in restriction to the cylinder $\Omega=\{(x,y,z) \,;\, x^2+y^2<1, \: 0<z<1\}$ and the result is displayed in Log-Log scale as a percentage of~$\operatorname{Vol}(\Omega)$.
Inlaid pictures represent  the vector field $u_h(t,x)$
for $t\in\{5\times 10^{-3},7\times10^{-2},10^2\}$.
}
\end{figure}

\paragraph{Examples}
\begin{enumerate}
\item
If $\Omega$ is bounded, the Cauchy-Schwarz inequality provides
for any $f\in L^2([T_0,T_1]\times\Omega$~:
$$\operatorname{Vol}(f;[T_0,T_1]) \leq \operatorname{Vol}(\Omega).$$
\item
If $f(t,x)=\mathbbm{1}_{\Omega'+\eta(t)}(x)$
is the characteristic function of a subset $\Omega'$ of $\Omega$ translated
by a vector~$\eta(t)$ such that $\Omega'+\eta(t)\subset \Omega$, then one has
$\operatorname{Vol}(f;[T_0,T_1])=\operatorname{Vol}(\Omega')$ and
$$\operatorname{Vol}(f+\epsilon(1-f);[T_0,T_1])=
\left(\frac{(1+q\epsilon )^2}{1+q\epsilon^2}\right)\operatorname{Vol}(\Omega')
\quad\text{with}\quad
q=\frac{\operatorname{Vol}(\Omega\backslash\Omega')}{\operatorname{Vol}(\Omega')}\cdotp
$$
\item
A simple computation on $\R^3$ gives~:
\begin{equation}\label{EXVOL}
\operatorname{Vol}(e^{\nu t \Delta}\delta_0;[T_0,T_1])=
8 \sqrt{2} \pi^{3/2}  \frac{\nu (T_1-T_0)}{\frac{1}{\sqrt{\nu T_0}}-\frac{1}{\sqrt{\nu T_1}}}.
\end{equation}
In particular the volume is 0 on $[0,T]$ because $\delta_0\not\in L^2$ and infinite on $[T,\infty]$.
On $[T,\lambda T]$ it is of the form~$C_\lambda (\nu T)^{3/2}$ 
 which conforms to the intuition that the heat kernel is around time $T$, mostly
 concentrated in a sphere of radius~$4\sqrt{\nu T}$.
\item
The velocity flow associated with an inviscid vortex line $\omega=
\delta_{x=0}\otimes \delta_{y=0}\otimes 1_{-1<z<1}$ behaves
as~$\frac{\pi}{\sqrt{x^2+y^2}}$ along $(0,0)\times[-1,1]$
hence belongs to $L^1([-1,1]^3)$ but is not square-integrable
on $[-1,1]^3$ ; one has therefore $\operatorname{Vol}(u_{\vert[-1,1]^3})=0$.
\item\label{OSEENEXAMPLE}
 The infinite viscous Oseen vortex line of direction $\mathbf{e}_3$ is the solution of \eqref{NS}
 given by $$u(t,x)=(u_h(t,x),0)\in\R^3$$ with
$$u_h(t,x)=\frac{1}{\sqrt{\nu t}} \:v \left(\frac{x_1}{\sqrt{\nu t}},\frac{x_2}{\sqrt{\nu t}}\right)\in\R^2 \quad\text{and}\quad
v(\xi)=\frac{\Gamma}{2\pi} \frac{\xi^\perp}{|\xi|^2}\left(1-e^{-|\xi|^2/4}\right),
\quad \Gamma\in\R.$$
The vorticity is the 2D heat kernel
$\omega(t) = (4\pi \nu t)^{-1} e^{-(x_1^2+x_2^2)/4\nu t}\,\mathbf{e}_3$.
Its characteristic scale is~$4\sqrt{\nu t}$.
As~$u$ is constant along the $z$-axis, it does not belong to any $L^p(\R^3)$.
One can however easily compute the volume function
in restriction to the cylinder $\Omega=\{(x,y,z) \,;\, x^2+y^2<1, \: |z|<1\}$.
The result is shown in Figure~\ref{figoseen}.

\smallskip
At the scale of $\Omega$, the vortex still appears concentrated around the $z$-axis
at $t=5\times 10^{-3}$. The peak of the volume function around $t=7\times10^{-2}$ occurs
when the characteristic scale of the vorticity matches that of $\Omega$.

\smallskip
Conversly, this simulation illustrates that for a given $t>0$, the length
$\lambda=\operatorname{Vol}(u:[t/2,t])^{1/3}$ (computed on a large enough domain)
determines the characteristic scale $4\sqrt{\nu t}$ of the vorticity.
Numerically, one has
$$\frac{\operatorname{Vol}(u:[t/2,t])}{\operatorname{Vol}(\Omega)}
\propto 1.17\left(\frac{4\sqrt{\nu t}}{\operatorname{Vol}(\Omega)^{1/3}}\right)^{1/4}$$
in the decades before the volume function reaches its peak (\textit{i.e.}~in physical
terms, when the characteristic scale of observation $\operatorname{Vol}(\Omega)^{1/3}$
exceeds the characteristic scale of the vorticity) .
\end{enumerate}

\subsection{Turbulence in the spectral sense of Obukhoff-Kolmogorov.}
\label{DEF}

This section provides the mathematical background
of Obukoff-Kolmogorov's
spectral theory of turbulence, known as ``K41 theory'' in reference to the publication
date of Kolmogorov's \cite{K41a}, \cite{K41b}, \cite{K41c} and Obukhoff's \cite{Obukov}
founding papers (see \cite[p.98]{Frish} for a precise chronology).

\subsubsection{K41-functions}

Let us start with an abstract mathematical definition.

\begin{dfn}
A function $u\in \mathcal{L}(\Omega) = L^\infty(\R_+;L^2(\Omega))\cap L^2(\R_+;\dot{H}^1(\Omega))$ is said to
be a \textsl{K41-function} on $[T_0,T_1]$ if there exists $C\in ]1,2[$ such that
\begin{equation}
\int_{\sigma(A^{1/2})}  K^2\,\bar{E}^\dagger(K)\, d\mu(K) \leq C
\int_{\Sigma(u;[T_0,T_1])}  K^2\,\bar{E}^\dagger(K)\, d\mu(K)
\end{equation}
where $\Sigma(u;[T_0,T_1])=\{K\in\sigma(A^{1/2}) \,;\, K\geq\operatorname{Vol}(u;[T_0,T_1])^{-1/3}\}$  and
$\operatorname{Vol}(u;[T_0,T_1])$ is defined by \eqref{VOL}.
\end{dfn}

On $\R^3$, if $u$ is a K41-function then for any $K_->0$ such that
\begin{subequations}
\begin{gather}
\label{HOMOGENEOUS}
(K_-)^{3}  \times \operatorname{Vol}(u;[T_0,T_1]) \leq 1
\\\intertext{one has for all $K_+>K_-$~:}
\label{ENSTROPHY}
\int_{0}^{\infty}  K^2\,\bar{E}^\ast(K)\, dK \leq (1+C) \int_{K_-}^{K_+}  K^2\,\bar{E}^\ast(K)\, dK.
\end{gather}
Any such interval $[K_-,K_+]$ is called an \textbf{inertial range} of $u$. The corresponding
\textbf{spectral Reynolds} number is~:
\begin{equation}\displaystyle \Re=\left(\frac{K_+}{K_-}\right)^{4/3}\end{equation}
One also defines a \textbf{spectral precision} parameter
(the reason of the fraction $5/3$ appears in the
next section)~:
\begin{equation}\label{CONSTANTS}
\gamma = \sup_{K\in[K_-,K_+]}\left| K \frac{d}{dK} (\log \bar{E}^\ast) +\frac{5}{3}\right|
\end{equation}
Obviously $\Re$ and $\gamma$ are dimensionless and there are infinitely many 
admissible quadruplets $(K_\pm,\Re,\gamma)$.
\end{subequations}

\medskip
On $\mathbb{T}^3$, the smallest non-vanishing frequency possible is $2\pi L^{-1}$ so \eqref{HOMOGENEOUS} is replaced by
\begin{subequations}
\begin{equation}\label{HOMOGENEOUSbis}
\begin{cases}
K_-=2\pi L^{-1} & \text{if}\quad
\displaystyle
\frac{\operatorname{Vol}(u;[T_0,T_1])}{\operatorname{Vol}(\mathbb{T}^3)}\geq (2\pi)^{-3},\\
2\pi L^{-1} \leq K_- \leq  \operatorname{Vol}(u;[T_0,T_1])^{-1/3} & \text{otherwise}.
\end{cases}
\end{equation}
Note that one always has $\operatorname{Vol}(u;[T_0,T_1])\leq \operatorname{Vol}(\mathbb{T}^3)$ by the Cauchy-Schwarz inequality.
The constant is $(2\pi)^{-3}\simeq 4\times 10^{-3}$.
In particular, one has~:
$$2\pi \leq L K_- \leq   \max\left\{2\pi;\left(\frac{\operatorname{Vol}(\mathbb{T}^3)}{\operatorname{Vol}(u;[T_0,T_1])}\right)^{1/3}\right\}.$$
Admissible frequencies $K_+$ are defined by~\eqref{ENSTROPHY} with an obvious
change of notations.
The discrete substitute for the definition of the spectral precision is~:
\begin{equation}\label{DISCRETEGAMMA}
\gamma = \max_{K_-\leq K_j<K_{j+1}\leq K_+} \left|\frac{\log\frac{\bar{E}^\dagger(K_{j+1})}{\bar{E}^\dagger(K_{j})}}{\log\frac{K_{j+1}}{K_j}}+\frac{5}{3}\right|
\end{equation}
where $\Sigma=\sigma(A^{1/2})=\{K_j\,;\,j\in\mathbb{N}\}$ with $K_j<K_{j+1}$.
\end{subequations}

\subsubsection{K41-turbulent flows}

Let us now turn back to fluid dynamics.
\begin{dfn}
Turbulence in the Kolomogorov-Obukov sense is the question to find
and describe solutions of~\eqref{NS} that are K41-functions on some
time interval $[T_0,T_1]$ and that possess at least one inertial range~$K_\pm$
in the asymptotic regime~:
\begin{equation}\label{K41}
\Re\gg1 \qquad\text{and}\qquad \gamma\log\Re \ll1.
\end{equation}
Such a solution  is  called a \textbf{K41-turbulent flow}.
\end{dfn}
In section \S\ref{RANGE}, the asymptotic~\eqref{K41} will be used to
recover the $K^{-5/3}$~law found in physics textbooks.

\medskip
\paragraph{Remarks}
\begin{enumerate}
\item
It is a mathematically open problem to construct exact
solutions of \eqref{NS} that possess this definite behavior.
However \eqref{K41} has been observed in numerous
experiments as the generic state of highly fluctuating flows (see \textit{e.g.}
\cite{BMONMF84}, \cite{KolmogorovConstant2},  \cite{BO97}, \cite{LA06}, \cite{KolmogorovConstant} and the references therein). Physics
textbooks usually don't mention the spectral precision $\gamma$
because they rely on Log-Log plots of the energy spectrum on which
the property $\gamma\log\Re\ll1$ is equivalent to having a substantial
amount of data concentrating along a straight line of slope $-5/3$
on $\log\Re$ decades of frequencies ; $\gamma$ is
the relative error on the slope of the line.
\item
In naive terms \eqref{ENSTROPHY} means that, by definition,
the K41-theory of turbulence is a spectral property of vortex structures because
at least half of the average enstrophy~$\langle{\norme[L^2]{\omega}^2}\rangle$
comes from the frequency range~$[K_-,K_+]$.
\item
The use of $\operatorname{Vol}(u;[T_0,T_1])$ at the foundation of K41-turbulence
is new. There is experimental evidence that K41-turbulence is generated on some
thin-structured subset of $\Omega$ deeply connected with the vorticity and
whose characteristic size $\ell_0=K_-^{-1}$ is
the largest scale involved in the inertial range.
Example~\ref{OSEENEXAMPLE} p.~\pageref{OSEENEXAMPLE}
has already shown a strong but subtle connection between the volume function and
the characteristic scale of vortex structures.
It will be shown (see Theorem~\ref{RANGETHM}) that
$$\ell_0\simeq\operatorname{Vol}(u;[T_0,T_1])^{1/3}$$
if $u$ is a K41-turbulent solution of~\eqref{NS}.
\end{enumerate}

\bigskip
Let us conclude this section with a short \textbf{dictionary} between mathematics and physics.
In physics textbooks, one can find statements like
\begin{equation}\label{DICTIONARY}
\ll\text{ Turbulent solutions $u$ of \eqref{NS} satisfy } F(u) \lesssim\, G(u)\gg
\end{equation}
where $F$ and $G$ are two functionals on $\mathcal{L}(\Omega)$
possibly depending on $T_0$ and $T_1$.
From a mathematical point of view, it
should be read as follows : there exists a function $C:\R^2\to]0,\infty[$ with
$$\displaystyle \lim_{\substack{\Re\to\infty\\\Re^\gamma\to1}} C(\Re,\gamma)
= C_0\in\:]0,\infty[$$
and constants $\Re_0,\gamma_0>0$
such that any solution of \eqref{NS} that is a K41-function on $[T_0,T_1]$ with a
Reynolds number $\Re\geq\Re_0$ and a spectral precision~$\gamma\leq\frac{\gamma_0}{\log\Re}$ satisfies
$$F(u)\leq C(\Re,\gamma) \,G(u).$$
In this case, any solution that admits
parameter~$(\Re,\gamma)$ in the asymptotic range \eqref{K41} will
 indeed satisfy $$F(u)\leq C_0' G(u) \qquad\text{with}\qquad C_0'/C_0\simeq1.$$
One should be aware that
that physics textbooks usually contain additional ``meta-assumptions'' such as
the isotropy or the homogeneity of the flow, which should then be translated
adequately and added to the assumptions of the mathematical statement.

\section{General properties of K41-functions}
\label{RANGE}

\subsection{Scale invariance and stability in Leray space}

The definition of K41-functions is invariant under \eqref{SCALING}~:
If $u$ is a K41-function on $[T_0,T_1]$ with parameters~$(K_\pm,\Re,\gamma)$
then $u_\lambda$~is a K41-function on $[\lambda^{-2} T_0, \lambda^{-2}T_1]$
with parameters $(\lambda K_\pm,\Re,\gamma)$.

\bigskip
The set of K41-functions is open in the Leray space $\mathcal{L}(\mathbb{T}^3)$. More
precisely, the following statement holds.
\begin{prop}
Assume that  $u\in \mathcal{L}(\mathbb{T}^3)$ is a
K41-function on $[T_0,T_1]$ with parameters~$(K_\pm,\Re)$. 
For any $\epsilon>0$ there exists $C_ \epsilon >0$
such that any $v\in \mathcal{L}(\mathbb{T}^3)$
with $$\norme[{L^2([T_0,T_1]\times\mathbb{T}^3)}]{u-v}^2
+\norme[{L^2([T_0,T_1]\times\mathbb{T}^3)}]{\nabla (u-v)}^2\leq C_\epsilon $$
is also a K41-function that satisfies \eqref{HOMOGENEOUS} and \eqref{ENSTROPHY}
with the same parameters $(K_\pm,\Re)$ but
with numerical constants~$1+\epsilon$ and~$1+C+\epsilon$.
\end{prop}
\Proof
Each term of \eqref{HOMOGENEOUS} and \eqref{ENSTROPHY} is
continuous on $H^1([T_0,T_1]\times\mathbb{T}^3)$.
\cqfd\endProof

Note however that the spectral precision $\gamma$ is not preserved, which
means that if a flow $u$ satisfies \eqref{K41}, any neighborhood
of $u$ in $\mathcal{L}(\Omega)$
will also contain functions that do not satisfy this asymptotic.
This illustrates a wise comment  by U.~Frish~\cite[pp.199-202]{Frish}:
``Questions (on turbulence) 
are likely to benefit from a close collaboration between mathematicians and
physicists but it will require more than better functional
analysis (...) ; some geometry is needed.'' 

\subsection{Kolmogorov's constant $\alpha$ and the $K^{-5/3}$~law}

Property \eqref{ENSTROPHY} involves the dissipation rate $\bar\varepsilon$
(which is the left-hand side) and frequencies.
There is only one way to define a quantity that has the same units as the
energy spectrum $[\text{length}]^3\cdot [\text{time}]^{-2}$ (\textit{i.e}~that
scales the same way under \eqref{SCALING})
and which is a power function of a dissipation rate
$\varepsilon=[\text{length}]^2\cdot [\text{time}]^{-3}$ and of a
frequency $K=[\text{length}]^{-1}$~:
$$\varepsilon(t)^{2/3} K^{-5/3}=[\text{length}]^3\cdot [\text{time}]^{-2}.$$
This fact makes the $\frac{5}{3}$ fraction in \eqref{CONSTANTS} a more obvious choice.
Note that even though $E^\ast(K,t)$ and $\varepsilon(t)$ depend on $\rho$, the
dimensionless fraction $$\frac{E^\ast(K,t)}{\varepsilon(t)^{2/3} K^{-5/3}}$$
is not ``missing'' a power $\rho^{1/3}$ because it is invariant under 
the transformation \eqref{SCALING}, which does not change the total mass.

\pagebreak[2]\bigskip
The following property of K41-functions is often mistaken for a definition of turbulence.
Its real meaning according to the dictionary \eqref{DICTIONARY} is that
the most valuable theorems concerning K41-functions will hold in the asymptotic
regime~\eqref{K41}.
\begin{prop}\label{K41b}
If $u$ is a K41-function on $[T_0,T_1]$,  the function $\alpha(K) =  \frac{\bar{E}^\ast(K)}{\bar{\varepsilon}^{2/3} K^{-5/3}}$ satisfies
\begin{equation}\label{K41bis}
\forall K,K'\in[K_-,K_+], \qquad
\Re^{-3\gamma/4} \leq \frac{\alpha(K)}{\alpha(K')} \leq
\Re^{3\gamma/4}
\end{equation}
\end{prop}
\Proof
The definition of $\gamma$ (namely \eqref{CONSTANTS} for $\R^3$ and
\eqref{DISCRETEGAMMA} for $\mathbb{T}^3$) implies~:
$$\left(\frac{K_1}{K_0}\right)^{-\frac{5}{3}-\gamma}
\leq \frac{\bar{E}^\ast (K_1)}{\bar{E}^\ast (K_0)} \leq
\left(\frac{K_1}{K_0}\right)^{-\frac{5}{3}+\gamma}$$
for any $K_0<K_1$ in $[K_-,K_+]$. Applying this inequality either to  $(K,K')=(K_0,K_1)$
or to $(K_1,K_0)$ one gets~:
$$\left(\frac{K_+}{K_-}\right)^{-\gamma}\leq
\left(\frac{K}{K'}\right)^{-\gamma}
\leq \frac{\alpha(K)}{\alpha(K')} \leq 
\left(\frac{K}{K'}\right)^{\gamma}\leq \left(\frac{K_+}{K_-}\right)^{\gamma}$$
hence the result.
\cqfd\endProof
\begin{dfn}
In the asymptotic regime~\eqref{K41}, the function $\alpha$ is essentially
constant on the inertial range and is called the \textbf{Kolmogorov constant}.
To fix further computations, one will choose from now on~:
\begin{equation}\label{ALPHA}
\alpha=\alpha(K_+)
\end{equation}
and \eqref{K41bis} reads $\displaystyle \left(\frac{\bar{E}^\ast(K)}{\alpha \bar{\varepsilon}^{2/3} K^{-5/3}}\right)^{\pm1}\leq \Re^{3\gamma/4}$ on $[K_-,K_+]$.
\end{dfn}

\paragraph{Remark} Proposition \ref{K41b} would also hold if in the definition of $\alpha(K)$
one would replace $\bar{\varepsilon}$ by another quantity having the dimension of
a dissipation rate, which in turn would change the value of the Kolmogorov constant~\eqref{ALPHA}. There is physical evidence that this definition leads to a universal numerical
value for $\alpha$ (see~\S\ref{UNIVERSALALPHA}) but mathematicians should question it.
The author is grateful to W.~Craig for this valuable remark.

\subsection{Bounds of the inertial range -- Expression of $K_\pm$}

Using dimensional analysis, there is only one way to define a frequency as a function
of $\bar\varepsilon$ and $\nu$~:
\begin{equation}\label{Kd}
K_d = \alpha^{-3/4} \left(\frac{\bar{\varepsilon}}{\nu^3}\right)^{1/4}.
\end{equation}
Since $\alpha$ is dimensionless, the power of $\alpha$ is arbitrary here, but it is the one
that provides the simplest statement in the following Theorem~\ref{RANGETHM}.
It also corresponds to the ``phantom'' homogeneity $\alpha\sim\rho^{1/3}$
mentioned above.
In physics, the length $\eta=K_d^{-1}$ is often referred to as
the ``Kolmogorov dissipation scale''.

Likewise, a frequency defined as a function of $\bar\varepsilon$ and $\bar{E}$ is~:
\begin{subequations}
\begin{equation}\label{Kc}
K_c= \alpha^{3/2} \frac{\bar\varepsilon}{\bar{E}^{3/2}}\cdotp
\end{equation}
In physics, the length $\ell_0 = K_c^{-1}$ is refered to as ``the size of large eddies''.
For $\mathbb{T}^3$, the energy spectrum misses the total impulsion (see \eqref{ENERGYbis}
and \S\ref{LOW}) so one substitutes the following definition for \eqref{Kc}~:
\begin{equation}\label{KcPERIODIC}
K_c= \alpha^{3/2}\bar\varepsilon \left(\bar{E} - \left[\int_{\mathbb{T}^3} \rho u_0(x)dx\right]^2\right)^{-3/2}
\end{equation}
\end{subequations}

Physics textbooks usually state that $K_+=K_d$ and $K_-=K_c$ by computing
$\bar{E}$ and $\bar{\varepsilon}$ for an idealized compactly supported energy
spectrum on $\R^3$~:
$$E^\ast(K) = \alpha \bar{\varepsilon}^{2/3} K^{-5/3} \mathbf{1}_{[K_-,K_+]}(K)$$
Converting this idea into a rigorous proof must be done carefully and actually
requires some additional assumptions to hold for $K_-$. Moreover, the computation
on $\mathbb{T}^3$ (which always seems to be dodged in physics textbooks)
happens to be extremely instructive (see also Remarks~\ref{RMKDISCRETESPECTRUM}
and~\ref{RMKDISCRETESPECTRUM2} p.~\pageref{RMKDISCRETESPECTRUM}).

\begin{thm}\label{RANGETHM}
The following inequalities hold.
\begin{description}
\item[1. Case of $\Omega=\R^3$.]
\begin{subequations}
If $u$ is a K41-function on $\Omega=\R^3$ with parameters $(K_\pm,\Re,\gamma)$, then~:
\begin{align}\label{K+}
\left(\frac{ \Re^{-3\gamma/4}}{6(1-\Re^{-1})}\right)^{3/4} &\leq \frac{K_+}{K_d} \leq \left(\frac{4\Re^{3\gamma/4}}{3(1-\Re^{-1})}\right)^{3/4},
\\
\label{K-}
\left(\frac{3}{2}\Re^{-3\gamma/4}(1-\Re^{-1/2})\right)^{3/2} &\leq
\frac{K_-}{K_c} \cdotp
\end{align}
Moreover, if $u$ is a solution of \eqref{NS} with initial data $u_0\in L^1\cap L^2(\R^3)$ then~:
 \begin{equation}
 \frac{K_-}{K_c}\leq \left(\frac{9\pi^2}{3\pi^2-4}\times\Re^{3\gamma/4}(1-\Re^{-1/2})\right)^{3/2}
 \end{equation}
and
\begin{equation}\label{K-Vol}
\frac{(3\pi^2-4)\Re^{-3\gamma/2}}{36(1-\Re^{-1/2})}
\leq (K_-)^{3}\times\operatorname{Vol}(u;[T_0,T_1]) \leq1.
\end{equation}
\end{subequations}
\item[2. Case of $\Omega=\mathbb{T}^3$.]
\begin{subequations}
If $u$ is a K41-function on $\Omega=\mathbb{T}^3$
with parameters $(K_\pm,\Re,\gamma)$, then~:
\begin{align}\label{K+bis}
\left(\frac{\Re^{-3\gamma/4}}{6}\right)^{3/4} &\leq \frac{K_+}{K_d}
\leq \left(\frac{\frac{1}{2}\,\Re^{3\gamma/4}
+ K_d^{-4/3}\times \underset{{\Re\to\infty}}{O(1)}}
{\frac{15}{16}  - \frac{3}{2}\Re^{-1}}\right)^{3/4},\\
\left( \Re^{-3\gamma/4} 
\left(\frac{15}{16}-3\Re^{-1/2}\right) \right)^{3/2}
&\leq \frac{K_-}{K_c}
\qquad \text{provided }K_-> 3(2\pi L^{-1}).
\end{align}
\end{subequations}
Moreover, if $u$ is a solution of \eqref{NS} with $\displaystyle \int_{\mathbb{T}^3} \rho u_0(x)dx=0$ and
\begin{equation}\label{NBTHASSUMPTION}
C(n_-)= \frac{\operatorname{Card}\{z\in\mathbb{Z}^3\,;\, |z|^2=n_-\}}{8\pi^3\sqrt{n_-}}<1
\quad\text{where}\quad n_-=\left(\frac{L K_-}{2\pi}\right)^2,
\end{equation}
then~:
\begin{subequations}
 \begin{equation}
 \frac{K_-}{K_c}\leq  \left(\frac{12\Re^{3\gamma/4}}{1-C(n_-)}\right)^{3/2}
 \end{equation}
and
\begin{equation}
\frac{1-C(n_-)}{12\Re^{3\gamma/2}}
\leq (K_-)^{3}\times\operatorname{Vol}(u;[T_0,T_1]) \leq
\max\left\{1;(2\pi)^3\:\frac{\operatorname{Vol}(u;[T_0,T_1])}
{\operatorname{Vol}(\mathbb{T}^3)}\right\}.
\end{equation}
\end{subequations}
\end{description}
\end{thm}

\begin{cor}
In the turbulent asymptotic $\Re\to\infty$ and $\Re^\gamma\to 1$, one has
\begin{equation}
C\leq \frac{K_c}{\operatorname{Vol}(u;[T_0,T_1])^{-1/3}} \leq C'
\end{equation}
for two numerical constants $C,C'$.
\end{cor}

\begin{figure}[hb!]\sf
\begin{minipage}{\textwidth}
\rule{\textwidth}{.1pt}
\begin{multicols}{2}
\begin{center}
\includegraphics[width=\columnwidth]{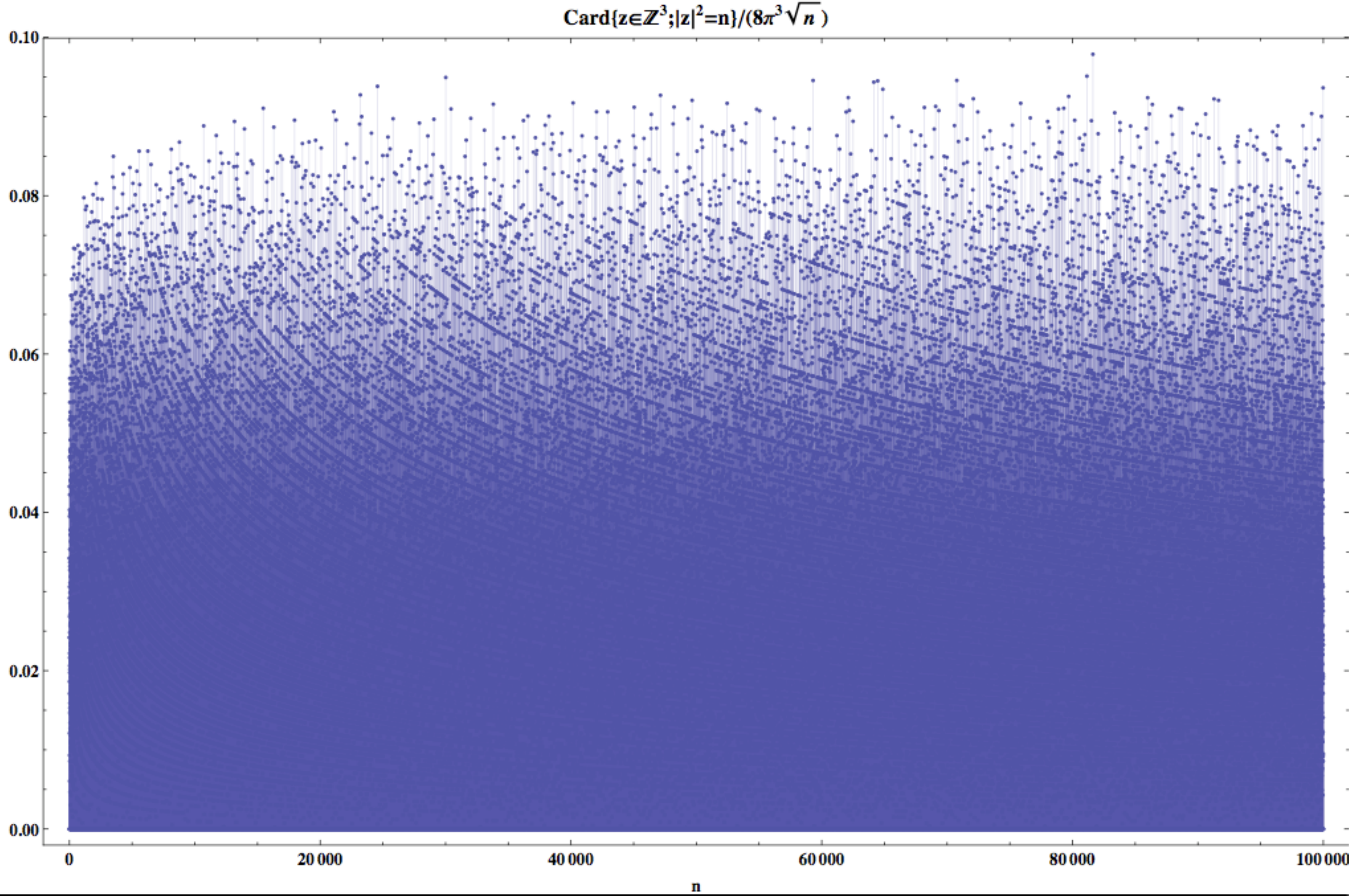}
\end{center}
\caption{\label{SUMOFSQUARES} Computation of $C(n)$.}\vspace*{1ex}
Numerical test that $r_3(n)/(8\pi^3\sqrt{n})<0.1$
for $n\leq 10^5$ where $r_3(n)=\operatorname{Card}\{z\in\mathbb{Z}^3\,;\, |z|^2=n\}$.
A similar test for $10^{10}\leq n \leq 10^{10}+100$ checks the same numerical bound
 $C(n)<0.1$ with a maximum of $9.37\times 10^{-2}$ for $n=10^{10}+1$.
It is obtained using \texttt{SquaresR[3, n]} with  \textit{Mathematica}\copyright.
The problem of computing the number of representations of an integer as the sum of
three squares has been addressed historically \textit{e.g.} in  \cite{Estermann} and \cite{bateman} but the asymptotic behavior of $r_3(n)$ cannot be read directly on
Hardy's explicit formula.
The series $\sum_{n\leq N} r_3(n)= \frac{4\pi}{3}N^{3/2}+O(N^{\epsilon +29/44})$
is the number of lattice points inside the sphere of radius~$\sqrt{N}$ (see \cite{ChamizoIwaniec}).
Similarly, a recent paper \cite{choi} shows that $\sum_{n\leq N} r_3(n)^2 = \frac{8\pi^4}{21\zeta(3)} N^2+O(N^{14/9})$ but again the remainder is
too large to prove the asymptotic behavior $r_3(n)=O(\sqrt{n})$.
Note that this conjectured asymptotic is extremely sharp because the sequence
$r_3(n)$ vanishes on the subsequence~$r_3(4^p(8q+7))=0$
thus cannot have a power law equivalent.
The best known estimate \cite{ChamizoIwaniec}
is $r_3(n)/\sqrt{n}= O(n^{\epsilon})$
for any $\epsilon>0$.
Thus it might be that $r_3(n)/\sqrt{n}$ is not bounded but grows extremely slowly.
However, from our numerical data, one can infer that even if it diverges as $\log n$,
then $r_3(n)$ will finally exceed $8\pi^3\sqrt{n}$ only for $n\sim 10^{50}$.
Such a scale is unrealistic because it would require the building of
a periodic torus of at least $L=6 \text{ light years}$ in order to ensure that
such a large wave number would investigate scales that exceed
the atomic one $(2\pi L^{-1}\sqrt{n})^{-1}=10^{-9} \text{m}$.
Therefore, in the range of validity of Navier-Stokes equations,
assumption \eqref{NBTHASSUMPTION} is numerically satisfied.
\end{multicols}
\end{minipage}
\end{figure}

Assumption $C(n_-)<1$ is still an open problem in number theory ;
the numerical test presented in Figure~\ref{SUMOFSQUARES}
ensures it is satisfied for most if not all practical purposes.

Numerically, the theorem reads for $\Omega=\R^3$~:
$$0.260  K_d\leq K_+\leq 1.25 K_d, \qquad 1.83 K_c \leq K_- \leq 6.46 K_c
\quad\text{and}\quad 0.711\leq (K_-)^{3}\times\operatorname{Vol}(u;[T_0,T_1])  \leq 1.$$

For $\Omega=\mathbb{T}^3$, the asymptotic \eqref{K41} implies $K_+\to\infty$
because $K_-\geq 2\pi L^{-1}$, thus \eqref{K+bis} provides $K_d\to\infty$.
One has $C(n_-)<0.1$ for $n_-\leq 10^5$ and
$$0.260  K_d\leq K_+\leq 0.625 K_d, \qquad 0.907 K_c \leq K_- \leq 48.7 K_c$$
and
$$8.33\times 10^{-3} \leq (K_-)^{3}\times\operatorname{Vol}(u;[T_0,T_1])  \leq 248.1.$$

\paragraph{Remark}
One could object that the numerical constants of the previous statement are not fundamental.
Indeed, for most of what will follow, one could just write $C(\Re,\gamma)$
with $C(\Re,\gamma)\to C>0$ in the asymptotic~\eqref{K41}. However, for Theorem~\ref{SMOOTHTURBULENCE}, which investigates the subtle relation between
smoothness and intermittency, the numerical values of the constants will mark
the difference between an empty statement and a meaningful result, which means that
one will have to show some discipline in each intermediary result.

\bigskip

\Proof[ (estimate of $K_+$)]
One can compute $\bar{\varepsilon}$ using \eqref{EPSILON}.
On the inertial range \eqref{K41bis}  leads to~:
$$ \bar\varepsilon \geq \nu \int_{K_-}^{K_+}  K^2\,\bar{E}^\ast(K)\, dK\geq
\frac{3\alpha\nu\bar{\varepsilon}^{2/3}}{4\Re^{3\gamma/4}} (K_+^{4/3}-K_-^{4/3})=
\frac{3\alpha\nu\bar{\varepsilon}^{2/3}}{4\Re^{3\gamma/4}} (1-\Re^{-1})K_+^{4/3}$$
hence the right-hand side of \eqref{K+}. 
Conversely, according to \eqref{ENSTROPHY} one can compute $\bar{\varepsilon}$
using only the inertial range where \eqref{K41bis} again provides~:
$$\bar\varepsilon \leq 4\nu \int_0^\infty  K^2\,\bar{E}^\ast(K)\, dK\leq 
8\alpha\nu\Re^{3\gamma/4}\bar{\varepsilon}^{2/3}\int_{K_-}^{K_+}  K^{2-5/3}\, dK =
6\alpha\nu\Re^{3\gamma/4}\bar{\varepsilon}^{2/3}(1-\Re^{-1})K_+^{4/3}$$
hence  the left-hand side. On $\mathbb{T}^3$, one gets instead~:
$$\frac{1}{2}\bar{\varepsilon}
\: \leq \:2\nu \sum_{K\in \Sigma\cap [K_-,K_+]} \left(\frac{2\pi}{K L}\right) K^2\,\bar{E}^\dagger(K)
\cdot (2\pi L^{-1}) \: \leq \: \bar{\varepsilon}$$
hence
$$\Re^{-3\gamma/4}\:\frac{\bar{\varepsilon}^{1/3}}{4\alpha\nu} \leq
\left(\frac{2\pi}{L}\right)^2\sum_{K\in \Sigma\cap [K_-,K_+]}K^{-2/3}
\leq \frac{\bar{\varepsilon}^{1/3}}{2\alpha\nu}\:\Re^{3\gamma/4}.$$
The next step is to show that the center term
is equivalent to $K_+^{4/3}$ (the apparently different game of powers reflects
the spectral asymptotic \eqref{SPECTRALASSYMPTOTIC} of the Stokes operator ; note
that physical dimensions are the same).
One has~:
$$\left(\frac{2\pi}{L}\right)^2\sum_{K\in \Sigma\cap [K_-,K_+]}K^{-2/3}
= \left(\frac{2\pi}{L}\right)^{4/3}\sum_{n\in\Box_3 \cap [n_-,n_+]}n^{-1/3}$$
where $\Box_3=\mathbb{N}\backslash\{4^p(8q+7) \,;\, p,q\in\mathbb{N}\}$ denotes the set
of integers that are the sum of three squares. Its complimentary 
is included in the subset of integers $n\equiv0,4 \text{ or }7 \text{ mod }8$, therefore~:
$$\sum_{n=1}^N n^{-1/3}\geq \sum_{n\in\Box_3 \cap [1,N]}n^{-1/3}\geq\sum_{n=1}^N n^{-1/3}
- \sum_{j\in\{0,4,7\}}\sum_{8k+j\leq N} (8k+j)^{-1/3}$$
which after comparison to an integral boils down to
$$\frac{3 N^{2/3}}{2} \geq \sum_{n\in\Box_3 \cap [1,N]}n^{-1/3} \geq
\frac{3}{2}\left((1+N)^{2/3}-1\right)
-\frac{3}{8^{1/3}} \left(\frac{3 (N/8)^{2/3}}{2}\right)
= \frac{15}{16} N^{2/3}+\underset{N\to\infty}{O(1)}.$$
One uses those estimates to compute the sum for $K\in\Sigma\cap[K_-,K_+]$~:
$$ \frac{3}{2} K_+^{4/3} \geq
\left(\frac{2\pi}{L}\right)^2\sum_{K\in \Sigma\cap [K_-,K_+]}K^{-2/3}
\geq K_+^{4/3} \left(\frac{15}{16}  - \frac{3}{2}\Re^{-1}\right) + \underset{{K_+\to\infty}}{O(1)} $$
Note that in the asymptotic \eqref{K41}, one has $K_+\to\infty$ because $K_-\geq 2\pi L^{-1}$.
One finally gets~:
$$\Re^{-3\gamma/4}\:\frac{\bar{\varepsilon}^{1/3}}{4\alpha\nu} \leq
\frac{3}{2} K_+^{4/3} \qquad\text{and}\qquad
K_+^{4/3} \left(\frac{15}{16}  - \frac{3}{2}\Re^{-1}\right) + \underset{{\Re\to\infty}}{O(1)}\leq \frac{\bar{\varepsilon}^{1/3}}{2\alpha\nu}\:\Re^{3\gamma/4}.$$

\cqfd\endProof
\Proof[ (estimate of $K_-$)]
To estimate $K_-$, one  compares the total energy with the
energy contained in the inertial range. 
If $u$ is a K41-function on $\R^3$, one has~:
$$\int_{K_+}^\infty E^\ast(K,t)dK \leq K_+^{-2} \int_{K_+}^\infty K^2 E^\ast(K,t)dK 
\leq K_+^{-2} \int_{K_-}^{K_+} K^2 E^\ast(K,t) dK \leq \int_{K_-}^{K_+} E^\ast(K,t)dK$$
and  \eqref{ENERGY} provides
\begin{subequations}
\begin{equation}\label{Epm}
\int_{K_-}^{K_+} \bar{E}^\ast(K) dK \leq
\bar{E} \leq 2 \int_{K_-}^{K_+} \bar{E}^\ast(K) dK + \int_0^{K_-}  \bar{E}^\ast(K) dK.
\end{equation}
On $\mathbb{T}^3$  the corresponding inequalities are~:
\begin{equation}\label{Epm2}
\begin{gathered}
\bar{E} - \left(\int_{\mathbb{T}^3}\rho u_0(x) dx\right)^2 \geq \left(\frac{2\pi}{L}\right) 
\sum_{K\in\Sigma\cap[K_-,K_+]} \left(\frac{2\pi}{KL}\right) \bar{E}^\ast(K)
\\
\bar{E} - \left(\int_{\mathbb{T}^3}\rho u_0(x) dx\right)^2 \leq
\left(\frac{2\pi}{L}\right) \left\{
\sum_{\substack{K\in\Sigma,\\ K< K_-}} \left(\frac{2\pi}{KL}\right)\bar{E}^\ast(K) +
2 \sum_{K\in\Sigma\cap[K_-,K_+]} \left(\frac{2\pi}{KL}\right)
\bar{E}^\ast(K) \right\}.
\end{gathered}
\end{equation}
\end{subequations}
Using \eqref{K41bis} on $[K_-,K_+]$, the lower bound of $\bar{E}$ provides \eqref{K-}~:
$$\bar{E} \geq \frac{\alpha \bar{\varepsilon}^{2/3}}{\Re^{3\gamma/4}}\int_{K_-}^{K_+}  K^{-5/3} dK= \frac{3\alpha \bar{\varepsilon}^{2/3}}{2\Re^{3\gamma/4}} (1-\Re^{-1/2}) K_-^{-2/3}.$$
The computation is similar on $\mathbb{T}^3$ ; one gets~:
$$ \left(\frac{K_-}{K_c}\right)^{2/3} \geq \Re^{-3\gamma/4} \times n_-^{1/3}C(n_-,n_+)
\qquad\text{with}\qquad C(n_1,n_2)=\sum_{n\in\Box_3\cap[n_1,n_2]} n^{-4/3}.$$
This time, one uses $C(n_-,n_+)=C(n_-,\infty)-C(1+n_+,\infty)$ with
$$3 (N-1)^{-1/3}\geq C(N,\infty) \geq \frac{15}{8}N^{-1/3}+O(N^{-4/3})
\geq \frac{15}{16}N^{-1/3},$$
the last inequality being valid if $N\geq10$.
The proof of the converse inequality and the last statement about $\operatorname{Vol}(u;[T_0,T_1])$ will be postponed until Proposition~\ref{K-PRECISED}. One will use the additional
assumptions and~\eqref{HOMOGENEOUS}--\eqref{HOMOGENEOUSbis}
to estimate the spectrum for low-frequencies in \eqref{Epm}--\eqref{Epm2}.
\cqfd\endProof

\section{Low-frequency spectrum}
\label{LOW}

In physics, two different low-frequency spectra have been described
(see \textit{e.g.} \cite[p.91 and p.358]{Davidson})~:
the Saffman spectrum behaves like $K^2$ and the Batchelor spectrum behaves
like $K^4$ as $K\to0$. Numerical simulations on $\mathbb{T}^3$ have confirmed that both
can exist and that the choice between $K^2$ or $K^4$ depends on the initial data~:
the initial behavior of the spectrum is preserved, though in case of $K^4$ the spectrum
does not seem totally stable in time (see \cite[fig.~5a and fig.~7]{lesieur00} and the references
therein where this subtle instability is referred to as backscatter).
For real-world experiments, the question is not properly
answered because the first Fourier coefficients reflect
the large scales in the experimental protocol more than turbulence itself.

\medskip\noindent%
In this section are proved the following points~:
\begin{itemize}
\item Theorem \ref{SPECTRUMFROML1}  : Saffman's estimate holds in general.
\item Theorem~\ref{SPECTRUMFROMLOCALIZATION} : On $\R^3$, the spatial localization of the initial data determines spectra in~$K^{2(1+\beta)}$ with possibly any~$\beta<1$ depending
on the exact decay at infinity.
\item
Batchelor's $K^4$ spectrum ($\beta =1$)
seems to occur only for unstable highly localized solutions. Furthermore, it was not
possible to obtain
any nontrivial estimate with $\beta>1$ even for highly localized flows.
\end{itemize}
Note that on $\mathbb{T}^3$, one can always artificially improve
the power of $K\in\Sigma^\ast$
of upper bounds by multiplying by~$K L/(2\pi)\geq1$, thus it is worth
pointing out that~\eqref{EbarLOW} below does not contain the length $L$ on the right-hand side.

\subsection{Estimate from $u_0\in L^1(\Omega)$}

From a physical point of view, it is reasonable to assume that
\begin{equation}\label{INT0}
\int_{\Omega} \rho u(t,x)dx=0
\end{equation}
which means that the fluid is globally at rest in the given coordinate system.
On $\Omega=\mathbb{T}^3$, Cauchy-Schwarz ensures that $u(t)\in L^1(\Omega)$
for any $t\geq0$ ; moreover, the total impulsion is constant in time so  \eqref{INT0} holds
for any $t\geq0$ if and only if it holds for the initial data.
On $\R^3$, the following statement justifies why \eqref{INT0} holds anyway.
\begin{thm}\label{PROPLOCALIZATION}
If $u$ is a Leray solution of \eqref{NS} with initial data $u_0\in L^1\cap L^2(\R^3)$, then
$u(t)\in L^1(\Omega)$ for any $t\geq0$ and~\eqref{INT0} holds for (at least) almost
every $t>0$.
\end{thm}
\Proof
On $\R^3$, the following inequality holds for Leray solutions of \eqref{NS}~:
\begin{equation}
\label{BOUND-L1}
\exists C>0, \quad \forall t\geq0,\qquad
\norme[L^1(\R^3)]{u(t)} \leq \norme[L^1(\R^3)]{u_0} + C\sqrt{\frac{t}{\nu}}\norme[L^2(\R^3)]{u_0}^2
\end{equation}
which ensures that $u(t)\in L^1(\R^3)$. Let us conclude first before proceeding
to the proof of~\eqref{BOUND-L1}.
Leray solutions are known to be smooth for almost every $t>0$ (see \cite{L34}).
At any such time, the following integration by part holds for any smooth bounded
subset $\Omega\subset\R^3$~:
$$\int_\Omega u_i=-\int_\Omega x_i(\div u) dx+\int_{\partial\Omega} x_i \:u \cdot dn.$$
For such times, the fact that $u(t)\in L^1\cap C^0(\R^3)$ also implies that~:
$$\lim_{R\to\infty} \int_R^\infty \left| \int_{|x|=r} x_i u(x) \cdot dn(x)\right| \frac{dr}{r} = 0.$$
One can choose a sequence of radii $r_k\to\infty$ such that
$\displaystyle \lim_{k\to\infty}\int_{|x|=r_k} x_i \:u(x) \cdot dn(x)=0$.
Combining both identities and $\div u = 0$ gives \eqref{INT0} for
any time $t\geq0$ such that $u(t)\in H^1\cap C^0$.

\medskip
Inequality \eqref{BOUND-L1} is part of the folklore of fluid mechanics and,
contrary to most estimates, it does not involve an exponential growth in time.
Let us recall briefly its derivation.
Rewriting \eqref{NS} as a heat-equation with a non-linear source term (see \cite[chap.~11]{Lemarie}), one gets~:
$$u(t,x)=e^{\nu t \Delta}u_0(x) + \int_0^t\int_{\R^3} (u\otimes u)(t',x') K(t-t',x-x') dx'dt'$$
with a convolution kernel that satisfies $|K(t,x)|\leq C(|x|+\sqrt{\nu t})^{-4}$ (see \textit{e.g.}~\cite{Vigneron}). 
Therefore, for any $\tau\in[0,t]$, one has~:
$$\norme[L^1]{u(\tau)} \leq \norme[L^1]{u_0} + C\norme[{L^\infty([0,t];L^1)}]{u\otimes u}
\int_0^t\int_0^\infty \frac{r^2 dr dt}{(r+\sqrt{\nu t})^4}.$$
But $\norme[L^1]{u\otimes u}\leq \norme[L^2]{u_0}^2$ and the last integral computes
down to $2\sqrt{t/\nu}$, hence \eqref{BOUND-L1}.
\cqfd\endProof

\bigskip
One can now deduce bounds on the lower end of the energy spectrum.
\begin{thm}\label{SPECTRUMFROML1}
If $u$ is a Leray solution of \eqref{NS} with $u_0\in L^2\cap L^1(\Omega)$
on~$\Omega=\R^3$ or $\mathbb{T}^3$, then for any $T_0<T_1$,
the energy spectrum on $[T_0,T_1]$ satisfies~:
\begin{equation}\label{EbarLOW}
 \begin{cases} \displaystyle
\bar{E}^\ast(K)\leq \operatorname{Vol}(u;[T_0,T_1]) K^2 \bar{E}
\times 
4/\pi^2  & \text{on }\R^3,\\[1ex] 
\bar{E}^\dagger(K)\leq \operatorname{Vol}(u;[T_0,T_1]) K^2 \bar{E}
\times 
\frac{1}{8\pi^3\sqrt{n}} \operatorname{Card}\{z\in\mathbb{Z}^3\,;\, |z|^2=n\}
& \text{on }\mathbb{T}^3,
\end{cases}
\end{equation}
where the volume $\operatorname{Vol}(u;[T_0,T_1])$ is defined by \eqref{VOL}.
For $\mathbb{T}^3$, one has $L=\operatorname{Vol}(\mathbb{T}^3)^{1/3}$
and $K=2\pi L^{-1}\sqrt{n}$ with  an integer $n\in\Box_3$ that is the sum of three squares.
\end{thm}
In the case of $\Omega=\mathbb{T}^3$, one does not require \eqref{INT0} to hold.
The numerical factor is illustrated by Figure~\ref{SUMOFSQUARES}.

\Proof 
Using $\norme[L^\infty]{\hat{u}}\leq\norme[L^1]{u}$,
Proposition~\ref{EXPLICITSPECTRUM} reads on $\R^3$~:
$$E^\ast(K,t)\leq  \rho \norme[L^1]{u(t)}^2 K^2 \times 
\frac{1}{2\pi^2}\int_0^\infty\psi(r)  r^2 dr $$
The integral is bounded by $8$
because $\supp\psi\subset[2^{-1},2]$ and $\int \psi(r)dr/r=1$. 
On $\mathbb{T}^3$, one has instead~:
$$E^\dagger(K,t)\leq  \rho \norme[L^1]{u(t)}^2 K^2 \times
\frac{1}{4\pi^2 LK} \operatorname{Card}\{k\in\hat{\mathbb{T}}^3\,;\, |k|=K\}.$$
Conclusion \eqref{EbarLOW} then follows immediately from
the definition \eqref{VOL} of the volume function.
\cqfd\endProof

\subsection{Estimates from localization norms of $u$ on $\R^3$}

One can improve the low-frequency estimation of the spectrum using localization
properties of the flow.
\begin{thm}\label{SPECTRUMFROMLOCALIZATION}
Given $\beta\in]0,1[$, there exists a constant $C_\beta>0$ such that
any Leray solution of \eqref{NS} on $\Omega=\R^3$ with initial data in $L^2\cap L^1$
satisfies for almost every $t\geq0$~:
\begin{equation}\label{LOW2}
E^\ast(t,K) \leq \rho C_\beta K^{2(1+\beta)}\norme[L^2(\R^3)]{(1+|x|)^{3/2+\beta} u(t)}^2.
\end{equation}
If the initial data $u_0$ satisfies \textit{e.g.}
\begin{equation}\label{LOCALIZ}
(1+|x|)^{3+\beta'} u_0 \in L^\infty(\R^3) \qquad\text{with}\qquad \beta'>\beta
\end{equation}
then the right-hand side of \eqref{LOW2} is finite  as long as $u$
is a smooth solution of $\eqref{NS}$ on $[0,t]$.
\end{thm}
\paragraph{Remarks.}
\begin{enumerate}
\item For $\beta=1$, the estimate \eqref{LOW2} holds with an extra multiplicative
factor $(1-\log K/K_0)$ for $K<K_0$ but the right-hand side
is infinite unless a non-generic (necessary but not sufficient) condition holds~:
\begin{equation}\label{DIAGONALENERGY}
\int_{\R^3} u_i(t,x)u_j(t,x) dx = \frac{1}{3}\norme[L^2]{u(t)}^2 \delta_{i,j}.
\end{equation}
This condition is not invariant by the flow of \eqref{NS}
and \cite{BRANDOLESE} contains examples of flows that will
check and violate \eqref{DIAGONALENERGY} at prescribed times.
Those smooth flows will satisfy \eqref{LOW2} for any $\beta<1$ and $t\geq0$ but
they will satisfy it for~$\beta=1$ (with the logarithmic correction) if and only if
$t\in\{t_0,t_1,\ldots\}$.
Examples of flows that satisfy \eqref{LOW2}  for all time with $\beta=1$
(with again the logarithmic correction) can be found in \cite{BRANDOLESE2}. 
\item
Contrary to~\eqref{BOUND-L1},
the best known short-time bound for weighted norms
is exponential in $t$.
Other localization norms of $u_0$ could be used instead of \eqref{LOCALIZ}.
For example the same theorem holds if~:
$$(1+|x|)^{3\left(1-\frac{1}{p}\right)+\beta'}u_0(x)\in L^p(\R^3)$$
with $p>3$ and $\beta'>\beta$ (see \cite{Vigneron}).
\item
On $\Omega=\mathbb{T}^3$, weighted norms  are
meaningless ; however, provided \eqref{INT0} is satisfied, one could get a similar result with
the weighted norms replaced by~:
$$\norme[C^\beta]{\hat{u}(t)} = \sup_{k\neq k'} \frac{|\hat{u}(t,k)-\hat{u}(t,k')|}
{|k-k'|^\beta}\cdotp$$
To the best of my knowledge, propagation of this semi-norm by the flow on $\mathbb{T}^3$ for $\beta\in]0,1[$ has not yet been studied. Neither has the propagation
of $u(t)\in\mathcal{F}^{-1}(C^\beta(\R^3))$ in the continuous case.
One  can however expect a generic failure of  the propagation of $\sup_{k\neq0}|k|^{-1}|\hat{u}(t,k)|$ corresponding to $\beta=1$.
\end{enumerate}

\Proof
One has only to prove \eqref{LOW2}. The rest of the statement follows from \cite{Vigneron}
for $\beta\in]0,1[$.
One has the following chain of continuous inclusions if $\beta\notin\mathbb{N}$~:
$$\norme[C^\beta]{\hat{u}}
\leq \norme[H^{\frac{3}{2}+\beta}]{\hat{u}} \leq \norme[L^2]{(1+|x|)^{3/2+\beta} u(t)}.$$
If $\beta\in ]0,1[$, one gets $|\hat{u}(t,\xi)| \leq C |\xi|^\beta \norme[L^2]{(1+|x|)^{3/2+\beta} u(t)}$ when $\hat{u}(t,0)=\int_{\R^3}u(t,x)dx=0$, which (according to Theorem~\ref{PROPLOCALIZATION}) holds for almost every $t\geq0$. Therefore~:
\begin{align*}
E^\ast(K,t) &= (2\pi)^{-3}\rho K^2\int_0^\infty\int_{\mathbb{S}^2}\psi(r)
|\hat{u}(t,Kr\vartheta)|^2 r^2 dr d\vartheta\\
& \leq C \rho K^{2(1+\beta)} \left(\int_0^\infty\psi(r)
 r^{2(1+\beta)} dr\right)  \norme[L^2]{(1+|x|)^{3/2+\beta} u(t)}^2
\end{align*}
which gives \eqref{LOW2} with a constant independent of $\psi$ because
$\supp\psi\subset[2^{-1},2]$ and
$\int \psi(r)dr/r=1$. 

\medskip
Let us briefly justify the notes that follow the Theorem.
The precise relation between~\eqref{DIAGONALENERGY} and the 
finiteness of the right-hand side of \eqref{LOW2}
is extensively studied in~\cite{BrandoleseVigneron}.
When $\beta=1$, the Sobolev space $H^{5/2}(\R^3)$ is not included in~$\text{Lip}(\R^3)$
but in Calderon's
space $C^1_\ast(\R^3)$ for which
$|f(\xi)-f(\eta)|\leq C |\xi-\eta|(1-\log(|\xi-\eta|))$ when $|\xi-\eta|<1$ ;
see \textit{e.g.}~\cite[p. 31]{Chemin}. When $\beta>1$, \textit{e.g.} for $\beta\in]1,2[$, the
following estimate holds~:
$$|\hat{u}(t,\xi)|=|\hat{u}(t,\xi)-\hat{u}(t,0)|\leq|\xi|
\left(\int_0^1 |\nabla \hat{u}(t,\sigma\xi)| d\sigma\right)
\leq |\xi| \left(|\xi|^{\beta-1}\norme[C^\beta]{\hat{u}(t)}\int_0^1 \sigma^{\beta-1}d\sigma+\norme[L^1]{|x| u(t)}\right) $$
so if $u$ is a highly localized flow like those  in \cite{BRANDOLESE2}, 
then \eqref{LOW2} holds with $\beta>1$ ; however,~$K^{2(\beta+1)}$ must be
replaced by~$K^4$ because $|\hat{u}(t,\xi)|\neq o(|\xi|)$, even for highly localized flows.
\cqfd\endProof

\subsection{Precisions on $K_-$}

Low-frequency bounds on the energy spectrum allow us  to 
conclude the proof of Theorem~\ref{RANGETHM}.
\begin{prop}\label{K-PRECISED}
Let $u$ be a Leray solution of \eqref{NS} and a K41-function on $[T_0,T_1]\times\Omega$.
\begin{enumerate}
\item
If $\Omega=\R^3$ and $u_0\in L^1\cap L^2(\R^3)$,  one has
\begin{equation}\label{K-2}
\frac{K_-}{K_c}\leq \left(\frac{9\pi^2}{3\pi^2-4}\times\Re^{3\gamma/4}(1-\Re^{-1/2})\right)^{3/2}
\end{equation}
and 
\begin{equation}\label{K-2bis}
\frac{3\pi^2-4}{36\Re^{3\gamma/2}(1-\Re^{-1/2})}
\leq (K_-)^{3}\times\operatorname{Vol}(u;[T_0,T_1])  \leq 1.
\end{equation}
\item
If $\Omega=\mathbb{T}^3$, let us assume that $\displaystyle \int_{\mathbb{T}^3} \rho u_0(x)dx=0$ and
$$
C(n_-)= \frac{\operatorname{Card}\{z\in\mathbb{Z}^3\,;\, |z|^2=n_-\}}{8\pi^3\sqrt{n_-}}<1
\quad\text{where}\quad n_-=\left(\frac{L K_-}{2\pi}\right)^2.
$$
Then similar inequalities to \eqref{K-2} and \eqref{K-2bis} hold, namely
\begin{equation}\label{PRECISEDPERIODIC}
\frac{K_-}{K_c}\leq \left(\frac{12\Re^{3\gamma/4}}{1-C(n_-)}\right)^{3/2}
\end{equation}
and
\begin{equation}
\frac{1-C(n_-)}{12\Re^{3\gamma/2}} \leq (K_-)^{3}\times\operatorname{Vol}(u;[T_0,T_1])  \leq 
\max\left\{1;(2\pi)^3\:\frac{\operatorname{Vol}(u;[T_0,T_1])}
{\operatorname{Vol}(\mathbb{T}^3)}\right\}.
\end{equation}
\end{enumerate}
\end{prop}
Assumption $C(n_-)<1$ 
is still an open problem in number theory
but the systematic numerical test presented in Figure~\ref{SUMOFSQUARES}
ensures that it is satisfied for at least $n_-\leq 10^5$, which should be sufficient for
most practical purposes. For example, in a periodic domain of size $L=1$m, it means that
the assumption is satisfied at least as soon as the size of large eddies $K_-^{-1}$
exceeds $0.5$mm. Note that there are no restrictions on~$K_+$, which means that turbulent
structures can develop details at much finer scales.

\Proof
On $\R^3$, the starting point is \eqref{Epm}. 
Using  \eqref{K41bis} on $[K_-,K_+]$, one has~:
$$\bar{E}\leq 3 \alpha \bar{\varepsilon}^{2/3} (K_-)^{-2/3}\Re^{3\gamma/4}(1-\Re^{-1/2})+
\int_0^{K_-} \bar{E}^\ast(K) dK.$$ Then \eqref{EbarLOW} gives
$\bar{E}^\ast(K) \leq \frac{4}{\pi^2}V \bar{E} K^2$ with $V=\operatorname{Vol}(u;[T_0,T_1])$,
hence~:
$$\int_0^{K_-} \bar{E}^\ast(K) dK \leq  \frac{4}{3\pi^2}V \bar{E} (K_-)^3.$$
Then \eqref{HOMOGENEOUS} reads $V (K_-)^3\leq 1$, so
$\bar{E}\leq 3 \alpha \bar{\varepsilon}^{2/3} (K_-)^{-2/3}\Re^{3\gamma/4}(1-\Re^{-1/2})
+  \frac{4}{3\pi^2}\bar{E}$
and \eqref{K-2} follows immediately.

\medskip
To get \eqref{K-2bis}, one checks the
compatibility between \eqref{K41bis} and \eqref{EbarLOW} at $K=K_-$, which provides~:
$$\Re^{-3\gamma/4}\alpha \bar{\varepsilon}^{2/3} (K_-)^{-5/3}\leq
\bar{E}^\ast(K_-)\leq \frac{4}{\pi^2} \operatorname{Vol}(u;[T_0,T_1]) (K_-)^2 \bar{E}.$$
On recognizes $\alpha\bar{\varepsilon}^{2/3}=\bar{E}K_c^{2/3}$, hence
using \eqref{K-2}~:
$$(K_-)^3\times\operatorname{Vol}(u;[T_0,T_1])  \geq 
\frac{\pi^2}{4\Re^{3\gamma/4}} \left(\frac{K_c}{K_-}\right)^{2/3}\geq
\frac{\pi^2}{4\Re^{3\gamma/4}} \left(\frac{9\pi^2}{3\pi^2-4}\times\Re^{3\gamma/4}(1-\Re^{-1/2})\right)^{-1}.
$$
The upper bound was already given by \eqref{HOMOGENEOUS}.

\medskip
For $\mathbb{T}^3$, the starting point is~\eqref{Epm2}.
As before one uses  \eqref{EbarLOW} on $[2\pi L^{-1},K_-]$~:
\begin{align*} \left(\frac{2\pi}{L}\right)
\sum_{\substack{K\in\Sigma^\ast,\\ K< K_-}} \left(\frac{2\pi}{KL}\right)\bar{E}^\dagger(K)
& \leq
C(n_-) \cdot \left(\frac{2\pi}{L}\right)^3 \left(\sum_{n=1}^{n_-}  \sqrt{n}\right)
\operatorname{Vol}(u;[T_0,T_1]) \bar{E} \\
&\leq C(n_-) \cdot  \left(\frac{2\pi}{L} \sqrt{n_-}\right)^3 \operatorname{Vol}(u;[T_0,T_1]) \bar{E}
\\& \qquad\text{with the exact same constant}\\
& = C(n_-) \cdot  (K_-)^3 \operatorname{Vol}(u;[T_0,T_1]) \bar{E}.
\end{align*}
If $\frac{\operatorname{Vol}(u;[T_0,T_1])}{\operatorname{Vol}(\mathbb{T}^3)}\geq (2\pi)^{-3}$  then \eqref{HOMOGENEOUSbis} gives $K_-=2\pi L^{-1}$ and the
sum on the left-hand side is empty so there is nothing to estimate. If not, then
\eqref{HOMOGENEOUSbis} gives $(K_-)^3 \operatorname{Vol}(u;[T_0,T_1])\leq 1$
and the sum is bounded by $C(n_-)  \bar{E}$.
Since it is assumed that $C(n_-)<1$,  one can bootstrap this term in the left-hand side.

On~$[K_-,K_+]$ one uses  \eqref{K41}, which gives as before~:
\begin{align*}
\left(\frac{2\pi}{L}\right)
\sum_{K\in\Sigma^\ast\cap[K_-,K_+]} \left(\frac{2\pi}{KL}\right) \bar{E}^\ast(K) 
& \leq \Re^{3\gamma/4} \left\{
\alpha \bar{\varepsilon}^{2/3} \left(\frac{2\pi}{L}\right)^{-2/3}
\sum_{n\in\Box_3\cap[n_-,n_+]} n^{-4/3}\right\}\\
& \leq 6\Re^{3\gamma/4} 
\alpha \bar{\varepsilon}^{2/3} K_-^{-2/3}
\qquad\text{if }n_-\geq2.
\end{align*}
For the last inequality, one has just estimated the sum on $\Box_3\cap [n_-,n_+]$ by the sum
on all integers greater than $n_-$ and then compared it to an integral.
Then \eqref{PRECISEDPERIODIC} follows from~:
$$\left(1-C(n_-)\right)\bar{E}\leq \left(\int_{\mathbb{T}^3} \rho u_0(x)dx\right)^2 +
12\Re^{3\gamma/4} \left\{\bar{E}-\left(\int_{\mathbb{T}^3} \rho u_0(x)dx\right)^2\right\}
\frac{K_-}{K_c}^{-2/3}.$$
Compatibility between \eqref{K41bis} and \eqref{EbarLOW} at $K=K_-$ then gives~:
$$(K_-)^3\times\operatorname{Vol}(u;[T_0,T_1])  \geq 
\frac{\Re^{-3\gamma/4}}{\max\{1;C(n_-)\}} \left(\frac{K_c}{K_-}\right)^{2/3}\geq
\frac{1-C(n_-)}{12\Re^{3\gamma/2}}$$
using again the assumption $C(n_-)<1$.
\cqfd\endProof

\section{High-frequency spectrum.}
\label{HIGH}

Physics textbooks often state that turbulent spectra have rapid decay at high-frequencies,
\textit{i.e.} that
\begin{equation}\label{HIGH-SPECTRUM}
\forall N\in\mathbb{N}, \qquad \sup_{K\geq0} K^N \bar{E}^\ast(K) <\infty.
\end{equation}
Property \eqref{HIGH-SPECTRUM} implies that $\bar{u}(x)$ is smooth.
Conversely, the following result relies on the best known smoothing effect for~\eqref{NS}
and implies that \eqref{HIGH-SPECTRUM} is automatically satisfied for
smooth solutions.

\begin{thm}\label{TH-HIGH}
There exist (dimensionless) numerical
constants $C$, $C_0>0$ and a family $c(K)\geq0$ with
$$\begin{cases}\displaystyle
\sum_{K\in\Sigma^\ast} c(K) = C
&\quad \text{if} \quad\Omega=\mathbb{T}^3\\[2em]\displaystyle
\forall K\in\R, \quad \sum_{n\in\mathbb{Z}} c(2^n K) \leq C
&\quad \text{if} \quad\Omega=\R^3
\end{cases}$$
that have the following properties.
For any smooth solution $u$ of \eqref{NS} on $[T_0,T_1]\times\Omega$ with
$u(T_0)\in H^1(\Omega)$, let us define
\begin{equation}\label{DEFDELTA}
\tau=\frac{\nu^3}{\sup_{[T_0,T_1]} \norme[L^2(\Omega)]{\nabla u}^4}
\qquad\text{and}\qquad \delta(t)=\frac{1}{2}\min\{\sqrt{\nu (t-T_0)};C_0\sqrt{\nu \tau}\}.
\end{equation}
Then for any $t\in [T_0,T_1]$ the energy spectrum satisfies~:
\begin{subequations}
\label{DECAYSPECTRUM}
\begin{itemize}
\item on $\Omega=\mathbb{T}^3$~:
\begin{equation}\label{Einf}
\forall K \in \Sigma^\ast, \qquad
E^\dagger(K,t) \leq
 c(K) \cdot \left(\frac{K}{L}\operatorname{Vol}(\mathbb{T}^3) \right) e^{-\delta(t) K} \: E_0,
\end{equation}
\item on $\Omega=\R^3$~:
\begin{equation}
\forall K\in\R_+^\ast, \qquad
\frac{1}{K}\int_K^{2K} E^\ast(t,k) dk
\leq  c(K) \cdot K^{-1}e^{-\delta(t) K} \: E_0,
\end{equation}
\end{itemize}
\end{subequations}
with $E_0=E(T_0)$ the initial kinetic energy.
\end{thm}
Exponential decay in~\eqref{DECAYSPECTRUM} with a uniform
constant $\delta$ on $[T_0+\tau,T_1]$ means that the analyticity radius of smooth
solutions remains uniformly bounded from below.
This question has raised some concerns in the asymptotic $\nu\to 0$
(see \cite[p.92-93]{Frish} and the references therein for a brief survey).
For the question of analyticity for a given $\nu>0$, the end of this section
contains detailed  bibliographical notes.

\bigskip
Even though the main interest of the statement is the high-frequency asymptotic,
it should be compared to the low-frequency inequality \eqref{EbarLOW} on $\mathbb{T}^3$.
\begin{cor}
If $u$ is a smooth solution of \eqref{NS} on $[T_0,T_1]\times\mathbb{T}^3$ then
for all $K\in\Sigma^\ast$~:
\begin{equation}
\bar{E}^\dagger(K) \leq C \left(\frac{K}{L}\operatorname{Vol}(\mathbb{T}^3)\right) E_0
\leq \frac{C}{2\pi} \operatorname{Vol}(\mathbb{T}^3) K^2 E_0.
\end{equation}
\end{cor}

\Proof[ of Theorem~\ref{TH-HIGH}]
The idea is to first prove a short-time analyticity estimate~\eqref{ENERGYANALYTIC}
valid for any Leray solution of~\eqref{NS} evolving from smooth initial data.
Then one uses the decay of kinetic energy $E(t)$ and the qualitative
assumption that $u$ remains smooth on $[0,\Delta]$ to propagate
this quantitative estimate along the time line.

Let us consider the following ODE with unknown function $\vartheta(t)$ and
a dimensionless constant $A$ that will be adjusted later on~:
\begin{equation}\label{ODE}
\vartheta(t)=\int_0^t \norme[L^2]{|\xi| e^{\Lambda(t')|\xi|}\hat{u}(t',\xi)}^2 dt'
\qquad\text{with}\qquad
\Lambda(t)=\sqrt{\nu t}-\frac{A\vartheta(t)}{\nu}\cdotp
\end{equation}
To be perfectly rigorous, one should consider a family of smooth approximations $u_n$
of $u$, e.g. Friedrichs approximation by low-pass filters. Then \eqref{NS} and \eqref{ODE}
are of Cauchy-Lipschitz type and can be solved for all $t\geq0$. The auxiliary
ODE~\eqref{ODE} will be used to prove an analytic estimate \eqref{ENERGYANALYTIC}
of $u_n$ that will pass to the limit $u_n\to u$ on some uniform time interval.
To keep the formulas reasonably sober, the index $n$ is dropped.

\bigskip
The image of \eqref{NS} by the Leray-Hopf projector $\mathbb{P}$ is
$$\partial_t u -\nu \Delta u +\mathbb{P}\div(u\otimes u)=0, \qquad u(0,x)=u_0(x).$$
Duhamel's formula reads
$\displaystyle \left|e^{\Lambda(t)|\xi|}\hat{u}(t,\xi)\right|\leq W_\text{L}(t,\xi)+\int_0^t
W_\text{NL}(t,t',\xi) dt'$
with
$$\displaystyle  W_\text{L}(t,\xi) = e^{\Lambda(t)|\xi|-\nu t|\xi|^2}|\widehat{u_0}(\xi)|$$
and
$$W_\text{NL}(t,t',\xi)=|\xi| \, e^{[\Lambda(t)-\Lambda(t')]|\xi|-\nu(t-t')|\xi|^2}
\left(e^{\Lambda(t')|\xi|} \int_{\R^3} |\hat{u}(t',\xi-\eta)| |\hat{u}(t',\eta)| d\eta\right).$$
To simplify notations of the time integral, one defines also~:
$$\displaystyle W_\text{NL}(t,\xi)=\int_0^t W_\text{NL}(t,t',\xi) dt'.$$
The same formula holds on $\mathbb{T}^3$ with the obvious modifications
in the notations of discrete spectra. Let us
first estimate the linear term. One has for any $t\geq0$~:
$$\Lambda(t)|\xi|-\nu t|\xi|^2 \leq \sqrt{\nu t}|\xi|-\nu t|\xi|^2\leq \frac{1}{2}-\frac{1}{2}\nu t |\xi|^2$$
therefore $\norme[L^2]{W_\text{L}(t,\cdot)}^2 \leq C \norme[L^2]{u_0}^2$
and
\begin{align}\notag
\nu \int_0^t \norme[L^2]{|\xi| W_\text{L}(t',\xi)}^2 dt' 
&\leq C \int_{\R^3} \left( \int_0^t \nu |\xi|^{2} e^{-\nu t'|\xi|^2} dt' \right)
|\widehat{u_0}(\xi)|^2 d\xi
\\\notag&= C \int_{\R^3} \left( 1-e^{-\nu t|\xi|^2} \right) |\widehat{u_0}(\xi)|^2 d\xi
\\\label{WL}&\leq C
\inf_{s\in[0,1]} (\nu t)^s \norme[L^2]{|\xi|^s\widehat{u_0}(\xi)}^2.
\end{align}
For the phase of the non-linear term, one uses the identity $(a-b)(1-\frac{a+b}{2})\leq 2$~:
\begin{align*}[\Lambda(t)-\Lambda(t')]|\xi|-\nu (t-t')|\xi|^2&=
-\frac{A|\xi|}{\nu}\left(\int_{t'}^t \frac{d\vartheta}{dt}\right) - \frac{1}{2}\nu(t-t')|\xi|^2\\
&\qquad+
|\xi|(\sqrt{\nu t}-\sqrt{\nu t'})\left( 1 - \frac{1}{2}(\sqrt{\nu t}+\sqrt{\nu t'})|\xi| \right) \\
&\leq -\frac{A|\xi|}{\nu}\left(\int_{t'}^t \frac{d\vartheta}{dt}\right) - \frac{1}{2}\nu(t-t')|\xi|^2 +2.
\end{align*}
The sub-linearity of $\xi\mapsto|\xi|$ provides~:
$$e^{\Lambda(t')|\xi|} \int_{\R^3} |\hat{u}(t',\xi-\eta)| |\hat{u}(t',\eta)| d\eta \leq 
\mathcal{F}\left[V(t',\cdot)^2\right](\xi)
\quad\text{with}\quad
\hat{V}(t',\xi)=e^{\Lambda(t')|\xi|}|\hat{u}(t',\xi)|.$$
One has $\norme[\dot{H}^1]{V(t)}^2=\norme[L^2]{\nabla V(t)}^2=\frac{d\vartheta}{dt}$.
A classical property of Sobolev-Besov
spaces is that the product $(f,g)\mapsto fg$ maps continuously $\dot{H}^1\times
\dot{H}^1$ to $\dot{H}^{1/2}$ ; therefore~:
$$0\leq
\mathcal{F}\left[V(t',\cdot)^2\right](\xi) \leq C  \frac{d\vartheta}{dt}(t') \cdot g(t',\xi) \, |\xi|^{-1/2}
\quad\text{with}\quad \norme[L^2]{g(t',\cdot)}\leq1.$$
Putting everything together, the following bound holds for the non-linear term~:
\begin{align}
\notag
\norme[L^2]{W_\text{NL}(t,\xi)}^2
 &\leq C \int_{\R^3}\left|\int_0^t
|\xi|^{1/2} e^{2-\frac{A|\xi|}{\nu}\int_{t'}^t\frac{d\vartheta}{dt}(\tau)d\tau - \frac{1}{2}\nu(t-{t'})|\xi|^2}\,\frac{d\vartheta}{dt}({t'})\,g({t'},\xi) \,d{t'}\right|^2 d\xi\\\notag
&\leq C \int_{\R^3}\left(\int_0^t \frac{d\vartheta}{dt}({t'}) |\xi| e^{-2\frac{A|\xi|}{\nu}\int_{t'}^t\frac{d\vartheta}{dt}(\tau)d\tau} \,d{t'}\right)\left(\int_0^t e^{- \nu(t-{t'})|\xi|^2}\,\frac{d\vartheta}{dt}({t'})\,g({t'},\xi)^2 \,d{t'}\right) d\xi\\\notag
&\leq \frac{C \nu}{A} \int_{\R^3}\left(\frac{1-e^{-2\frac{A|\xi|}{\nu}\vartheta(t)}}{2}\right)\left(\int_0^t e^{- \nu(t-{t'})|\xi|^2}\,\frac{d\vartheta}{dt}({t'})\,g({t'},\xi)^2 \,d{t'}\right) d\xi
\\\label{WNL}&\leq \frac{C \nu}{A}\,\vartheta(t).
\end{align}
Along the same lines, one has~:
\begin{align}\notag
\nu \int_0^t \norme[L^2]{|\xi| W_\text{NL}(t',\xi)}^2  dt'  &\leq
\frac{C\nu^2}{A}\int_0^t 
\int_{\R^3}\left(1-e^{-\frac{2A|\xi|}{\nu}\vartheta(t')}\right)\left(\int_0^{t'}  |\xi|^2 e^{- \nu(t'-t'')|\xi|^2}\,\frac{d\vartheta}{dt}(t'')\,g(t'',\xi)^2 \,dt''\right) d\xi dt'
\\\notag&\leq \frac{C\nu^2}{A} \int_0^t \int_{\R^3} \left(\int_{t''}^t
|\xi|^2 \,e^{- \nu(t'-t'')|\xi|^2}\,dt'\right) \frac{d\vartheta}{dt}(t'')\,g(t'',\xi)^2 d\xi dt''
\\\notag&\leq \frac{C\nu}{A}\int_0^t \int_{\R^3} \left(
1-e^{- \nu(t-t'')|\xi|^2}\right) \frac{d\vartheta}{dt}(t'')\,g(t'',\xi)^2 d\xi dt''
\\\label{WNLT}
&\leq  \frac{C\nu}{A}\,\vartheta(t).
\end{align}
Estimations \eqref{WL} and \eqref{WNLT} imply the following bootstrap~:
$$\vartheta(t) \leq \frac{C}{\nu} \inf_{s\in[0,1]} (\nu t)^s \norme[L^2]{|\xi|^s\widehat{u_0}(\xi)}^2
+\frac{C'}{A}\vartheta(t) \leq
\frac{2C}{\nu} \inf_{s\in[0,1]} (\nu t)^s \norme[L^2]{|\xi|^s\widehat{u_0}(\xi)}^2$$
provided the numerical constant $A$ is chosen sufficiently large.
To conclude the analytic estimate of $u$, let us now focus on the time interval on
which $\Lambda(t)\geq \frac{1}{2}\sqrt{\nu t}$. Let us therefore define~:
$$T^\ast = \inf\left\{t>0 \,;\, \forall t'\in[0,t], \enspace
\displaystyle \vartheta(t')\leq \frac{\nu^{3/2} (t')^{1/2}}{2A}\right\}.$$
The bootstrap inequality provides\footnote{In particular,
one can pass to the limit $u_n\to u$ on $[0,T^\ast]$ and this method guarantees smoothness. Note that $T^\ast=+\infty$
if $\norme[\dot{H}^{1/2}]{u_0}\leq \nu\sqrt{C_0}$ so
this method also provides a proof of Kato's theorem \cite{FujitaKato}.
}~:
$$T^\ast \geq \sup_{s\in]\frac{1}{2},1]} \left( \frac{ C_0 \nu^{\frac{5}{2}-s}}{
\norme[\dot{H}^s]{u_0}^{2}} \right)^{\frac{1}{s-\frac{1}{2}}} \geq
\frac{ C_0^2 \nu^3}{\norme[L^2]{\nabla u_0}^4}.$$
On $[0,T^\ast]$, one has $\frac{1}{2}\sqrt{\nu t}\leq \Lambda(t) \leq \sqrt{\nu t}$  hence
\eqref{WNL} and again the bootstrap inequality gives~:
\begin{equation}\label{ENERGYANALYTIC}
\int_{\R^3}{e^{\sqrt{\nu t}|\xi|} |\hat{u}(t,\xi)}|^2d\xi \leq 
\int_{\R^3} \left|e^{\Lambda(t)|\xi|}\hat{u}(t,\xi)\right|^2 d\xi \leq
\left\{\norme[L^2]{W_\text{L}(t,\cdot)}+ \norme[L^2]{W_\text{NL}(t,\cdot)}\right\}^2
\leq C \norme[L^2]{u_0}^2.
\end{equation}
This inequality proves that the radius of analyticity of $u(t,\cdot)$
exceeds $\frac{1}{2}\sqrt{\nu t}$ on $[0,T^\ast]$. The constant $C$ is
purely numerical. This estimates also holds with an obvious change of notations
on $\mathbb{T}^3$.

\bigskip
Let us now assume that $u$ is a smooth solution
of \eqref{NS} on $[0,T]\times\R^3$ or $[0,T]\times\mathbb{T}^3$
with $T$ possibly much larger than~$T^\ast$,
one can translate \eqref{ENERGYANALYTIC} along the time line in the following way.
One defines~:
$$\tau=\frac{\nu^3}{\sup_{[0,T]}\norme[L^2(\Omega)]{\nabla u}^4} 
\qquad\text{and}\qquad
\delta_0(\tau) = \sqrt{\nu \tau}.$$
Then for any time $t_0\in[0,T]$,  the estimate \eqref{ENERGYANALYTIC} holds
on $[t_0,t_0+C_0^2\tau]$ with the same constant $C$~:
\begin{equation}\label{GLOBALENERGYANALYTIC}
\forall t_0\in[0,T],\enspace \forall t\in[t_0,t_0+C_0^2\tau], \qquad 
C^{-1}\norme[L^2]{e^{\frac{1}{2}\sqrt{-\nu (t-t_0)\Delta}} u(t,\cdot)}^2 \leq \norme[L^2]{u(t_0,\cdot)}^2
\leq  \norme[L^2]{u_0}^2.
\end{equation}
The best estimate at time $t$ is therefore obtained by choosing $t_0=t-C_0^2\tau$ provided $t>C_0^2\tau$. When $t\leq C_0^2\tau$, one can only rely on the initial estimate with $t_0=0$.

\bigskip
One can now conclude the proof of Theorem~\ref{TH-HIGH} from a straighforward computation.
As the Bessel-Parseval formula is slightly different let us detail it~:
\begin{equation}\label{PARSEVAL}
\norme[L^2(\Omega)]{f}^2 = \begin{cases}
\displaystyle
\frac{1}{(2\pi)^3}
\int_{\R^3} |\hat{f}(\xi)|^2 d\xi & \text{on }\R^3,
\\[1em]\displaystyle
\frac{1}{\operatorname{Vol}(\mathbb{T}^3)}\sum_{k\in\hat{\mathbb{T}}^3} |\hat{f}(k)|^2 \quad
& \text{on }\mathbb{T}^3.
\end{cases}
\end{equation}
To simplify notations, let us denote by $\delta^\ast(t)=\min\{\delta_0(t),C_0\delta_0(\tau)\}$.
On $[0,T]\times\mathbb{T}^3$, the inequality \eqref{GLOBALENERGYANALYTIC} reads~:
$$\frac{1}{\operatorname{Vol}(\mathbb{T}^3)}\sum_{k\in\hat{\mathbb{T}}^3} e^{\delta^\ast |k|} |\hat{u}(t,k)|^2 \leq C\norme[L^2]{u_0}^2$$
and \eqref{E*PERIODIC} gives a sequence of dimensionless coefficients $c_K\geq 0$ such that~:
\begin{equation}
E^\ast(K,t) = (2\pi)^{-2}
\rho \left(\frac{K}{L}\right) e^{-\delta^\ast K} \sum_{|k|=K}e^{\delta^\ast |k|}  |\hat{u}(t,k)|^2 \leq
C \operatorname{Vol}(\mathbb{T}^3) \left(\frac{K}{L}\right) e^{-\delta^\ast K} E(0) \cdot c_K  \quad \text{with}\quad \sum_{K} c_K\leq1.
\end{equation}
On $[0,T]\times\R^3$, the inequality \eqref{GLOBALENERGYANALYTIC} does not
directly provide a pointwise estimate of $E^\ast(t,K)$ so one considers instead~:
$$\widetilde{E}^\ast(K,t) = \frac{1}{K}\int_K^{2K} E^\ast(t,k) dk$$
The spectrum is given by \eqref{E*}. One gets~:
$$\widetilde{E}^\ast(K,t) = (2\pi)^{-3}\frac{\rho}{K}\int_{\R^3}\left(\int_K^{2K}
\frac{e^{-\delta^\ast |\xi|}}{k}\,\psi\left(\frac{|\xi|}{k}\right) dk\right)
e^{\delta^\ast |\xi|} |\hat{u}(t,\xi)|^2 \frac{d\xi}{(2\pi)^3}.$$
One checks easily, using the properties of the cut-off function $\psi$ (uniform in the limit $\psi(r)\rightharpoonup \delta_{r=1}$), that
$$\int_K^{2K} \frac{e^{-\delta^\ast |\xi|}}{k}\,\psi\left(\frac{|\xi|}{k}\right) dk \leq 
e^{-\frac{1}{2}\delta^\ast K}
\int_{|\xi|/2K}^{|\xi|/K} \psi(r) dr \leq 2 c(K) e^{-\frac{1}{2}\delta^\ast K}$$
with $\sum_{n\in\mathbb{Z}} c(2^n K) \leq 1$.
 Inequality \eqref{GLOBALENERGYANALYTIC}  gives~:
 \begin{equation}
\widetilde{E}^\ast(K,t)\leq 2(2\pi)^{-3}e^{-\frac{1}{2}\delta^\ast K}\frac{c(K)}{K}\rho
\int_{\R^3} e^{\delta^\ast |\xi|} |\hat{u}(t,\xi)|^2 \frac{d\xi}{(2\pi)^3}
\leq C K^{-1}e^{-\frac{1}{2}\delta^\ast K} E(0)\cdot c(K).
\end{equation}
This concludes the proof of \eqref{DECAYSPECTRUM}
and of Theorem~\ref{TH-HIGH}. 
\cqfd\endProof

\paragraph{Bibliographical note on the analytic estimate \eqref{GLOBALENERGYANALYTIC}.}
For physical reasons, one was interested only in $L^2$ norms, therefore a pointwise
majoration in Fourier variables was sufficient. One should notice however the
following idea suggested by \cite{L04}. Instead of using $|\xi|\leq |\xi-\eta|+|\eta|$, one can
rely on the exact formula~:
$$\forall\alpha\in\R^3,\qquad e^{i\alpha\cdot\nabla}(uv)=
\left(e^{i\alpha\cdot\nabla} u\right)\left(e^{i\alpha\cdot\nabla} v\right).$$
One could therefore replace \eqref{ODE} by a family of ODEs
$$\vartheta_\alpha(t) = \int_0^t \norme[L^p]{\nabla e^{i\Lambda_\alpha(\tau)\,\alpha\cdot\nabla}u(\tau)}^2 d\tau
\quad\text{with}\quad \Lambda_\alpha(t)=\sqrt{\nu t}-\lambda \vartheta_\alpha(t)$$
that would provide uniform bounds on $\norme[L^p]{e^{\sqrt{\nu t}\,i\alpha\cdot\nabla}u(t)}^2$
for $t\in[0,T^\ast]$ and $\alpha\in B=\{\beta\in\R^3\,;\,|\beta|\leq C\}$. One can then deduce the local-in-time
$L^p$-analytic estimate of $u$ because an elementary exercise
in Littlewood-Paley theory states that~:
$$\norme[L^p]{e^{\sqrt{-\nu t\Delta}}f}\leq C \sup_{\alpha\in B}
\norme[L^p]{e^{\sqrt{\nu t}\,i\alpha\cdot\nabla}f}^2.$$
However, due to the lack of a uniform a-priori bound of $\norme[L^p]{u(t)}$ when $p\neq2$,
one can only propagate the local estimate along the time-line by stating that
\begin{equation}\label{GLOBALLpANALYTIC}
\text{$u$ is smooth, namely }u\in L^\infty([0,T]\times\Omega)  \quad\Longrightarrow\quad
e^{\sqrt{- \min\{\nu t,\nu t_\ast\}\Delta}} u(t)\in L^p
\end{equation}
but unlike \eqref{GLOBALENERGYANALYTIC}, the corresponding inequality
involves both $u_0$ and $\sup_{[0,T]}\norme[L^p(\Omega)]{u}$ to compute $t_\ast$
and the~$L^p$-analytic norm.
Statements similar to~\eqref{GLOBALLpANALYTIC} can be found \textit{e.g.} in
\cite{L00}, \cite{L04}, \cite{Lemarie} or \cite{GPS09}.
However, as these authors make the wise choice not to state the
corresponding estimate,
the status of \eqref{GLOBALENERGYANALYTIC} was unclear.
In particular, estimates hidden in \eqref{GLOBALLpANALYTIC} cannot usually be translated 
in time with uniform constants.
For the convenience of the reader I decided to provide an
elementary proof of~\eqref{GLOBALENERGYANALYTIC}.
The method was greatly inspired by \cite{C04},
where a statement similar to \eqref{ENERGYANALYTIC} was proved
for the~$L^2([0,T]\times\R^3)$-norm of $e^{\sqrt{-\nu t\Delta}}u(t)$.

Some subtle mathematical problems are closely related to Theorem~\ref{TH-HIGH}.
Analyticity in time and Gevrey classes have been studied in \cite{FT89}. The
difficult question of analyticity for general domains seems to have
been studied only once, in \cite{M67},  and the question of external forces is dealt
with in \cite{GK98}.
Readers interested in analyticity questions for 2D-turbulence should refer
\textit{e.g.}~to \cite{BCFM-95}.

\section{Necessary conditions satisfied by turbulent flows}
\label{IFTHM}

In this section, one investigates the necessary conditions satisfied by K41-turbulent
flows. One proves the relation between Integral and Taylor
Scale Reynolds numbers (Theorem \ref{REYNOLDSTHM}). On~$\mathbb{T}^3$,
this relation is shown to characterize a time-scale on which free turbulence
can be observed (Theorem \ref{TIMESCALETHM}).
Finaly one investigates whether smooth
solutions can be K41-turbulent (Theorems \ref{SMOOTHTURBULENCE}
and \ref{SMOOTHTURBULENCE2}).

\subsection{Reynolds numbers}\label{RN}

The Reynolds number $\Re=\left(\frac{K_+}{K_-}\right)^{4/3}$
is the so-called
\textsl{Integral Scale} Reynolds number  \cite[\S7.18 p.107]{Frish}.
In experiments, one uses mostly the \textsl{Taylor Scale} Reynolds
number $\frac{U_\text{rms} \lambda_\text{T}}{\nu}$
with $U_\text{rms}(t) = \sqrt{E(t)}$ and the Taylor scale defined by $\frac{1}{\lambda_\text{T}^2(t)}=\frac{\norme[L^2]{\nabla u(t)}^2}{\norme[L^2]{u(t)}^2}$ ; in other words~:
$$\frac{U_\text{rms} \lambda_\text{T}}{\nu} =
\frac{\rho^{1/2} \norme[L^2]{u(t)}^2}{\nu \norme[L^2]{\nabla u(t)}}\cdotp$$
Usually  $R_\lambda$ is ``tailored'' to the needs of each experiment to get a
time-independent number.

\begin{dfn}
If $u$ is a Leray solution of \eqref{NS}, 
the Taylor-Scale Reynolds number of $u$ on $[T_0,T_1]$ is~:
\begin{equation}
R_\lambda = \frac{\rho^{1/2} \langle \norme[L^2]{u}^2\rangle}
{\alpha^{3/2} \nu \langle\norme[L^2]{\nabla u}^2\rangle^{1/2}}
= \sqrt{\frac{\bar{E}^2}{\alpha^3 \nu \bar{\varepsilon}}}\cdotp
\end{equation}
\end{dfn}
Observations \cite[7.17 p.107]{Frish} suggest that $R_\lambda \simeq \Re^{1/2}$.
This is indeed a rigorous fact.

\begin{thm}\label{REYNOLDSTHM}
Assume that $u$ is a Leray solution of \eqref{NS} and a K41-function.
\begin{enumerate}
\item For $\Omega=\R^3$, assume also that $u_0\in L^1\cap L^2$. Then~:
\begin{equation}
\frac{\Re^{-9\gamma/4}}{216(1-\Re^{-1})(1-\Re^{-1/2})^2}
\leq \frac{\Re}{R_\lambda^2} \leq  \frac{16\Re^{9\gamma/4}}{27(1-\Re^{-1})(1-\Re^{-1/2})^2}.
\end{equation}
In particular~:
\begin{equation}\label{EPSILONANDRE}
\bar\varepsilon  \simeq \frac{\bar{E}^2}{\alpha^3\nu\Re} \simeq
\frac{\bar{E}^{3/2}}{\alpha^{3/2} \operatorname{Vol}(u;[T_0,T_1])^{1/3}}
\quad\text{and}\quad
\Re \simeq \frac{\bar{E}^{1/2}  \operatorname{Vol}(u;[T_0,T_1])^{1/3}}{\alpha^{3/2}\nu}\cdotp
\end{equation}
According to \eqref{DICTIONARY}, the symbol $\simeq$ means that both inequalities hold with
constants depending on $\Re$ and $\Re^\gamma$ but that the constants
have a purely numerical limit in the asymptotic \eqref{K41}.
\item
For $\mathbb{T}^3$, let us assume that $\displaystyle \int_{\mathbb{T}^3} \rho u_0(x)dx=0$
and that \eqref{NBTHASSUMPTION} is satisfied.
Then similar inequalities hold. The numerical values are~:
\begin{equation}\label{RonRlambda2}
\frac{(1-C(n_-))^2}{864} + o(\Re;\gamma)\leq \frac{\Re}{R_\lambda^2}
\leq \frac{2048}{3375} + o(\Re;\gamma)
\end{equation}
where $o(\Re;\gamma)\to0$ in the asymptotic \eqref{K41}.
\end{enumerate}
\end{thm}
Note that the right-hand side of~\eqref{RonRlambda2} requires only a lower bound
on $K_-/K_c$, thus it holds regardless of~\eqref{NBTHASSUMPTION}.
\Proof
Using the definition of $\Re$ and \eqref{Kd}, \eqref{Kc} for $K_d$ and $K_c$, one gets~:
$$\Re=\left(\frac{K_+}{K_-}\right)^{4/3} \quad\text{and}\quad
\left(\frac{K_d}{K_c}\right)^{4/3} = \frac{\bar{E}^2}{\alpha^3\nu\bar{\varepsilon}}
= R_\lambda^2$$
thus Theorem~\ref{RANGETHM} implies $\Re\simeq R_\lambda^2$.
Proposition \ref{K-PRECISED} provides $\operatorname{Vol}(u;[T_0,T_1])\simeq (K_-)^{-3}$
hence $$\Re \simeq  \left(\operatorname{Vol}(u;[T_0,T_1])^{1/3} K_d\right)^{4/3}
= \frac{\bar{\varepsilon}^{1/3}}{\alpha \nu}\:\operatorname{Vol}(u;[T_0,T_1])^{4/9}.$$
Using $\Re\simeq R_\lambda^2$ provides
$$\bar\varepsilon \simeq
\left(\frac{\bar{E}^2}{\alpha^2 \operatorname{Vol}(u;[T_0,T_1])^{4/9}}\right)^{3/4}
= \frac{\bar{E}^{3/2}}{\alpha^{3/2} \operatorname{Vol}(u;[T_0,T_1])^{1/3}}\cdotp$$
Substitution of this formula in the previous expression of $\Re$ yields the
last formula in~\eqref{EPSILONANDRE}.
In the case of~$\Omega=\mathbb{T}^3$, the proof is similar.
\cqfd\endProof

\subsection{Time-scale on which free turbulence can be observed}
\label{PAR:TIMESCALE}

\begin{dfn}
Given a Leray solution $u$ of \eqref{NS} on $\Omega=\mathbb{T}^3$
with $\displaystyle \int_{\mathbb{T}^3} \rho u_0=0$, let us define the following time scale,
that one could call for example the ``transfer time from $u_0$ to
$\omega=\operatorname{rot}u$ on $[T_0,T_1]$''~:
\begin{equation}
\mathcal{T}(u_0;\omega) =\frac{\alpha^3\nu^2}{\rho}\times
\frac{\displaystyle \int_{T_0}^{T_1} \norme[L^2]{\omega(t)}^2 dt}{\norme[L^2]{u_0}^4}\cdotp
\end{equation}
\end{dfn}
According to \eqref{LERAYINEQ}, one has
 $\mathcal{T}(u_0,\omega)\leq \frac{\alpha^3\nu}{E_0}\cdotp$
The following statement proves that $\Re\,\mathcal{T}(u_0,\omega)$ is the precise
time-scale on which turbulence can be observed~:
for a two short observation time one will not see K41-properties and
for a too long observation time the time-average will describe the fluid as being mostly at rest.
This time-scale appears also naturally in the computation of $\bar{\varepsilon}$
for a smooth flow~(see~\eqref{EPSILONSMOOTH} below). 

\begin{thm}\label{TIMESCALETHM}
If $u$ is a K41-function on $[T_0,T_1]\times\mathbb{T}^3$ and
a Leray solution of \eqref{NS} with $\displaystyle \int_{\mathbb{T}^3}\rho u_0(x)dx=0$,
the following inequality holds for numerical factors $C_{j}(\Re,\gamma)\to 1$
in the asymptotic \eqref{K41}~:
\begin{equation}\label{TIMESCALEEQ}
\frac{3375}{2048}C_1(\Re,\gamma) \: \Re \:\mathcal{T}(u_0;\omega)
\leq  T_1-T_0 \leq 
\frac{128 C_2(\Re,\gamma)}{3375 \pi^4}\:
\frac{1}{\Re\:\mathcal{T}(u_0;\omega)}\times
\frac{\operatorname{Vol}(\mathbb{T}^3)^{4/3}}{\nu^2}\cdotp
\end{equation}
Moreover, if~\eqref{NBTHASSUMPTION} is satisfied, then~:
\begin{equation}\label{DeltaTupperbound}
\frac{3375}{2048}C_3(\Re,\gamma) \leq 
\left(\frac{\bar{E}}{E_0}\right)^2  \frac{T_1-T_0}{\Re \:\mathcal{T}(u_0;\omega)}
\leq \frac{864 C_4(\Re,\gamma)}{(1-C(n_-))^2}.
\end{equation}
\end{thm}
In particular, \eqref{TIMESCALEEQ} gives~: 
\begin{equation}
\frac{\Re}{\Re_c}\leq \frac{512 C_5(\Re,\gamma)}{3375\pi^2}
\qquad\text{with}\qquad
\Re_c = \frac{\operatorname{Vol}(\mathbb{T}^3)^{2/3}}
{\nu\mathcal{T}(u_0;\omega)}\cdotp
\end{equation}

\Proof
Poincar\'{e}'s inequality reads~:
$$\norme[L^2(\mathbb{T}^3)]{u(t)-\int_{\Omega} u(t,x)dx} \leq \frac{\operatorname{Vol}(\mathbb{T}^3)^{1/3}}{2\pi} \norme[L^2(\mathbb{T}^3)]{\nabla u(t)}.$$
On $\mathbb{T}^3$, the total impulsion is preserved so \eqref{INT0} holds for any $t\geq T_0$.
Combined with \eqref{ENERGYCONSERVATION} one gets~:
$$\norme[L^2]{u(t)}^2 + 2\nu \left(\frac{2\pi}{\operatorname{Vol}(\mathbb{T}^3)^{1/3}}\right)^2\int_{T_0}^t \norme[L^2]{u(t')}^2dt' \leq \norme[L^2]{u(T_0)}^2$$
and Gronwall's lemma provides (see \cite{KO92} for a numerical confirmation)~:
$$\norme[L^2(\mathbb{T}^3)]{u(t)} \leq \norme[L^2(\mathbb{T}^3)]{u(T_0)}
\exp\left(-\left(\frac{2\pi}{\operatorname{Vol}(\mathbb{T}^3)^{1/3}}\right)^2 \nu(t-T_0)\right).$$
In turn, this gives~:
$$\bar{E}(K) \leq E_0 \Phi\left(
\left(\frac{2\pi}{\operatorname{Vol}(\mathbb{T}^3)^{1/3}}\right)^2 \nu(T_1-T_0)
\right)\qquad\text{with}\qquad
\Phi(s)=\frac{1-e^{-s}}{s}\cdotp$$
Combining the definition of $\bar{\varepsilon}$ and $R_\lambda$ gives~:
\begin{equation}\label{COMPUT:intomega2}
\int_{T_0}^{T_1}\norme[L^2]{\omega(t)}^2 dt =
\frac{(T_1-T_0)\bar{E}^2}{\alpha^3\rho\nu^2\Re} \times
\frac{\Re}{R_\lambda^2}\cdotp
\end{equation}
Theorem~\ref{REYNOLDSTHM} on $\mathbb{T}^3$
ensures $\Re\leq(\frac{2048}{3375}+o(1)) R_\lambda^2$,
which now reads
$$\int_{T_0}^{T_1}\norme[L^2]{\omega(t)}^2 dt \leq \left(\frac{2048}{3375}+o(1)\right) 
\frac{E_0^2}{\alpha^3\rho\nu^2\Re} \times (T_1-T_0)\Phi^2\left(
\left(\frac{2\pi}{\operatorname{Vol}(\mathbb{T}^3)^{1/3}}\right)^2 \nu(T_1-T_0)
\right)$$
and $\Phi^2(s)\leq \min\{1;1/s^2\}$ provides both upper and lower
estimates on $T_1-T_0$ in~\eqref{TIMESCALEEQ}.
Note that one does not even require~\eqref{NBTHASSUMPTION}, for this assumption
is only needed to get a more precise upper-bound of $T_1-T_0$.

\medskip
Conversely, if~\eqref{NBTHASSUMPTION} holds,
then $\Re\geq ((1-C(n_-))^2/864+o(1))R_\lambda^2$ and~\eqref{COMPUT:intomega2}
boils down to~\eqref{DeltaTupperbound}.
\cqfd\endProof

\bigskip
\paragraph{Compatibility with the time scale on which smoothness is guaranteed.}
Mathematicians can guarantee the smoothness
of the solution of \eqref{NS} on at least $[T_0,T_0+\Theta]$ with \textit{e.g.}
$$\Theta =\frac{C_0 \nu^3}{\norme[L^2]{\omega_0}^4}\cdotp$$
Is such an interval long enough for the observation of free turbulence ?
According to Theorem~\ref{TIMESCALETHM},
the answer is \textbf{yes provided $\Theta \gtrsim \Re\,\mathcal{T}(u_0,\omega)$}.
Indeed, if $T_1-T_0\leq \Theta \leq C\,\mathcal{T}(u_0,\omega)$
then the only possible Reynolds number on $[T_0,T_1]$ is~$\Re\lesssim C$ and the turbulent
asymptotic~\eqref{K41} cannot be achieved.

\medskip
Inequality $\Theta\gtrsim \Re\,\mathcal{T}(u_0,\omega)$ is equivalent to
$$\rho \nu\norme[L^2]{u_0}^4 \gtrsim \alpha^3 \Re \norme[L^2]{\omega_0}^4
\int_{T_0}^{T_1} \norme[L^2]{\omega}^2,$$
which, combined with Poincar\'{e}'s inequality and $\rho=\operatorname{Vol}(\Omega)^{-1}$, implies~:
\begin{equation}\label{BIZARE}
\frac{\alpha^3 \Re}{\rho\nu}\,(T_1-T_0) \bar{\varepsilon}
=\alpha^3 \Re \int_{T_0}^{T_1} \norme[L^2]{\omega}^2
\lesssim \nu \operatorname{Vol}(\Omega)^{1/3}.
\end{equation}
For large initial data, this makes it impossible for the initial energy to dissipate
almost completely on $[T_0,T_1]$ because the left-hand side would then
be equivalent to~$\alpha^3 \Re(\rho\nu)^{-1} E_0\gg \nu \operatorname{Vol}(\Omega)^{1/3}$.

\bigskip
To put it simply, the time-scale $\Theta$ on which regularity is guaranteed 
is too short to observe a fully developed free turbulence and thus it
might not have a deep physical meaning.
One should however refrain from jumping to the conclusion
that K41-turbulence is an obstruction to smoothness.
If one keeps~\eqref{HIGH-SPECTRUM} in mind (which is accepted in
every physics textbook), it indicates instead
that \eqref{NS} can develop a very specific dynamic in Fourier space and that
the techniques used to prove local smoothness have failed to capture it.
We will see in \S\ref{parTROUBLING}, that K41-turbulence is
troublingly close to the best known local smoothing effect.

\subsection{Two necessary conditions satisfied by smooth turbulent flows}
\label{NECESSARY}

The goal of this last section is to investigate the necessary conditions that occur
when a smooth solution~$u$ of~\eqref{NS} happens to be a K41-function.
Two Theorems can be stated.

\subsubsection{Necessity of intermittency}

In section \S\ref{TAI}, temporal intermittency was defined as a substantial deviation
between $\varepsilon(t)$ and $\bar{\varepsilon}$.
If $u$ is smooth on $[T_0,T_1]$, one can compute $\bar{\varepsilon}$
with~\eqref{SMOOTH-EPSILON}~:
$$\bar{\varepsilon} = \frac{E_0-E_1}{T_1-T_0}\cdotp$$
The conservation of energy reads~:
$$E_0-E_1 = \rho\nu \int_{T_0}^{T_1} \norme[L^2]{\omega(t)}^2dt
= \frac{E_0^2}{\alpha^3\nu}  \: \mathcal{T}(u_0;\omega)$$
hence 
\begin{equation}\label{EPSILONSMOOTH}
\bar{\varepsilon}  = \frac{E_0^2}{\alpha^3\nu}  \: \frac{\mathcal{T}(u_0;\omega)}{T_1-T_0}\cdotp
\end{equation}

\begin{thm}\label{SMOOTHTURBULENCE}
If $u$ is a \textbf{smooth} solution of \eqref{NS} \textit{i.e.} $u\in L^\infty([T_0,T_1]\times\Omega)$ with $\Omega=\R^3$ or $\mathbb{T}^3$ and a K41-function on~$[T_0,T_1]$,
the following condition must be satisfied~:
\begin{equation}\label{NINTERM}
\frac{E_0+E_1}{2E_0}+\int_{T_0}^{T_1}\frac{|\varepsilon(t)-\bar\varepsilon|}{E_0} dt \geq
C(\Re,\gamma)\sqrt{\frac{\Re \: \mathcal{T}(u_0;\omega)}{T_1-T_0}}
\times \begin{cases}
\frac{3\sqrt{3}}{4}& \text{if }\Omega=\R^3,\\[1ex]
\frac{15\sqrt{15}}{32\sqrt{2}} & \text{if }\Omega=\mathbb{T}^3,
\end{cases}
\end{equation}
with  $C(\Re,\gamma)\to1$ in the asymptotic \eqref{K41}.
The numerical constant is bounded from below by $1.299$ on $\R^3$ and by~$1.283$
on $\mathbb{T}^3$.
\end{thm}\Proof
Theorem~\ref{REYNOLDSTHM}  provides in the asymptotic \eqref{K41}~:
$$\bar{\varepsilon} = \frac{\bar{E}^2}{\alpha^3\nu \Re} \times \frac{\Re}{R_\lambda^2}
\leq \frac{\bar{E}^2}{\alpha^3\nu \Re} \times \begin{cases}
\frac{16}{27}+o(1) & \text{if }\Omega=\R^3,\\[1ex]
\frac{2048}{3375}+o(1) & \text{if }\Omega=\mathbb{T}^3.
\end{cases}$$
Estimate \eqref{INTERMITTENCY} reads~:
$$\bar{E}^2\leq 
\left(\frac{E_0+E_1}{2}+\int_{T_0}^{T_1}|\varepsilon(t)-\bar{\varepsilon}|dt\right)^2$$
and \eqref{NINTERM} follows immediately from the comparison
with~\eqref{EPSILONSMOOTH}.
\cqfd\endProof

\bigskip
One can wonder whether~\eqref{NINTERM} actually provides a lower
bound on intermittency.
For $\Omega=\mathbb{T}^3$, the answer is subtle but positive.
According to Theorem~\ref{TIMESCALETHM} and a careful track of constants,
the right-hand side of~\eqref{NINTERM} cannot asymptotically exceed $1$.
On the other hand, the energy estimate gives
$\frac{E_0+E_1}{2E_0}\in [1/2;1].$
Therefore, if one assumes that~:
\begin{itemize}
\item
most of the initial energy has been dissipated, \textit{i.e.} $E_1\ll E_0$,
then $\frac{E_0+E_1}{2E_0}\simeq \frac{1}{2}$
\item
and that the observation time is the minimal, \textit{i.e.} $T_1$ is chosen short
enough such that the right-hand side of~\eqref{NINTERM} belongs to $]\frac{1}{2},1]$
\end{itemize}
then the left-hand side can be bootstrapped in the right-hand side and
inequality \eqref{NINTERM} indeed provides a lower bound of intermittency~:
$$\int_{T_0}^{T_1}|\varepsilon(t)-\bar\varepsilon| dt \geq \frac{C}{2} E_0$$
with $C\in]0,1]$.
Given smooth data $u_0$, it is not clear whether one can find $T_1$ such that
both conditions are simultaneously satisfied (the problem is to prove
 that $E(T_0+\Re \:\mathcal{T}(u_0;\omega)) \leq c E(T_0)$ with a sufficiently
small numerical constant $c<1$).
However, physical intuition suggests that this is the case since the
best time-scale to observe free turbulence is to wait till most of the initial energy has been
dissipated but not any longer.

\subsubsection{Compatibility of the $K^{-5/3}$ law with the analytic smoothing effect}

The question dealt with in this section is the following :
are the finer scales of turbulent vortex structures
limited by the analyticity radius of the solution ? In other words, if $\delta$ denotes the
analyticity radius of a K41-turbulent solution $u$, what is the possible range
of $\delta K_+$ ?
In the regime $\delta K_+\leq C$ the finer scale of turbulent structures is~$K_+^{-1}$
and is limited from below by $C^{-1}\delta$, which drastically limits 
the possible Reynolds numbers~:
$$\Re \leq \left(\frac{C}{\delta K_-}\right)^{4/3}\quad\text{and in particular}\quad
\Re\leq \left(\frac{C}{2\pi} \frac{L}{\delta}\right)^{4/3}
\quad\text{on }\mathbb{T}^3.$$
Conversely, the regime $\delta K_+\gg1$ means that turbulent structures exist
at much finer scales than the analyticity radius and the asymptotic~\eqref{K41}
remains possible (see also~\S\ref{PAR:23LAW} below).

\medskip
The following result denotes the
compatibility at $K=K_+$ between K41-property \eqref{K41bis} and 
the high-frequency bound on the spectrum given by Theorem~\ref{TH-HIGH}.
There are two cases depending on whether the initial data
is supposed to have no additional smoothness to being $L^\infty\cap H^1(\Omega)$
or if on the contrary, one considers a flow already ``well prepared'' by
the analytic smoothing effect.

\begin{thm}\label{SMOOTHTURBULENCE2}
Let $u$ be a smooth solution of \eqref{NS} on $[T_0,T_1]\times\Omega$
with $\Omega=\R^3$ or $\mathbb{T}^3$ and a K41-function on~$[T_0,T_1]$.
There exists a numerical constant $C>0$ such that the following conditions are satisfied.
\begin{description}
\item[1. Unprepared data on $\mathbb{T}^3$ --]
If $u(T_0)\in H^1(\mathbb{T}^3)$ with $\displaystyle \int_{\mathbb{T}^3} \rho u_0(x)dx=0$
and if $u\in L^\infty([T_0,T_1]\times\Omega)$, then~:
\begin{equation}\label{POLYREYNOLDS}
\frac{1-C(n_-)}{\Re^{3\gamma/2}} \frac{\bar{E}}{E(T_0)} \leq
C \Re^2 \left\{ E^{-\delta_0 K_+} + \frac{\delta_0 K_+}{1+(\delta_0 K_+)^3}\:
\frac{3C_0^2 \tau}{T_1-T_0} \right\}
\end{equation}
provided~\eqref{NBTHASSUMPTION} holds and where
$\displaystyle \delta_0 =  \frac{C_0}{2}\sqrt{\nu \tau} =
\frac{C_0\nu^2}{2\sup_{[T_0,T_1]}\norme[L^2]{\omega}^2}$
is the analyticity radius guaranteed by Theorem~\ref{TH-HIGH}  on $[T_0+\tau,T_1]$.
\item[2. Well prepared data on $\mathbb{T}^3$ --]
If moreover one assumes that the initial data $u_0\in H^1(\mathbb{T}^3)$ was given at $t=0$
and that the K41-property holds on $[T_0,T_1]$ with
\begin{equation}\label{DELAY}
T_0\geq \frac{\nu^3}{\sup_{[0,T_1]}\norme[L^2]{\omega}^4}
\end{equation}
then Theorem~\ref{TH-HIGH} guarantees an analyticity radius of at least
$\displaystyle \delta_1 = \frac{C_0\nu^2}{2\sup_{[0,T_1]}\norme[L^2]{\omega}^2}$
on $[T_0,T_1]$. In this case, \eqref{POLYREYNOLDS} can be improved to read~:
\begin{equation}\label{EXPREYNOLDS}
\delta_1 K_+ \leq 
\log\left(\frac{C\Re^{3\gamma/2}}{1-C(n_-)} \right) + \log\left(\frac{E(0)}{\bar{E}}\right) + 2\log\Re.
\end{equation}
\item[3. Case of $\R^3$ --]
Similar results are also valid for $\Omega=\R^3$ except that the assumptions
$\int_{\mathbb{T}^3} \rho u_0(x)dx=0$ and~\eqref{NBTHASSUMPTION} 
must be dropped and replaced by $u_0\in L^1\cap L^2(\R^3)$. Estimate~\eqref{POLYREYNOLDS} reads~:
\begin{equation}
\frac{1-\Re^{-1/2}}{\Re^{3\gamma/2}}\:\frac{\bar{E}}{E(T_0)} \leq C \Re^{1/2}
\left\{ E^{-\delta_0 K_+} + \frac{\delta_0 K_+}{1+(\delta_0 K_+)^3}\:
\frac{3C_0^2 \tau}{T_1-T_0} \right\}
\end{equation}
and \eqref{EXPREYNOLDS} becomes instead~:
\begin{equation}
\delta_1 K_+ \leq 
\log\left(\frac{C\Re^{3\gamma/2}}{1-\Re^{-1/2}} \right) + \log\left(\frac{E(0)}{\bar{E}}\right) + 
\frac{1}{2}\log\Re.
\end{equation}
\end{description}
\end{thm}
\paragraph{Remarks}
\begin{enumerate}
\item
According to the comments at the end of \S\ref{PAR:TIMESCALE},
the ``well-preparedness'' assumption \eqref{DELAY} is likely
to be satisfied if the initial data is large enough.
\item
In the unprepared case, one almost gets the asymptotic
$\displaystyle \delta_0 K_+ \lesssim (E_0/\bar{E})\,\Re$
on $\mathbb{T}^3$ and 
$(E_0/\bar{E})^{1/2}\,\Re^{1/4}$ on $\R^3 $. In the well prepared case, it improves rigorously to~:
$$\delta_1 K_+ \leq C' + \log (E_0/\bar{E}) + \begin{cases}
2\log \Re &\text{on }\Omega=\mathbb{T}^3\\
\frac{1}{2}\log\Re  &\text{on }\Omega=\R^3.
\end{cases}$$
It will be shown in~\S\ref{PAR:23LAW} that there is an experimental
hint towards $\delta K_+\gtrsim 1$ and that $\delta K_+\simeq1$
must hold when physicists claim to observe fully developed turbulence.
\item
One can compute $\delta_i K_+$ ($i=0,1$) using only norms of the vorticity~:
$$\delta_i K_+ \simeq
\frac{\nu^{5/4}\bar\varepsilon^{1/4}}{\alpha^{3/4}\sup_{t}\norme[L^2]{\omega(t)}^2}\cdotp$$
Thus, using $\bar{E}^2\simeq\alpha^{3}\nu\bar{\varepsilon} \,\Re $
from~\eqref{EPSILONANDRE},
the well prepared case can be reformulated in the following form~:
\begin{equation}
\alpha^3 \nu \bar\varepsilon\cdot
\exp\left(\frac{c\nu^{5/4}\bar\varepsilon^{1/4}}{\alpha^{3/4}\sup_{[0,T_1]}\norme[L^2]{\omega}^2}\right)
\leq C \, E_0^2 \times \begin{cases}
\Re^{3}  &\text{on }\Omega=\mathbb{T}^3\\
1 \qquad &\text{on }\Omega=\R^3.
\end{cases}
\end{equation}
\end{enumerate}

\Proof
Let us first investigate the case of $\mathbb{T}^3$.
Comparison of Theorem~\ref{TH-HIGH} with \eqref{K41bis} for $K=K_+$ provides~:
$$\Re^{-3\gamma/4} \alpha \bar{\varepsilon}^{2/3} K_+^{-5/3} \leq 
\bar{E}^\ast(K_+)
\leq C  K_+ L^2 \langle e^{-\delta(t) K_+}\rangle \: E(T_0) $$
where $\delta(t)$ is defined in Theorem~\ref{TH-HIGH}. Let us write it down as~:
$$\Re^{-3\gamma/4} \alpha \bar{\varepsilon}^{2/3} K_+^{-2/3} 
\leq C  (K_+ L)^2 \langle e^{-\delta(t) K_+}\rangle \: E(T_0).$$
According to definition~\eqref{KcPERIODIC} and Theorem~\ref{RANGETHM}, one has~:
$$\alpha \bar{\varepsilon}^{2/3}K_+^{-2/3} = \left(\frac{K_c}{K_+}\right)^{2/3} \bar{E} \geq \frac{1-C(n_-)}{12\Re^{3\gamma/4}}\:\frac{\bar{E}}{\Re^{1/2}}\cdotp$$
Similarly, $(K_+ L)^2 = (K_- L)^2 \Re^{3/2} \leq (2\pi)^2 \Re^{3/2}$.
A direct computation of the time-average of $e^{-\delta(t)K_+}$ reads~:
$$\left\langle e^{-\delta(t)K_+} \right\rangle =
 e^{-\delta_0 K_+} +  \frac{\tau}{T_1-T_0} \Psi(\delta_0 K_+)$$
with $\displaystyle \tau = \frac{\nu^3}{\sup_{[T_0,T_1]}\norme[L^2]{\omega}^4}$ and
$\displaystyle \delta_0 = \frac{C_0}{2}\sqrt{\nu \tau}$ and the function
$$\Psi(s)=C_0^2 \left(\frac{2(1-(1+s)e^{-s})}{s}-e^{-s}\right)$$
that satisfy $\displaystyle \frac{1}{4}\left(\frac{C_0^2 s}{1+s^3}\right)\leq \psi(s)\leq 3\left(\frac{C_0^2 s}{1+s^3}\right)$ for any $s\in\R_+$.
This gives \eqref{POLYREYNOLDS}.
In the case of ``well prepared'' data, Theorem~\ref{TH-HIGH} gives $\delta(t)=\delta_0$
for any $t\in[T_0,T_1]$ thus $\left\langle e^{-\delta(t)K} \right\rangle = e^{-\delta_0 K_+}$
and $$\frac{1-C(n_-)}{\Re^{3\gamma/2}} \frac{\bar{E}}{E(0)} \leq
C \Re^2 e^{-\delta_0 K_+}$$
from which~\eqref{EXPREYNOLDS} follows immediately.

On $\R^3$ and provided $u_0\in L^1$, a comparison of
Theorem~\ref{TH-HIGH} and \eqref{K41bis} implies instead~:
$$\frac{1-\Re^{-1/2}}{\Re^{3\gamma/2}}\:\frac{\bar{E}}{E_0}
\leq C \left\langle e^{-\delta(t)K_+} \right\rangle \Re^{1/2}$$
which explains the different game of powers of the Reynolds number $\Re$.
\cqfd\endProof

\section{Final remarks and some open problems}
\label{INTERACT}

Let us end this article with some comments on the questions at stake and especially the
compatibility between the spectral theory presented above and the experimentally
accessible quantities called structure functions. I will conclude with a striking numerical
fact that suggests that turbulent flows might actually saturate the current
estimation of analytic regularity.

\subsection{Finding examples of turbulent flows}

The burning question is to find examples of K41-turbulent flows.
According to our definition, this problem can be split in two. The first step is
to find K41-functions that are solutions of \eqref{NS}. The second one is
to find, among those functions, some that satisfy the asymptotic~\eqref{K41}.

Proving that a given function is K41 requires to compute a lower bound on
$\operatorname{Vol}(u;[T_0,T_1])$ and then check that a substantial amount
of enstrophy is contained in the scales below this bound. It can at least be tested
numerically as this property is stable in Leray space.
Satisfying~\eqref{K41} is more subtle and is likely to require additional assumptions. The 
most obvious one will be a form of isotropy (with a proper definition,
still to be found) because the energy spectrum was defined with isotropic spectral cutoffs.
The question of an anisotropic theory is also widely open and
could be of interest in concrete situations. For example, atmospheric turbulence in the
jet stream (rapid air flows at high altitude, well known to jet pilots) is anisotropic
because the ratio height/width is a tiny parameter.

\medskip
For bounded domains, one could for example investigate
solutions whose vorticity satisfies~:
$$\left|\int_{\Omega} \frac{\omega(t,x)}{|\omega(t,x)|} \, dx\right| \leq c_0.$$
The cheapest conjecture is that for $c_0$ small enough (which means that the vortex
lines are somehow distributed isotropically),
the asymptotic \eqref{K41} should hold.

This is also where probabilities might prove handy.
The corresponding conjecture is that \eqref{K41}
holds in an average sense over a statistical ensemble of solutions, for some ergodic
probability measure.

However, these conjectures should be balanced by taking into account 
substantial fluctuations of the dissipation rate~$\varepsilon$ both in the temporal and 
in the spatial domain, \textit{i.e.}~intermittency (see \S\ref{PAR:INTERM}).

\subsection{Geometrical structures of turbulence, localization and intermittency}
\label{PAR:INTERM}

Understanding the geometrical structures involved in turbulence is a major challenge
and the core of modern research on turbulence (see \textit{e.g.}~\cite{brachet},
\cite{Constantin97}, \cite{Sreenivasan99} and the references therein).
The main idea is that substantial fluctuations of the dissipation rate exist both at large
scales (caused by the mechanisms of agitation) and at small scales (caused by
the stretching of vortex lines). In consequence, a refined theory of turbulence
cannot rely solely on the average value~$\bar{\varepsilon}$.
In this article the focus was on ``large scale'' turbulence, \textit{i.e.}~in the case where
the fluctuations of $\varepsilon$ are small compared to $\bar{\varepsilon}$.

\medskip
The properties of  $\operatorname{Vol}(u;[T_0,T_1])$ and $\mathcal{T}(u_0;\omega)$
in regard to the characteristic scales of the geometric structures in turbulent
flows call for a closer look.
For example, a starting point for studying intermittency is the following definition.
\begin{dfn}
Given $\Omega'(t)\subset\Omega$ a smooth family of smooth subsets of $\Omega$,
a function $u\in\mathcal{L}(\Omega)$ is said to be
a \textit{local K41-function} on $\Omega'\times[T_0,T_1]$ if $$u'(t,x)=u(t,x)\chi(t,x)$$
is a K41-function on $\Omega\times[T_0,T_1]$ where $\chi$ is a smooth cutoff
function such that
$$ \chi(t,x)=1  \enspace\text{if}\enspace x\in\Omega'(t) \qquad\text{and}\qquad
\operatorname{Vol}(\chi;[T_0,T_1])\simeq
\operatorname{Vol}(\mathbf{1}_{\Omega'};[T_0,T_1]).$$
\end{dfn}
The first question of localization and intermittency is to study
which sub-domains of $\Omega_0\times [T_0,T_1]$ are admissible
if $u$ is a local K41-function on $\Omega_0\times [T_0,T_1]$ and to determine
the corresponding parameters $K_\pm(\Omega')$, $\Re(\Omega')$
and $\gamma(\Omega')$. Using \eqref{K-Vol}, an obvious restriction
is $\operatorname{Vol}(u')\lesssim
\operatorname{Vol}(\chi)\lesssim \operatorname{Vol}(u)$.

\medskip
Related to the volume function is the following inequality, which
(with Leray's inequality) is one of the only estimates that
has a sub-linear growth in time~:
$$\norme[L^1(\R^3)]{u(t)} \leq \norme[L^1]{u_0} + C\sqrt{\frac{t}{\nu}}\norme[L^2]{u_0}^2.$$
Improving this inequality for a given flow on $\R^3$ or finding examples of saturations on
some time intervals would definitely be interesting.

\subsection{Modelisation versus PDEs}

When dealing with the spectral description of turbulence, two natural questions occur.
Which solutions of \eqref{NS} satisfy the spectral asymptotic \eqref{K41} ? 
And in that case, what are their qualitative properties ?

Historically, Kolmogorov addressed both questions simultaneously by proposing a model
based on strong probabilistic assumptions that cannot be checked directly, namely by
assuming that ``the difference $u(t,x)-u(t',x')$ is a probabilistic variable  
whose law does not depend on $t,t',x,x'$ and is not affected either by rotations
of the coordinate system'' \cite{K41a}.
Obviously the consequences of his predictions have been extensively studied and
most of them where found valid, but sometimes with not as good a precision
as one could hope for~: structure functions $S_p(\ell)$ (see \S\ref{PAR:STRUCT})
are predicted to scale as $\ell^{p/3}$ but for $p\geq4$ one observes that the
exponent $p/3$ must be corrected by some small negative term.
These mishaps  are known under the generic name of ``intermittency'' and
various models have been sought to explain them, including a second
probabilistic model by Kolmogorov~\cite{K41d}.

The main contribution of this article is to show that \textsl{one can study the
qualitative properties of turbulent flows independently from the
a-priori models of the structure of such a flow
and that it can be done using deterministic tools}. This path should
allow colleagues to concentrate on the core mathematical problems.

\subsection{Spectral problems for general domains}

One cannot deny that the tradition of probabilistic models in turbulence
is a convenient way to avoid dealing with
the spectrum of the Stokes operator on domains\ldots{}
However, the lack of precision on the true energy spectrum (the function denoted by
$E^\dagger$, defined on $\sigma(A^{1/2})$ where $A$ is the Stokes operator)
considerably darkens the foundations of some experimental protocols that focus
on the spectrum $E^\ast$ defined with Fourier coefficients and
assumed without proof to be the shell averages of~$E^\dagger$.

More precisely, as the spectral theory of the Stokes
operator is not known, the naive protocol is to collect data in some subregion,
then compute Fourier coefficients as if the flow was periodic,
using an FFT-type algorithm with ad-hoc anti-aliasing techniques,
and then finally compute the energy spectrum with the
formula~\eqref{EXPERIMENTAL} valid for $\mathbb{T}^3$.
The question is~: does one really look at the spectrum of $u$ ?

For example, sorting the spectrum of the Stokes operator
on $\prod(\R/L_i\mathbb{Z})$ and checking~\eqref{SPECTRALASSYMPTOTIC}
raises non-trivial problems of rational approximation as soon as $L_i\neq L_j$.
However, \eqref{SPECTRALASSYMPTOTIC} played a crucial role in
the proof of~\eqref{EXPERIMENTAL} on $\mathbb{T}^3$. Worse, if one assumes
generically simple eigenvalues, Weyl's asymptotic for the Dirichlet problem in
dimension $n$ gives $K_j^2\sim C(n;|\Omega|) j^{2/n}$ which suggests
$$K_{j+1}-K_j \lesssim \frac{ C'(n;|\Omega|) }{K_j^{n-1}}$$
with no general lower bound if the domain is a small perturbation of one with
multiple eigenvalues.

Of course, other less naive experimental protocols are used but they also have
their share of mathematical problems (see \S\ref{PAR:23LAW} and \S\ref{PAR:45LAW}
on structure functions).

\bigskip
The author conjectures that 
the properties of $E^\dagger(K,t):\sigma(A^{1/2})\times\R_+\to\R_+$
reflect a general arithmetic of the specific nonlinearity $\mathbb{P}((u\cdot\nabla) u)$
and thus could be truly universal.
However, discrepancies between the various possible spectra of the Stokes operator
might explain why $E^\ast$ defined by~\eqref{EXPERIMENTAL} is not always
the appropriate experimental quantity.

This might be a first key explaining why there seem to be ``many''
theories of turbulence, depending on what kind of flow one deals
with ; it is likely that the spectral theory for the Stokes
operator in an outer domain
has little in common with the spectral theory in a wind tunnel
or that of a swimming pool\ldots{}

\subsection{Universality of the Kolmogorov constant}
\label{UNIVERSALALPHA}

Conversely, to give some credit to the protocols based on computing the sum
of Fourier coefficients on spherical shells of frequencies, one could investigate
the properties of sub-domains of turbulent flows and show that ``away from
anything'' the rule of thumb is that of $\R^3$ or $\mathbb{T}^3$.

For example,
experimental evidence suggests that the Kolmogorov constant $\alpha$ 
in  \eqref{K41bis} is  universal and that, with reasonable precision, it does not
even depend on the flow or the shape of the domain  (provided turbulence
is homogeneous and isotropic). Published values
are in the range $\alpha\in [0.45,2.4]$ with a common agreement
\cite{KolmogorovConstant2} around $\alpha=0.5$.
See also \cite{KolmogorovConstant}, \cite{Heinz} and
the references therein. The proof of the universality of a small range for
$\alpha$ would be very instructive\footnote{As a ``per unit of mass'' theory,
the transformation $(u,p,\rho)\mapsto(u,\mu p,\mu\rho)$
was disregarded before. However, to the best of my ignorance,
the following experiment has never been done~:
compare the Kolmogorov constant $\alpha$ for flows of
fluids having similar viscosities but a fixed mass ratio between molecules,
when the the same number of molecules and similar velocity fields are involved.
For example, a water flow and a flow of heavy water where hydrogen atoms
are replaced with its heavier isotope deuterium. Would the Kolmogorov constant
be identical ?}.

\subsection{Scale-by-scale balance of energy and energy cascade}

The following identity describes the scale-by-scale energy budget and follows
directly from an energy estimate of \eqref{NS}.
Let us denote by $S_K=\chi\left(\frac{A^{1/2}}{K}\right)$ the low-pass filter and
$S_K^\ast$ its adjoint on $L^2(\Omega)$~:
\begin{align*}\frac{d}{dt}\left( \int_0^K E^\ast(k,t) dk \right) &=
\frac{1}{2}\frac{d}{dt} \left(\rho\norme[L^2]{S_K u}^2\right) \\&= 
-2\rho\nu\norme[L^2]{\nabla(S_K u)}^2-2\rho\left(S_K^\ast S_K((u\cdot\nabla) u)\vert u\right)_{L^2}\\
&\qquad+2\rho\left([S_K^\ast S_K,\nu \Delta]u-[S_K^\ast S_K,\nabla]p\vert u\right)_{L^2}
\end{align*}
On $\R^3$ or $\mathbb{T}^3$, all the commutators vanish and one gets~:
\begin{equation}\label{IDENTITYSK}
\frac{d}{dt}\left( \int_0^K E^\ast(k,t) dk \right) =- 2\rho\nu
\norme[L^2]{\nabla(S_K u)}^2 - 2\rho \left( S_K^\ast S_K(u\otimes u) \vert
 \nabla u \right)_{L^2}.
\end{equation}
This identity is the base of the famous ``energy cascade'' interpretation.
The left-hand side denotes the rate of change of the energy contained
at low frequencies,
\textit{i.e.}~``at scales larger than $K^{-1}$''. It is balanced by the energy dissipation
 at such scales (first term on the right-hand side) plus the so called ``energy flux to smaller 
 scales'' (second term). If this term is negative, energy is indeed transferred away
 from $[0,K]$ hence goes towards lower scales.
 There is no ``energy injection'' term as in our case the external
 force is zero.
\begin{dfn}
A solution of \eqref{NS} on $\Omega=\R^3$ or $\mathbb{T}^3$ has the ``energy cascade'' property at time $t\geq0$ if there exists $K_-^\ast<K_+^\ast$ such that~:
\begin{equation}
\forall K\in [K_-^\ast,K_+^\ast], \qquad
\sum_{i,j} \int_{\Omega} S_K^\ast S_K(u_i\otimes u_j) \cdot \partial_i u_j(t,x) dx>0.
\end{equation}
In the limit of non-smooth frequency cutoff, $S_K^\ast S_K\to \mathcal{F}^{-1}\circ
\mathbbm{1}_{|\xi|\leq K}\circ\mathcal{F}$.
\end{dfn}

There is yet no mathematical evidence that this property holds for a general
class of solutions of \eqref{NS}, not even for the few explicit solutions known.
A rigorous connection between this definition of the energy cascade and K41-turbulence is
an open problem, in particular the relation between $K_\pm$ and $K_\pm^\ast$.

\bigskip
Let us however mention this interesting property that was pointed out to
the author by~\cite{Paicu}.
\begin{thm}
Assume that $u$ is a Leray solution on $\R^3$ and that for all $t\in [T_0,T_1]$,
one has a ``reverse'' cascade~:
\begin{equation}\label{REVCASC}
\forall K\geq0,\qquad 
\sum_{i,j} \int_{\Omega} S_K^\ast S_K(u_i\otimes u_j) \cdot \partial_i u_j(t,x) dx\leq0.
\end{equation}
Then the energy equality $\displaystyle
\norme[L^2]{u(T_0)}^2=\norme[L^2]{u(T_1)}^2+2\nu\int_{T_0}^{T_1}
\norme[L^2]{\nabla u(t)}^2dt$ holds.
\end{thm}
\Proof
For all $K\geq0$, assumption \eqref{REVCASC} and the identity \eqref{IDENTITYSK}
provide~:
$$\norme[L^2]{S_K u(T_0)}^2 \leq \norme[L^2]{S_K u(T_0)}^2 - 2\int_{T_0}^{T_1}\left( S_K^\ast S_K(u\otimes u) \vert  \nabla u \right)_{L^2}
= \norme[L^2]{S_K u(T_1)}^2+2\nu \int_{T_0}^{T_1}\norme[L^2]{\nabla(S_K u)}^2.$$
Letting $K\to\infty$, one gets~:
$$\norme[L^2]{u(T_0)}^2\leq \norme[L^2]{u(T_1)}^2+2\nu\int_{T_0}^{T_1}
\norme[L^2]{\nabla u}^2 \leq \norme[L^2]{u(T_0)}^2$$
the right-hand side being the classical  Leray inequality.
\cqfd\endProof

\subsection{Structure functions $S_p(\ell)$}
\label{PAR:STRUCT}

As the energy spectrum is not always available,
other universal properties have been looked for as substitutes to~\eqref{K41}
and are commonly accepted as experimental evidence of turbulence.
However, the connection between these properties and the spectral definition of
K41-turbulence is not as rigorous as one could
hope for and constitutes an immediate source of interesting
mathematical problems.

\subsubsection{The second-order structure function}

\begin{dfn}
Let us introduce the correlation function~:
\begin{equation}\label{GAMMA}
\Gamma(t,y) = \rho \int_{\R^3} u(t,x+y)u(t,x) \,dx =E(t) - S_2(t,y)
\end{equation}
with
\begin{equation}
S_2(t,y)=\frac{\rho}{2}\int_{\R^3}|u(t,x+y)-u(t,x)|^2 dx.
\end{equation}
The following quantity is called the second-order structure function of $u$ on $[T_0,T_1]$~:
\begin{equation}
\forall \ell\geq0 \qquad
S_2(\ell) = \int_{\mathbb{S}^2} \bar{S}_2(\ell\theta) d\theta.
\end{equation}
Recall that the $\bar{S}_2(y)$ denotes the time average \eqref{TAVERAGE}
of $S_2(t,y)$ on $[T_0,T_1]$.
\end{dfn}

The following result is sometimes called the Wiener-Khinchin formula.
\begin{thm}\label{WIENERKINCHINE}
For any function $u\in L^1([T_0,T_1];L^2(\R^3))$, one has~:
\begin{equation}
\frac{S_2(\ell) }{4\pi}= \bar{E} -
\int_0^\infty  \frac{ \sin(\ell K)}{\ell K} \bar{E}^\dagger(K)dK
= \int_0^\infty \left(1-\frac{ \sin(\ell K)}{\ell K}\right) \bar{E}^\dagger(K)dK.
\end{equation}
\end{thm}
It is usually suggested in physics textbooks that ``the energy spectrum is the
Fourier transform of the correlation function'' but the formula is systematically left unstated.
As far as proof is concerned, it is usually claimed to be a consequence of various
probabilistic assumptions.
One can prove instead that it is a perfectly determinist statement that
relies on the following property~:
the Fourier transform in $\R^3$ commutes with the process of replacing a given function
by the radial one whose values are the averages of the initial function on each sphere.
\begin{prop}\label{FRADIAL}
For a function $f\in\mathcal{S}(\R^3)$, let us define~:
\begin{equation}
F(r) = \int_{\mathbb{S}^2} f(r\theta) d\theta \qquad\text{and}\qquad
G(\lambda) = \int_{\mathbb{S}^2} \hat{f}(\lambda\theta) d\theta
\end{equation}
and $S(\sigma) =  \sigma\sin \sigma$. Then, one has~:
\begin{subequations}
\begin{equation}\label{FOURIER}
\lambda^2 G(\lambda)=4\pi \int_0^\infty  F(r)\:S(r\lambda) dr.
\end{equation}
The  inversion formula reads~:
\begin{equation}\label{INVERSEFOURIER}
r^2 F(r) = \frac{1}{2\pi^2} \int_0^\infty G(\lambda) \: S(r\lambda) d\lambda.
\end{equation}
\end{subequations}
Moreover, if $f(x)=\frac{1}{4\pi}F(|x|)$ is radial, then $\hat{f}$ is also a radial function
thus $\hat{f}(\xi)=\frac{1}{4\pi}G(|\xi|)$ can be computed
with~\eqref{FOURIER}-\eqref{INVERSEFOURIER}.
\end{prop}
The radial statement is classical (see \cite[Chap. 4, Theorem 3.3]{Stein}).
The rest is implicit to the stability under Fourier transform of the
orthogonal decomposition of $L^2(\R^3)$
into spherical harmonics \cite[Chap. 4, Lemma 2.18]{Stein}.
Indeed, joining $L^2=\bigoplus \mathfrak{H}_j$ and
$\mathcal{F}:\mathfrak{H}_j\to \mathfrak{H}_j$ implies
that $\mathcal{F}$ commutes with the orthogonal projection on each~$\mathfrak{H}_j$. 
This abstract argument allows one to claim that $\frac{1}{4\pi}F$ is defined for
any $f\in L^2(\R^3)$ as its orthogonal projection on $\mathfrak{H}_0$
and that the statement still holds in this case.
\Proof Let us first establish the statement for a radial function $f(|x|)$ on $\R^3$.
Its Fourier transform is~:
$$|\xi|^2\hat{f}(\xi) = \int_0^\infty S\left(r|\xi| ;  \frac{\xi}{|\xi|}\right) f(r) dr
\quad\text{with}\quad
S(\rho,\theta')= \rho^2\int_{\mathbb{S}^2} e^{-i \rho \: \theta\cdot \theta'} d \theta.$$
The function $S:\R_+\times\mathbb{S}^2\to\mathbb{C}$ is invariant by rotations
of the second variable~:
$$S(\rho, \theta)=S(\rho,e_1).$$
As $\rho^{-1}S(\rho,e_1)$ has the same value and derivative at the
origin as $4\pi\sin\rho$, one has $S(\rho,e_1)=4\pi\rho\sin\rho$ (one could also
directly compute the integral in polar coordinates)
and $$|\xi|^2\hat{f}(\xi)=4\pi \int_0^\infty  f(r)\:S(r|\xi|) dr.$$
This proves the theorem in the case of a radial function.

\medskip
Let us now turn to the general case. One defines a function $f_0$ by~:
$$f(x)= \frac{1}{4\pi}F(|x|) + f_0(x)$$
which ensures that~:
$$\forall r\geq0, \qquad \int_{\mathbb{S}^2}f_0(r\theta) d\theta = 0.$$
Applying the theorem for the radial part gives~:
$$G(\lambda) = \frac{4\pi}{\lambda^2} \int_0^\infty F(r) \: S(r\lambda)dr + 
\int_{\mathbb{S}^2}\widehat{f_0}(\lambda\theta) d\theta .$$
The last term boils down easily using Fubini's theorem~:
\begin{align*}
\int_{\mathbb{S}^2}\widehat{f_0}(\lambda\theta) d\theta &= 
\frac{1}{\lambda^3} \iint_{\R^3\times\mathbb{S}^2}
e^{-iy\cdot\theta} f_0\left(\frac{y}{\lambda}\right) dyd\theta
\\&=\frac{1}{\lambda^3} \iiint_{\R_+\times\mathbb{S}^2\times\mathbb{S}^2}
e^{-i\rho\theta \cdot\theta'} f_0\left(\frac{\rho}{\lambda}\theta\right) \rho^2 d\rho d\theta d\theta'
\\&=\frac{4\pi}{\lambda^3} \iint_{\R_+\times\mathbb{S}^2}
\frac{\sin \rho}{\rho} \: f_0\left(\frac{\rho}{\lambda}\theta\right) \rho^2 d\rho d\theta
\\&=0.
\end{align*}
The second to last identity is the previous computation of $S(\rho,\theta)=4\pi\rho\sin\rho$
and the last one is the fact that the sphere averages of $f_0$ vanish.
\cqfd\endProof

\bigskip
\Proof[ of Theorem~\ref{WIENERKINCHINE}]
The first step is~:
$$S_2(\ell)=4\pi \bar{E}- \int_{\mathbb{S}^2}\bar{\Gamma}(\ell \theta) d\theta$$
which results immediately from the definitions.
Next, one observes that~:
\begin{equation}\label{PoorMansWIENERKINCHINE}
 \rho |\hat{u}(t,\xi)|^2 = \int_{\R^3} e^{-i y\cdot\xi} \Gamma(t,y) dy
\end{equation}
which express the duality between convolution and regular multiplication.
Taking the time average on~$[T_0,T_1]$, then sphere averages, and finally
substituting the definition~\eqref{EdaggerR3}, one gets~:
$$\int_{\mathbb{S}^2}
\widehat{\bar{\Gamma}}(K\theta) d\theta = \frac{(2\pi)^3}{K^2} \bar{E}^\dagger(K).$$
The conclusion now follows directly from Proposition~\ref{FRADIAL}.
\cqfd\endProof

\subsubsection{Range of validity of the ``2/3 law'' and radius of analyticity of $u$}
\label{PAR:23LAW}

Experimental evidence \cite[p.57-61]{Frish} suggests that~:
\begin{equation}\label{S2}
\forall \ell \in [\eta,\ell_0], \qquad
S_2(\ell) \quad \simeq \quad \beta ( \bar{\varepsilon} \ell)^{2/3}
\end{equation}
where $\eta\simeq K_+^{-1}$ is the dissipation scale and
$\ell_0\simeq K_-^{-1}$ is the size of large eddies.
Usually, physics courses state that this so called 2/3 law is equivalent
to the $K^{-5/3}$ decay of the energy spectrum. This ``equivalence'' calls for a closer look.

 \begin{figure}[ht!]\sf
 \begin{center}
\includegraphics[width=.89\textwidth]{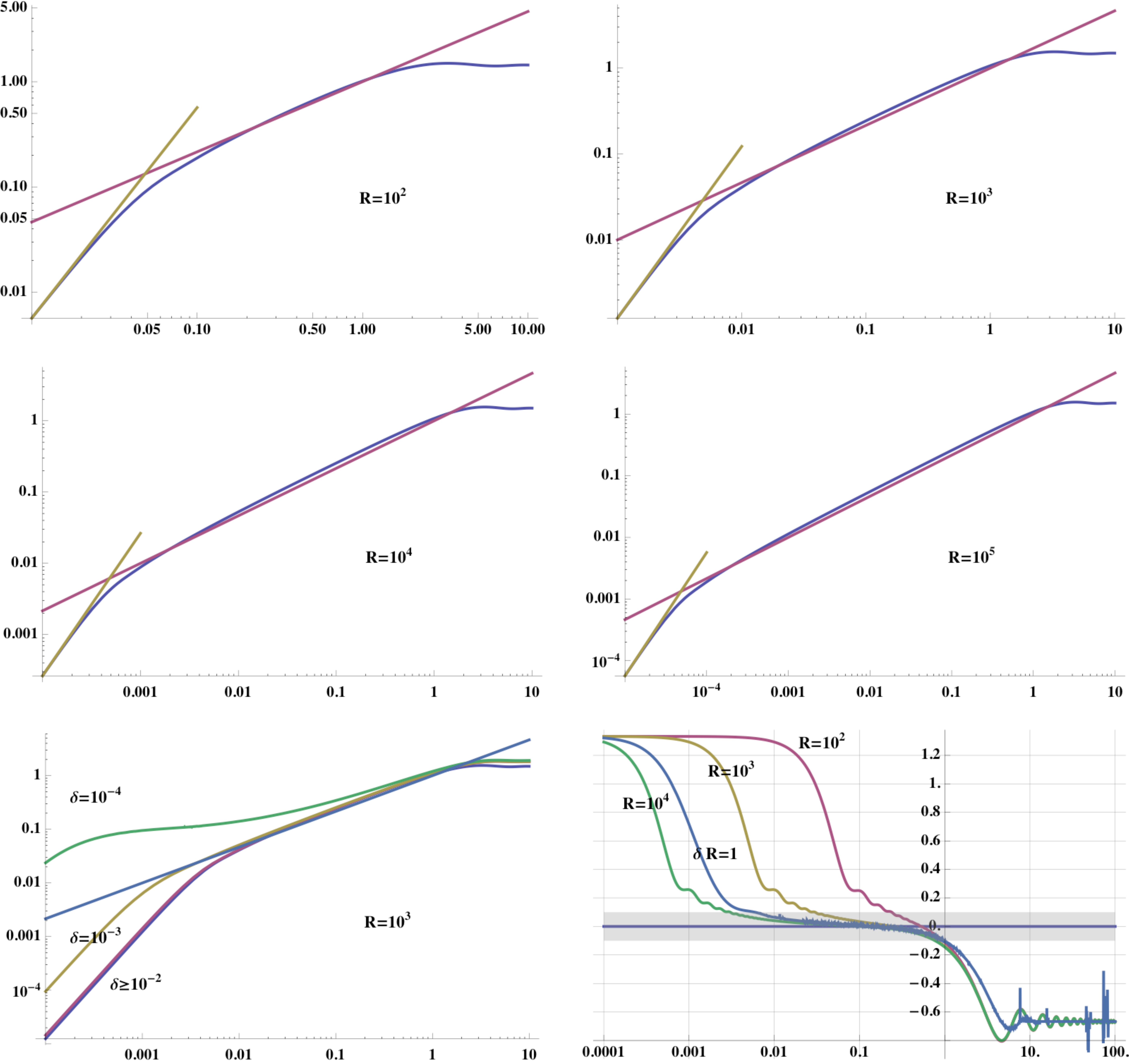}
\caption{\label{FIGWK} Numerical investigation of the domain of validity of~\eqref{S2}.}
\end{center}
\textsl{Top four} -- $S_2(\ell)$ is computed by Theorem~\ref{WIENERKINCHINE}
using an idealized spectrum $\bar{E}^\dagger(K)=K^{-5/3} \mathbf{1}_{[1,R]}(K)$ for
different values of $R$. Plot in Log-Log scale on $[R^{-1},10]$. Also represented is $(R\ell)^2 S_2(R^{-1})$ for $\ell<10/R$ and $\ell^{2/3}$. One observes a close fit of $S_2(\ell)$
and $\ell^{2/3}$ on $[10/R,1]$. The graph was obtained by formal integration with
\textit{Mathematica}\copyright.\\
\textsl{Bottom left} --  $S_2(\ell)$ computed for a ``real'' energy spectrum
$$\bar{E}^\dagger(K)=\begin{cases}
K^2 & K\leq 1\\
K^{-5/3} & K\in [1,R]\\
R^{-5/3}e^{-\delta(K-R)} & K\geq R
\end{cases}$$
with $R=10^3$ and a radius of analyticity $\delta\in\{10^{-2},10^{-3},10^{-4}\}$
(this time with numerical integration).
The range of validity of~\eqref{S2} is maximal for~$R\delta=1$  but
drops drastically when $R\delta\ll1$.
\\
\textsl{Bottom right} -- Spectral precision $\ell \frac{d}{d\ell}\log S_2(\ell)-\frac{2}{3}$
for an idealized spectrum and $R\in\{10^{2},10^{3},10^{4}\}$. When the precision is
in the gray band, \eqref{S2} is satisfied with a relative error of less than 10\%.
Also represented is the precision for the previous
``real'' spectrum with $R=10^{3}$ and $\delta=10^{-3}$. For $\ell\geq10$, the
function $S_2(\ell)$ is oscillatory, which reduces the precision of the
numerical integration in the case of the ``real'' spectrum.
\end{figure}

\medskip
The poor man's argument is the following. Let us consider
an ideal case where $u$ would be a radial power function whose
energy spectrum is exactly $K^{-5/3}$ (this means $|\hat{u}(t,\xi)|=|\xi|^{-11/6}$
even though this function is not a solution of~\eqref{NS} nor even a square integrable one).
Then, for any $\ell>0$ and $\theta\in\mathbb{S}^2$, \eqref{PoorMansWIENERKINCHINE} gives~:
$$\bar{\Gamma}(\ell \theta)=\ell^{2/3}\bar{\Gamma}(\theta).$$
However in this case $\bar{E}=+\infty$ and $S_2(\ell)$ are undefined. 
Even excellent physics textbooks do not further justify the ``equivalence'' between
the 2/3 law and the $K^{-5/3}$ except by a vague reference to a probabilistic version of
Theorem~\ref{WIENERKINCHINE}\ldots

Slightly more careful physicists \cite[p.~87]{Frish} state
that the 2/3 law is only an asymptotic property that holds provided the limits are taken in the
proper order~: $T_1-T_0\to\infty$, then $\nu\to0$, then $\ell\to0$, and
that taking the limits in any other order will lead to trouble.
From a mathematical point of view,
a rigorous upper bound for $S_2(\ell)$ can be found in \cite{CNT99}.
However, the relation with \eqref{K41} is not established.

\bigskip
One can check numerically that \eqref{K41}
and the 2/3 law indeed share a strong connection  (see Figure~\ref{FIGWK}).
The conclusion of this computation is that the domain of validity of \eqref{S2}
is $[C K_+^{-1},K_-^{-1}]$ if the analyticity radius of $u$ exceeds $K_+^{-1}$.
However, if the analyticity radius becomes smaller than $K_+^{-1}$, then the
range of validity of~\eqref{S2} shrinks dramatically.

This fact should be put in perspective with a common experimental fact~:
physicists observe the range of validity of~\eqref{S2} to be often of much
smaller amplitude than the inertial range.  Figure~\ref{FIGWK} suggests
that  the analyticity radius $\delta$ of such a flow is \textsl{smaller} than Kolmogorov's
dissipation scale $K_+^{-1}$. Conversely, a fully developed turbulence for
the structure function $S_2(\ell)$ indicates that $\delta\geq K_+^{-1}$.

\bigskip
Even though this numerical study is encouraging, it leaves the rigorous connexion
between K41-turbulence and the structure function $S_2(\ell)$ on the list of open
problems. Moreover, as the numerical observation suggests
that $S_2(\ell)=O(\ell^2)$ around $\ell\to0$ thus
one cannot expect to prove \eqref{S2} by a finite expansion of $S_2(\ell)$.
Instead one will have to prove directly that there exists a smooth function $\gamma_0$ with
$\gamma_0(s)\ll1$ if $s\geq1$ and such that
$$\sup_{\ell\in[C K_+^{-1}; K_-^{-1}]}\left|\ell \frac{d}{d\ell}\log S_2(\ell)-\frac{2}{3}\right|
\log\left(\frac{K_-^{-1}}{C K_+^{-1}}\right)\leq \gamma_0(\delta K_+)$$
where $\delta$ denotes the radius of analyticity of a K41-function $u$.
This subtle exercise in harmonic analysis should be an excellent warm-up
round before tackling the question of finding examples of turbulent flows.
This also explains why physicists dodge the problem of accessing the ``intermediary''
regime of~\eqref{S2} by taking suitable asymptotics that will push its domain
of validity all the way to $\ell\to0$.

\subsubsection{Higher order structure functions and the ``4/5 law''}
\label{PAR:45LAW}

Other structure functions play a central role in experimental protocols.
The most celebrated is the so called ``4/5 law at inertial scales'' for the third-order function
(this function was historically related to the measurement mechanisms used
to acquire experimental data in real flows)~:
\begin{equation}\label{S3}
S_3^{||}(\ell) =\left\langle
\iint_{\Omega\times\mathbb{S}^2} \left((u(t,x+\ell\theta)-u(t,x))
\cdot \theta \right)^3 \frac{dx d \theta}{\operatorname{Vol}(\Omega)}
\right\rangle
\quad \simeq \quad -\frac{4}{5} \bar{\varepsilon} \ell
\end{equation}
Recall that the brackets $\langle\cdot\rangle$ denote time average.
This fact is claimed as being ``rigorously established'' in most physics textbooks
(see \cite[chap. 2]{brachet} or \cite[p.76-86]{Frish}). It was indeed addressed in the first
paper of Kolmogorov \cite{K41a}. However, the rigorous path from~\eqref{K41}
to~\eqref{S3} still requires some enlightenment. Worse, one can conjecture that
the property $S_3(\ell)=C\ell+O(\ell^2)$ holds independently from the asymptotic~\eqref{K41}
and that it is only a consequence of the rapid decay at infinity of the energy spectrum.
Let us follow a proof step by step and point out the dark spots.

\medskip
The starting point is the following identity, called the Karman-Howarth-Monin relation (let us
recall that $\Gamma$  has been defined by~\eqref{GAMMA})~:
\begin{equation}
\partial_t \Gamma - 2\nu \Delta_y \Gamma = \sigma
\quad\text{with}\quad
\sigma(t,y)=\frac{\rho}{2}\int_{\R^3}
\operatorname{Tr} \left\{^{t}(u\otimes u)(t,x) \cdot \left[\nabla u(t,x+y)+\nabla u(t,x-y)\right]
\right\} \,dx.
\end{equation}
Even though it is often presented as a consequence of the probabilistic assumptions,
it is a perfectly deterministic relation that
follows immediately from equation \eqref{NS}, the definition of~$\Gamma$  and
the obvious fact that~:
$$\Delta \Gamma(t,y) = -\frac{\rho}{2}\int_{\R^3} \nabla u (x+y)\nabla u(x)dx.$$
One can rewrite the right-hand side in the following way~:
\begin{equation}
\sigma(t,y) = \frac{\rho}{4}\div_y
\left(\int_{\R^3} |\delta(t,x,y)|^2 \delta(t,x,y) \,dx\right).
\end{equation}
with $\delta(t,x,y)=u(t,x+y)-u(t,x)$.

\medskip
Let us now compute the Fourier transform and integrate over a sphere of radius $K$.
One gets~:
\begin{equation}\label{EVOLSPECTRUM}
\partial_t E^\dagger(K,t) + 2\nu K^2 E^\dagger(K,t) = \partial_K \Pi(K,t).
\end{equation}
One can compute $\Pi(K,t)$ directly from $\sigma(t,y)$ using
Proposition~\ref{FRADIAL}~:
\begin{align}
\notag
\Pi(K,t) &=\frac{2}{(2\pi)^3}\int_0^K \int_{\mathbb{S}^2}\int_{\R^3} 
e^{-ki y\cdot\theta}\sigma(t,y) k^2 dk d\theta dy\\
\notag
&=\frac{1}{\pi^2}\int_{\R^3} \frac{\sin(K|y|)-K|y|\cos(K|y|)}{|y|^3} \sigma(t,y)dy\\
\notag
&=\frac{1}{2\pi}\int_0^\infty \frac{2}{\pi}\frac{\sin (K\ell)}{\ell}\times
 (1+\ell\partial_\ell)\left[\sigma_0(t,\ell)\right] d\ell
\end{align}
where 
$\displaystyle \sigma_0(t,\ell) = \int_{\mathbb{S}^2} \sigma(t,\ell\theta) d\theta$
and $\displaystyle \int_0^\infty \frac{2}{\pi} \frac{\sin (K\ell)}{\ell}  d\ell =1$.
Next, one computes the time averages on $[T_0,T_1]$.
The Fourier inversion formula in Proposition~\ref{FRADIAL} then reads~:
\begin{equation}\label{SIGMA0}
 (1+\ell\partial_\ell)\bar{\sigma}_0=
 \frac{\ell}{8\pi} \int_0^\infty \frac{\sin K\ell}{K} \, \bar{\Pi}(K)dK.
 \end{equation}
The time average of \eqref{EVOLSPECTRUM} over $[T_0,T_1]$
then integrated over $[K_-,K]$ reads~:
$$\bar{\Pi}(K)-\bar{\Pi}(K_-)=
2\nu \int_{K_-}^K \kappa^2 \bar{E}^\dagger(\kappa) d\kappa -
\int_{K_-}^{K}\frac{E^\dagger(\kappa,T_0)-E^\dagger(\kappa,T_1)}{T_1-T_0}d\kappa.$$
If the spectrum decays sufficiently fast,
the right hand-side is almost constant in the range $K\geq K_+$.
Each term is also equivalent to $\bar{\varepsilon}$ if the solution is smooth.
Thus, for $K\geq K_+$, one gets $\bar{\Pi}(K)\simeq \bar{\Pi}(K_-)$. 
Usually, this claim is the first dark spot justified by some ``suitable'' asymptotic
like letting $t\to\infty$, then $\nu\to0$ with fixed $\bar{\varepsilon}$. It seems however
to be merely a consequence of the analytic regularity.

The last step is to use  \eqref{SIGMA0} to
convert the constancy of $\Pi_\infty=\bar{\Pi}(K)=\bar\Pi(K_-)$ for $K\geq K_+$.
The right-hand side can be developed as follows for $\ell\to0$~:
$$ (1+\ell\partial_\ell)\bar{\sigma}_0=\frac{\Pi_\infty}{8}+O(\ell).$$
This singular ODE admits only one bounded
solution as $\ell\to0$, namely~$\bar{\sigma}_0(\ell) = \frac{1}{8}\Pi_\infty+O(\ell)$.
Substituting the definition of $\sigma(t,x)$, one gets~:
$$2\bar{\sigma}_0(\ell) =
\frac{1}{\ell^2}\frac{d}{d\ell}\left[ \ell^2
\overline{\iint_{\R^3\times\mathbb{S}^2} |\delta(t,x,\ell\theta)|^2 \delta(t,x,\ell\theta) \,dxd\theta}
\:\right]
= \frac{\Pi_\infty}{4} + O(\ell).$$
It is a weak form of~\eqref{S3} that implies by integration of the finite expansion~:\begin{equation}
\overline{\iint_{\R^3\times\mathbb{S}^2} |\delta(t,x,\ell\theta)|^2 \delta(t,x,\ell\theta) \,dxd\theta}=C\ell +O(\ell^2)
\end{equation}
for some numerical constant $C=\frac{1}{12}\Pi_\infty$.
As the remainder terms have been neglected, this identity is not rigorously established but
one can conjecture that it will hold for any smooth
solution of~\eqref{NS} whose spectrum decays sufficiently fast.

\medskip
Using a probabilistic assumption of homogeneity and isotropy in the region
of observation~$\Omega$, physicists claim that~:
$$\bar{\sigma}_0(\ell) = -\frac{\operatorname{Vol}(\Omega)}{96}
(3+\ell\partial_\ell)(5+\ell\partial_\ell)
\left[\frac{S_3^{||}(\ell)}{\ell}\right].$$
This identity is the second dark spot because it is not clear how to get a similar
formula independently of any a-priori model of the flow.
Substitution in the previous equation for $\sigma_0$ then gives~:
$$-\frac{1}{12}(1+\ell\partial_\ell)(3+\ell\partial_\ell)(5+\ell\partial_\ell)\left[\frac{S_3^{||}(\ell)}{\ell}\right]=\Pi_\infty$$
and this singular ODE admits only one bounded solution as $\ell\to0$, namely
$S_3^{||}(\ell)=-\frac{4}{5}\Pi_\infty\ell$. One can check immediately that
$\Pi_\infty$ has the dimensions of a dissipation and physicists claim indeed (and this
is the third dark spot) that $\Pi_\infty=\bar\varepsilon$.

\subsection{Extremal properties of turbulent flows regarding the optimality of analytic estimates.}
\label{parTROUBLING}

Let us put an end to this section and this article with a striking observation inspired by a
comment of Claude Bardos~:
``Personally I do not believe that the solutions of the incompressible Euler or Navier-Stokes equations blow up, but it may well be that there are no other general estimates than the one presently found'' \cite{Bardos01}.

\begin{figure}[ht!]\sf
\begin{minipage}{\textwidth}
\begin{multicols}{2}
\begin{center}
\includegraphics[width=.5\textwidth]{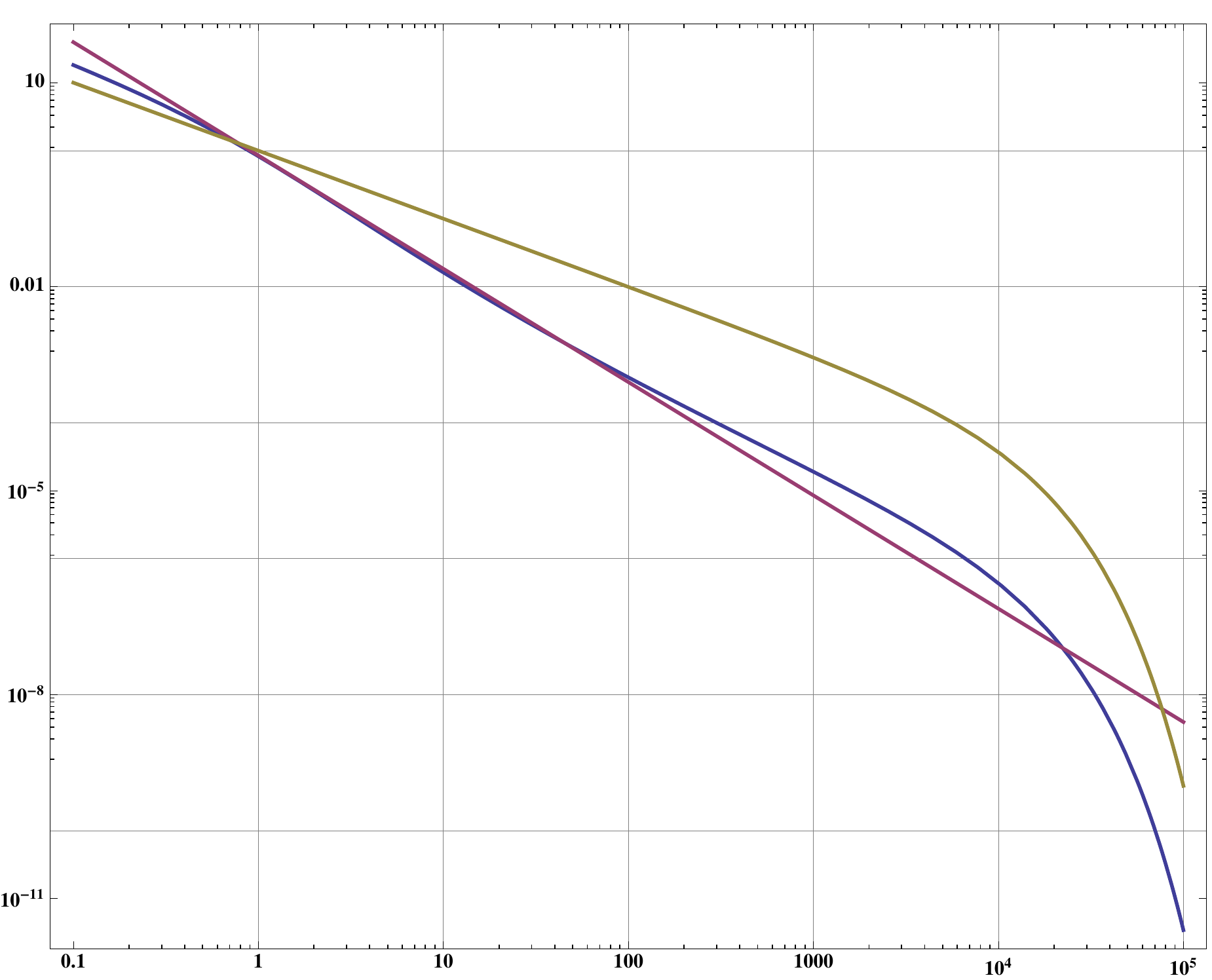}
\end{center}
\columnbreak
\phantom{.}
\includegraphics[width=.5\textwidth]{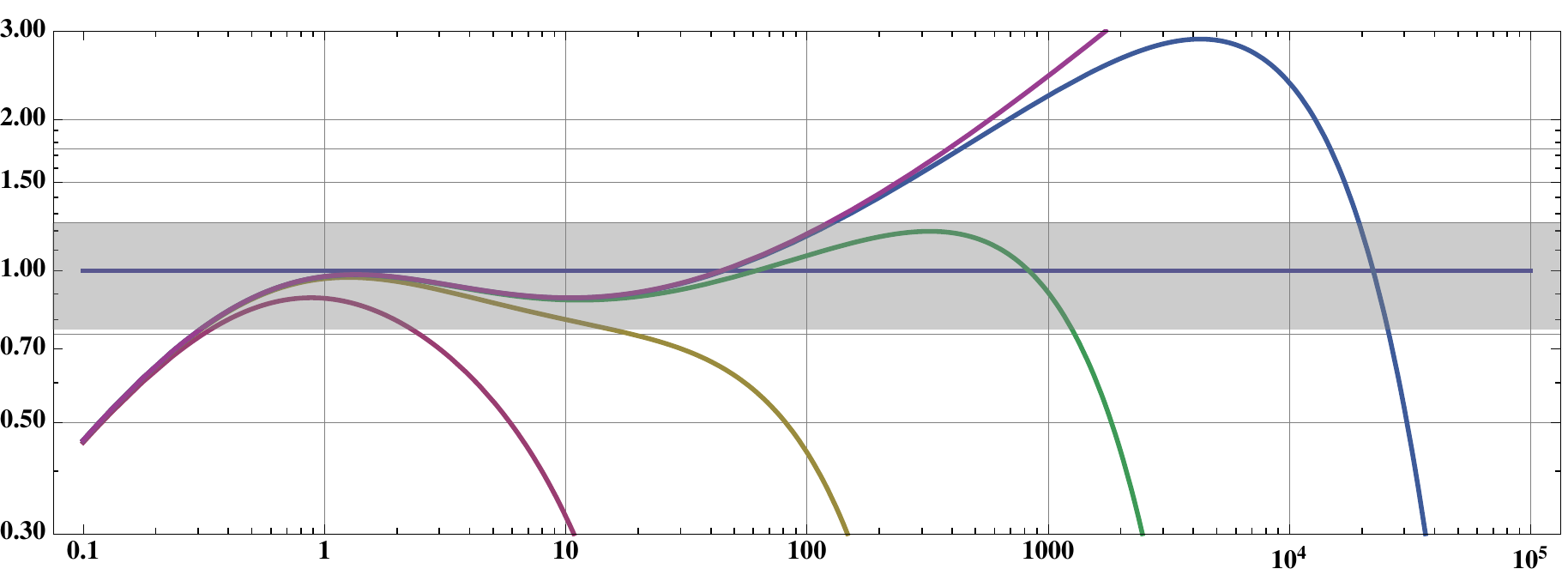}
\textsl{Left} --
Plot of $\chi_\delta(K)$ in Log-Log scale for $\delta=10^{-4}$. The straigh
line is $0.85 K^{-5/3}$. The upper graph is the rigorous bound $K^{-1} e^{-\delta K}$
of the energy spectrum.\\
\textsl{Top} --
Relative error $\chi_\delta(K)/(0.85 K^{-5/3})$ in Log-Log scale for $\delta\in\{10^{-1}, 10^{-2},
10^{-3}, 10^{-4}, 10^{-5}\}$. For $\delta=10^{-3}$, it shows that $\chi_\delta(K)$ matches~$0.85K^{-5/3}$ with less than 25\% of relative error on almost 4 decades.
\end{multicols}
\end{minipage}
\caption{\label{STRIKING} A striking observation : a logarithmic correction to the 
rigorous analytic estimates of the energy spectrum has a definite $K^{-5/3}$ behavior on
a large range of $K$.}
\end{figure}

\bigskip
Let us recall our upper bound of the energy spectrum for smooth ``old'' solutions~:
$$\bar{E}^\ast(K) \lesssim K^{-1}e^{-\delta K} E(T).$$
This estimate has been proved by a technique that seems to be 
the cutting edge of quantitative smoothness estimates
for parabolic equations. It is natural to ask whether flows
exist for which this inequality is  optimal on some
large range of $K$. If this is the case, then the energy
spectrum of such flows would exceed any substantial correction
to this estimate.

Let us plot therefore the following logarithmic correction (see Figure~\ref{STRIKING})~:
$$\chi_\delta(K)=\frac{K^{-1}e^{-\delta K}}{\log^2(2+K)}\cdotp$$
An extremely troubling fact is that, on a log-log diagram, this corrector
shows a definite $K^{-5/3}$ behavior !
For example when $\delta\in[10^{-5},10^{-1}]$, the
$K^{-5/3}$ behavior appears for roughly $K\in[1,\frac{1}{\delta}]$.

This observation is a powerful suggestion that K41-turbulent flows might exist
among smooth solutions of \eqref{NS} and that these flows are responsible
for the failure of extending local regularity methods. They will nonetheless
provide examples saturating the classical inequalities of
fluid mechanics.

\end{document}